\documentclass[a4paper, 11pt]{article}

\newcommand{\anote}[1]{}
\renewcommand{\anote}[1]{\textcolor{brown}{\textbf{[Aras: #1]}}}

\usepackage[english]{babel}
\usepackage[T1]{fontenc}
\usepackage[utf8]{inputenc}
\usepackage{amssymb}
\usepackage{amsthm}
\usepackage{array,amsmath, mleftright, thmtools}
\usepackage{setspace}
\usepackage{tcolorbox}
\usepackage[dvipsnames]{xcolor}
\definecolor{linkgreen}{HTML}{009A33}
\usepackage{blindtext}
\usepackage{mathrsfs}

\usepackage{mathtools}
\usepackage[breaklinks=true]{hyperref}
\hypersetup{
    colorlinks,
    linkcolor={linkgreen},
    citecolor={blue!50!black},
    urlcolor={blue!50!black}
}
\usepackage{bbding}
\usepackage{framed}
\usepackage{xcolor,colortbl}
\usepackage{hhline}
\usepackage[normalem]{ulem}

\renewcommand{\sout}{\bgroup\markoverwith{\textcolor{red}{\rule[0.5ex]{2pt}{0.4pt}}}\ULon}

\usepackage{dsfont}
\usepackage{cases}
\usepackage{pdfpages}
\usepackage{empheq}

\DeclareMathOperator{\dom}{dom}

\newcommand{\bQ}{\mathbb{Q}}
\newcommand{\bR}{\mathbb{R}}
\newcommand{\bS}{\mathbb{S}}

\newcommand{\cP}{\mathcal{P}}

\newcommand{\cZ}{\mathcal{Z}}
\let\hat\widehat

\colorlet{revpurple}{Purple}

\newcommand{\psout}{\bgroup\markoverwith{\textcolor{revpurple}{\rule[0.5ex]{2pt}{0.4pt}}}\ULon}

\newcommand\argmax{\operatorname{argmax}}
\newcommand\argmin{\operatorname{argmin}}

\usepackage{bm}
\usepackage[inkscapearea=page]{svg}

\usepackage{caption}[font=small,skip=0pt]

\newtheorem{lemma}{Lemma}
\newtheorem{assumption}{Assumption}
\newtheorem{proposition}{Proposition}
\newtheorem{corollary}{Corollary}

\newtheorem{example}{Example}

\newtheorem{definition}{Definition}

\newtheorem{remark}{Remark}

\allowdisplaybreaks

\makeatletter
\newcommand{\support}[1]{\aslv@support#1\@nil}
\def\aslv@support#1\\#2|#3\@nil{%
  \mleft(
  \left( \genfrac..{0pt}{}{#1}{#2}
  \mskip\aslv@supportskip \right)
  \middle|
  \mskip\aslv@supportskip
  #3
  \mright)
}
\newmuskip\aslv@supportskip
\makeatother
\makeatletter
\newcommand{\pushright}[1]{\ifmeasuring@#1\else\omit\hfill$\displaystyle#1$\fi\ignorespaces}
\makeatother

\makeatletter
\newcommand{\prr}[1]{\ifmeasuring@\else\omit\hfill\displaystyle#1\fi\ignorespaces}
\makeatother

\usepackage{graphicx}
\usepackage{tikz}
\usetikzlibrary{arrows.meta}
\usetikzlibrary{bending}
\usetikzlibrary{decorations.pathreplacing}
\usetikzlibrary{shapes}
\usetikzlibrary{patterns}
\usetikzlibrary{decorations.pathmorphing} 
\usetikzlibrary{positioning}
\usepackage{pgfplots}
\pgfplotsset{compat=1.15}
\usepackage{adjustbox}
\usepackage{subcaption}
\usepackage{wrapfig}
\usepackage{makecell}

\usepackage{natbib}
 \bibpunct[, ]{(}{)}{,}{a}{}{,}%
\usepackage{multirow}
\usepackage{booktabs}
\usepackage{longtable}

\usepackage{tabu}
\usepackage{pdflscape}
\makeatletter
\newbox\LT@firstfoot
\def\endfirstfoot{\LT@end@hd@ft\LT@firstfoot}
\newdimen\LT@footdiff
\def\LT@start{%
  \let\LT@start\endgraf
  \endgraf\penalty\z@
  \vskip\LTpre
  \LT@footdiff-\ht\LT@foot
  \advance\LT@footdiff\ht\LT@firstfoot
  \dimen@\pagetotal
  \advance\dimen@ \ht\ifvoid\LT@firsthead\LT@head\else\LT@firsthead\fi
  \advance\dimen@ \dp\ifvoid\LT@firsthead\LT@head\else\LT@firsthead\fi
  \advance\dimen@ \ht\ifvoid\LT@firstfoot\LT@foot\else\LT@firstfoot\fi
  \dimen@ii\vfuzz
  \vfuzz\maxdimen
  \setbox\tw@\copy\z@
  \setbox\tw@\vsplit\tw@ to \ht\@arstrutbox
  \setbox\tw@\vbox{\unvbox\tw@}%
  \vfuzz\dimen@ii
  \advance\dimen@ \ht
      \ifdim\ht\@arstrutbox>\ht\tw@\@arstrutbox\else\tw@\fi
  \advance\dimen@\dp
      \ifdim\dp\@arstrutbox>\dp\tw@\@arstrutbox\else\tw@\fi
  \advance\dimen@ -\pagegoal
  \ifdim \dimen@>\z@\vfil\break\fi
  \global\@colroom\@colht
  \ifvoid\LT@firstfoot
    \ifvoid\LT@foot
    \else
      \advance\vsize-\ht\LT@foot
      \global\advance\@colroom-\ht\LT@foot
      \dimen@\pagegoal\advance\dimen@-\ht\LT@foot\pagegoal\dimen@
      \maxdepth\z@
    \fi
  \else
    \advance\vsize-\ht\LT@firstfoot
    \global\advance\@colroom-\ht\LT@firstfoot
    \dimen@\pagegoal\advance\dimen@-\ht\LT@firstfoot\pagegoal\dimen@
    \maxdepth\z@
  \fi
  \ifvoid\LT@firsthead\copy\LT@head\else\box\LT@firsthead\fi\nobreak
  \output{\LT@output}%
}
\def\LT@output{%
  \ifnum\outputpenalty <-\@Mi
    \ifnum\outputpenalty > -\LT@end@pen
      \LT@err{floats and marginpars not allowed in a longtable}\@ehc
    \else
      \setbox\z@\vbox{\unvbox\@cclv}%
      \ifdim \ht\LT@lastfoot>\ht\LT@foot
        \dimen@\pagegoal
        \advance\dimen@-\ht\LT@lastfoot
        \ifdim\dimen@<\ht\z@
          \setbox\@cclv\vbox{\unvbox\z@\copy\LT@foot\vss}%
          \@makecol
          \@outputpage
          \setbox\z@\vbox{\box\LT@head}%
        \fi
      \fi  
      \global\@colroom\@colht
      \global\vsize\@colht   
      \vbox
        {\unvbox\z@\box\ifvoid\LT@lastfoot\LT@foot\else\LT@lastfoot\fi}%
    \fi
  \else
    \ifvoid\LT@firstfoot
      \setbox\@cclv\vbox{\unvbox\@cclv\copy\LT@foot\vss}%
      \@makecol
      \@outputpage
      \global\vsize\@colroom
    \else
      \setbox\@cclv\vbox{\unvbox\@cclv\box\LT@firstfoot\vss}%
      \@makecol
      \@outputpage
      \global\advance\@colroom\LT@footdiff
      \global\vsize\@colroom
    \fi
    \copy\LT@head\nobreak
  \fi
}
\makeatother

\usepackage{lipsum}

\usepackage[shortlabels]{enumitem}

\usepackage{todonotes}

\usepackage[ruled, vlined]{algorithm2e}
\SetKwRepeat{Do}{do}{while}

\usepackage{algpseudocode}
\usepackage{float}

\usepackage[top=25mm, bottom=25mm, left=25mm, right=25mm]{geometry}

\newcommand \Template[2]{
\usepackage{fancyhdr}
\pagestyle{fancy}
\fancyhf{}
\fancyhead[RO]{#1}
\fancyhead[LO]{#2}
\fancyfoot[CO]{\thepage}

\renewcommand{\headrulewidth}{1pt}
\renewcommand{\footrulewidth}{1pt}}

\usepackage{titling}
\usepackage{authblk}
\newcommand \linedabstract[1]{
  \renewcommand{\maketitlehookd}{%
    \mbox{}\medskip\par
    \centering
    \hrule\medskip
    \begin{minipage}{0.9\textwidth} 
    {\small\textbf{Abstract}\\ #1}
    \end{minipage}\medskip\hrule
    }      
}

\definecolor{darkred}{rgb}{0.8,0,0}



\usepackage{soul}

\definecolor{airforceblue}{rgb}{0.36, 0.54, 0.66}
\definecolor{darkbrown}{rgb}{0.4, 0.26, 0.13}

    \title{\textbf{Efficient Data Fusion in \\ Distributionally Robust Optimization}}
\author[$\dagger$]{Buse \c{S}en}
\author[$\ddagger$]{Aras Selvi}
\author[$\dagger\dagger$]{Jean Pauphilet}
\author[$\dagger$]{Daniel Kuhn}
\affil[$\dagger$]{Risk Analytics and Optimization Chair, EPFL, Lausanne, Switzerland \protect\\
\texttt{buse.sen@epfl.ch}, \texttt{daniel.kuhn@epfl.ch}}
\affil[$\ddagger$]{Operations and Technology, UCL School of Management, London, UK \protect\\
\texttt{a.selvi@ucl.ac.uk}}
\affil[$\dagger\dagger$]{Management Science and Operations, London Business School, London, UK \protect\\
\texttt{jpauphilet@london.edu}}

\linedabstract{Decisions under uncertainty increasingly rest on fragmented evidence, such as samples from related but different populations, privacy-protected summary statistics, or complementary statistical descriptions of the same data. Distributionally robust optimization (DRO) can fuse such evidence by constructing one ambiguity set per source and hedging against all distributions in the intersection of these sets. This fusion approach is less conservative than aggregating the sources into a single ambiguity set, but it is computationally demanding. Exact reformulations grow exponentially with the number of sources, and even deciding whether the intersection is nonempty is NP-hard. We show that fusion DRO admits an exact reformulation in which the ambiguity sets are coupled only through functional Lagrange multipliers. Restricting these multipliers to affine decision rules yields upper bounds that decompose across the ambiguity sets, are exact for quadratic losses over Wasserstein balls around Gaussian references, and provably dominate aggregation-based relaxations such as Minkowski sums and balls around mixtures. The approximations admit finite conic reformulations in important special cases and subgradient methods for weakly convex losses, which makes fusion DRO applicable to neural network training. In portfolio, newsvendor, and binary classification experiments, the affine approximation stays close to the exact value wherever the latter is computable, scales far better than exact methods, and outperforms aggregation-based benchmarks out-of-sample. 
}

\Template{}{Efficient Data Fusion in DRO}{}
\begin{document}

\maketitle

\section{Introduction}
\label{sec:intro}
Decision makers rarely have enough data of their own and therefore increasingly draw on additional sources. A lender may complement its own loan records with the noisy aggregate statistics released by a credit bureau, and a small community hospital may complement its few local records with data from a large academic medical center in its network. Such sources are informative but imperfect, as the statistics of the bureau hide the individual records and the patients of the medical center differ from those of the community hospital. No single source thus describes the distribution that governs the decision, yet each one constrains it. Such fragmented evidence takes many forms, including samples from related but different distributions \citep{taskesen2021sequential,rychener2024wasserstein,lixianyu2026robust}, privacy-preserving or noisy surrogates \citep{selvi2025distributionally}, complementary datasets \citep{awasthi2022datajoin}, and a single dataset that admits several statistical descriptions, such as distances and moment constraints \citep{gao2017distributionally,tanoumand2023intersecting,wang2024covariate}. This paper asks how to make reliable decisions from such evidence.\looseness=-1

The decision maker seeks a decision $\bm{\beta}$ in a nonempty closed set $\mathcal{B}\subseteq\mathbb{R}^n$ that minimizes the expected loss $\mathbb{E}_{\mathbb{P}}[\ell(\bm{\beta};\bm{\xi})]$, where the uncertain parameter $\bm{\xi}$ ranges over a nonempty closed set $\Xi\subseteq\mathbb{R}^m$ and follows a distribution $\mathbb{P}\in\mathcal{P}_0(\Xi)$. As $\mathbb{P}$ is unknown, distributionally robust optimization (DRO) replaces it with an ambiguity set $\mathcal{P}$ of plausible distributions and solves
\begin{align}\label{eq:classical_dro}
\inf_{\bm{\beta}\in\mathcal{B}}\sup_{\mathbb{Q}\in\mathcal{P}}\mathbb{E}_{\mathbb{Q}}[\ell(\bm{\beta};\bm{\xi})].
\end{align}
Classical DRO constructs~$\mathcal{P}$ as a Wasserstein ball~\citep{mohajerin2018data,zhao2018data} or a $\phi$-divergence ball~\citep{ben2013robust,duchi2021learning} around a single reference distribution, typically the empirical distribution of the training data, or defines it through moment constraints~\citep{delage2010distributionally}. Fragmented evidence does not fit these templates: each piece of evidence naturally induces its own ambiguity set $\mathfrak{B}_k\subseteq\mathcal{P}_0(\Xi)$, $k\in[K]$, and the distributions consistent with all of them form the \emph{fusion ambiguity~set}\begin{align}\label{eq:intersection-set}
    \mathcal{P}=\cap_{k=1}^K\mathfrak{B}_k.
\end{align}
We refer to the sets~$\mathfrak{B}_k$, $k\in[K]$, as \emph{elementary ambiguity sets} and to~\eqref{eq:classical_dro} with their intersection as the \emph{fusion DRO problem}. Typically, each elementary set is a neighborhood of the empirical distribution of one data source, whose type is matched to the imperfections of that source. Samples corrupted by bounded measurement noise suggest an $\infty$-Wasserstein ball whose radius bounds the noise. Samples contaminated by outliers suggest a total variation ball whose radius reflects the fraction of outliers, and samples from a shifted distribution suggest a $2$-Wasserstein ball whose radius reflects the size of the shift relative to the target distribution~$\mathbb{P}$. Other sets may constrain moments or marginals directly. 
Taking the intersection of all elementary ambiguity sets yields less conservative decisions. 

However, the price of data fusion is computational. Tractable reformulations of worst-case expectation problems are known for many elementary ambiguity sets, but not for  fusion ambiguity sets.
Suppose, for example, that the elementary ambiguity sets are Wasserstein balls around empirical distributions with $N_1,\ldots,N_K$ atoms, respectively. Known exact reformulations need to  enumerate all combinations of atoms across the empirical distributions, and involve a  number of  decision variables and constraints that grows as $\prod_{k=1}^KN_k$~\citep{rychener2024wasserstein}. Worse, for $2$-Wasserstein balls on~$\Xi=\mathbb{R}^m$, deciding whether their intersection is nonempty is NP-hard when both $K$ and $m$ are part of the input, so that even evaluating the worst-case expectation for a fixed decision is intractable in general~\citep[theorem~4]{rychener2024wasserstein}. 

As fusion DRO is intractable in general, we seek instead safe approximations that retain the fusion structure. We first derive an exact reformulation whose objective is a sum of $K$ worst-case expectations, one over each elementary ambiguity set, while all interaction between the sets is moved into a single constraint involving functional Lagrange multipliers. Restricting these multipliers to affine or static decision rules simplifies this constraint, leaving $K$ tractable worst-case expectation problems that are linked only by finitely many linear constraints on the decision rule coefficients. The resulting upper bounds on the worst-case expected loss are exact in canonical settings, provably tighter than aggregation-based alternatives, and amenable to first-order methods even for nonconvex losses involving neural networks. In our experiments, the affine bound stays close to the exact value wherever the latter can be computed and remains tractable far beyond that point, and it roughly halves the suboptimality of models trained over a single elementary ambiguity set when learning a classifier from noisy moments and histograms. 

Our main contributions can be summarized as follows.

\begin{itemize}
    \item In Section~\ref{sec:main}, we recast fusion DRO exactly as a problem in which the $K$ ambiguity sets are coupled only through functional Lagrange multipliers. Restricting them to affine or static decision rules yields tractable upper bounds. For $2$-Wasserstein balls around references from a common elliptical family, such as Gaussians, the affine bound is exact for quadratic losses and the static bound for affine losses if the references share their covariance matrices.
    \item In Section~\ref{sec:comparison}, we show that our bounds are provably tighter than aggregation-based relaxations that replace the intersection by a Minkowski sum or by a single ball around a mixture of the references, and that a barycenter-based construction admits no such guarantee. For $1$-Wasserstein balls around a common reference, the Minkowski relaxation reduces to elastic net regularization, which our approximations improve on.
    \item In Section~\ref{sec:sgd}, we derive finite conic reformulations for fusion DRO problems with discrete reference distributions and develop subgradient methods with convergence guarantees for weakly convex losses over Wasserstein balls. We also explain how these subgradient methods extend to $\phi$-divergence and moment ambiguity sets.
    \item In Section~\ref{sec:experiments}, we demonstrate on portfolio optimization, multi-item newsvendor and neural network classification problems that our approximations are tight as well as scalable and that they outperform aggregation-based benchmarks out of sample.
\end{itemize}

\paragraph{Related Work.}

When several sources of information bear on a decision, the literature combines them in four ways. One can intersect the elementary ambiguity sets they induce as in fusion DRO, average or enclose the sources in a single ambiguity set, let separate ambiguity sets constrain different but linked distributions, or seek good performance across all sources. We review these alternatives in turn and then discuss the methodological roots of our approach.

\textit{Intersecting ambiguity sets.} Fusion DRO has so far only been studied for specific types of elementary ambiguity sets and for loss functions that allow for exact convex reformulations. Intersections of $K=2$ Wasserstein balls have been used for robust data joining \citep{awasthi2022datajoin} and for mitigating robust overfitting \citep{selvi2025distributionally}, both with logistic loss, and for contextual optimization under covariate shift with piecewise concave loss \citep{wang2024covariate}. A Wasserstein ball has been intersected with a moment ambiguity set in portfolio selection \citep{gao2017distributionally} and resource allocation \citep{luo2023wasserstein}, goodness-of-fit confidence regions have been intersected with each other or with Wasserstein balls when $K=2$ or the distributions are univariate \citep{tanoumand2023intersecting}, and two moment-based sets around heterogeneous empirical distributions have been intersected for domain adaptation \citep{taskesen2021sequential}. Closest to our work, \citet{rychener2024wasserstein} and \citet{li2026harnessing} intersect $K\geq 2$ optimal transport balls around the empirical distributions of heterogeneous data sources. \citet{rychener2024wasserstein} derive exact convex reformulations for piecewise concave losses, whose sizes grow as $\prod_{k=1}^KN_k$, establish the NP-hardness discussed above, and show that fusion DRO is stable under perturbations of the data, unlike approaches that center a single ball at a Wasserstein barycenter of the sources. Our work complements these studies in three respects. It handles any number of elementary ambiguity sets of possibly different types within one framework. It replaces exact reformulations of exponential size by approximations that are exact in canonical settings, provably tighter than aggregation-based relaxations, and of a size that grows only linearly in the total number of atoms. Finally, it extends to losses that are nonconvex in the decision, such as the training loss of a neural network, for which our subgradient methods compute approximate stationary points.

\textit{Aggregating evidence into one ambiguity set.} A second school of thought combines the sources into a single ambiguity set by averaging or enclosing them, whereas fusion DRO intersects the elementary ambiguity sets they induce. Averaging can take place at the level of samples, distributions or ambiguity sets. The simplest form is pooling, which \citet{gupta2022data} apply without robustification by shrinking the empirical distributions of many related problems toward their pooled distribution. \citet{keehan2025anger} and \citet{guo2026optimization} instead center a Wasserstein ball at a weighted empirical distribution, which downweights older samples to account for temporal drift or less accurate sources, respectively. \citet{ma2024bayesian} average at the level of ambiguity sets, combining Wasserstein balls through a weighted Minkowski sum, which contains all mixtures of distributions from the individual balls. When no samples from the target distribution~$\mathbb{P}$ are available, \citet{han2026robust} and \citet{lixianyu2026robust} instead enclose the source distributions in a single ambiguity set, namely a Wasserstein ball with estimated center and radius and a ball of kernel mean embeddings in a reproducing kernel Hilbert space, respectively. Section~\ref{sec:comparison} shows that the Minkowski sum and a ball around the mixture of the sources are provably more conservative than our decision rule approximations of fusion DRO.

\textit{Linking several distributions.} A third school of thought uses several ambiguity sets that constrain different, linked distributions rather than a single one. To protect against data corruption, \citet{bennouna2023certified,bennouna2025holistic} chain L\'evy-Prokhorov and Kullback-Leibler balls through an intermediate distribution, and \citet{nietert2023outlier} let the center of a Wasserstein ball vary within a total variation neighborhood of the data. Similarly, \citet{vanparys2024efficient} links noisy observations to candidate clean distributions through a known noise process, \citet{estebanperez2022distributionally} center Wasserstein balls at reweighted empirical distributions to infer conditional distributions, \citet{nguyen2022robust} place separate balls around class-conditional distributions, and \citet{shapiro2023bayesian} average the worst-case losses over balls around candidate models with posterior weights. These models describe structured relationships between distributions, whereas fusion DRO describes several pieces of evidence about one distribution.

\textit{Robustness across sources and transfer learning.} A fourth school of thought seeks good performance across several source distributions or subpopulations simultaneously, whereas fusion DRO uses all sources to narrow down a single target distribution~$\mathbb{P}$. Examples include optimizing the worst case over mixtures of sources \citep{mohri2019agnostic,liu2026distributionally} or over subpopulations \citep{hashimoto2018fairness,sagawa2020group}, and combining source predictors so that they perform well under any mixture of the sources \citep{mansour2008domain}. Transfer learning and domain adaptation, in contrast, target a single distribution but typically require samples from it, which are used to reweight source data or learn shared representations. The underlying theory bounds the target risk in terms of a source risk and a discrepancy between the two distributions \citep{ben2010theory}. Fusion DRO can be viewed as a robust counterpart of these paradigms that hedges against all target distributions within prescribed discrepancies of every source simultaneously, and it applies even when the sources provide only summary statistics. Meta-learning \citep{finn2017model} requires many related tasks to learn how to adapt to a new one, whereas fusion DRO works with a handful of sources. Finally, the term data fusion is also used in causal inference for combining heterogeneous studies \citep{bareinboim2016causal}, where the focus is on identification rather than robust decision making.

\textit{Methodological foundations.} Our approximations restrict functional Lagrange multipliers to affine or static decision rules, an idea that originates from adjustable robust optimization and stochastic programming. There, affine policies \citep{bental2004adjustable}, primal and dual linear decision rules \citep{kuhn2011primal} and lifted decision rules \citep{georghiou2015generalized} trade optimality for tractability, and DRO has adopted them for two-stage problems over moment-based ambiguity sets \citep{goh2010distributionally,bertsimas2019adaptive}. In all these works, decision rules approximate recourse decisions that adapt to the uncertain parameters. Fusion DRO has no recourse. Instead, the decision rules approximate the multipliers of the constraints that link the elementary ambiguity sets, and restricting them decouples the sets. The resulting approximations are tractable whenever, for each elementary ambiguity set, the worst-case expectation of the loss plus a quadratic function is tractable. Similar to the liftings of \citet{bertsimas2019adaptive}, which introduce auxiliary variables for nonlinear functions appearing in the ambiguity set, our decision rules may also depend on functions whose expectations the elementary sets constrain.

Our exact reformulation rests on a strong duality result for fusion DRO with convex elementary ambiguity sets, which complements the duality theory of \citet{rychener2024wasserstein} for Wasserstein balls; see \citet{kuhn2025distributionally} for a general survey of DRO. For Wasserstein balls around empirical distributions, the exact dual problem has one semi-infinite constraint for each tuple consisting of one atom per empirical distribution. To the best of our knowledge, this problem cannot be cast as an empirical risk minimization problem, and it is unclear whether stochastic gradient methods can address it. Our decision rule approximations involve no such constraints. Their objective function is a sum of expectations, one per elementary ambiguity set, linked only by finitely many linear constraints, which makes them amenable to subgradient methods for weakly convex losses \citep{namkoong2016stochastic,davis2019stochastic}.

Conceptually, intersecting the sets of distributions consistent with several sources corresponds to the conjunctive combination of credal sets in imprecise probability theory, and minimizing the worst-case expected loss over the resulting set corresponds to its $\Gamma$-minimax decision criterion \citep{walley1991statistical}, which is closely related to maxmin expected utility \citep{gilboa1989maxmin}. Fusion DRO thus shares its modeling principle with this literature but focuses on computation.

\paragraph{Notation.} For $K\in\mathbb N$, we write $[K]:=\{1,\ldots,K\}$ and define $\Delta(K):=\{\bm{t}\in\bR_+^K: \sum_{k=1}^{K}t_k=1\}$ as the probability simplex in $\mathbb{R}^K$. For a Borel set $\Xi\subseteq\mathbb{R}^m$, we write $\mathcal{P}_0(\Xi)$ for the set of Borel probability measures supported on $\Xi$, and $C(\Xi)$ for the space of continuous real-valued functions on $\Xi$. The expectation operator with respect to $\mathbb{Q}\in\cP_0(\Xi)$ is denoted by $\mathbb{E}_{\mathbb{Q}}[g(\bm{\xi})]$ for any Borel function $g:\Xi\to\mathbb{R}$. Throughout the paper, we adopt the inf-addition convention whereby $\infty-\infty:=\infty$. This implies, in particular, that $\mathbb{E}_{\mathbb{Q}}[g(\bm{\xi})]= \infty$ whenever both the positive and the negative parts of~$g$ integrate to~$\infty$. For a Borel map $T:\mathbb R^m\to\mathbb R^m$, we use $T_{\#}\mathbb{P}$ to denote the pushforward of $\mathbb{P}$ under $T$. For functions $f,g:\mathbb{R}_{++}\to\mathbb{R}_+$, we write $f(\delta)=\mathcal{O}(g(\delta))$ if there exist $C,\delta_0>0$ such that $f(\delta)\leq Cg(\delta)$ for all $\delta\in(0,\delta_0)$. For a function $h:\mathbb{R}^d\to\mathbb{R}\cup\{+\infty\}$, we define its convex conjugate  $h^*:\mathbb R^d\to\mathbb R\cup\{-\infty,\infty\}$ via $h^*(\bm{z})\coloneqq\sup_{\bm{x}\in\mathbb{R}^d}\{\bm{z}^\top\bm{x}-h(\bm{x})\}$.

\section{Decision Rule Approximations for Fusion DRO}
\label{sec:main}
The central challenge in the DRO problem~\eqref{eq:classical_dro} with fusion ambiguity set~\eqref{eq:intersection-set} lies in the inner worst-case expectation. Indeed, evaluating the inner supremum in~\eqref{eq:classical_dro} amounts to solving an infinite-dimensional optimization problem over probability measures subject to multiple complicated constraints. Our approach is to approximate the dual of this inner problem tractably and then to embed the result back into the outer minimization problem over $\bm{\beta}$. To this end, we fix $\bm{\beta}$ and absorb it into the loss, reducing the analysis to the (primal) \emph{uncertainty quantification problem}
 \begin{align}\label{prob:sup-primal} \tag{P}
    v_\mathrm{P} \coloneqq \sup_{\mathbb{Q}\in\mathcal{P}}\,\mathbb{E}_{\mathbb{Q}}[\ell(\bm{\xi})],
\end{align}
where $\ell$ is a fixed Borel measurable loss function. In order to rule out pathological situations of no practical interest, we assume throughout that each elementary ambiguity set~$\mathfrak{B}_k$, $k\in [K]$,  contains at least one distribution~$\mathbb{Q}$ with $\mathbb{E}_{\mathbb{Q}}[\ell(\bm{\xi})]>-\infty$. In the remainder, we develop tractable approximations of~\eqref{prob:sup-primal} when each $\mathfrak B_k$ in~\eqref{eq:intersection-set} belongs to a (possibly different) family of standard ambiguity sets. A leading example is the family of discrepancy-based ambiguity sets of the form
\begin{align}\label{eq:discrepancy-ball}
    \mathfrak{B}_k = \bigl\{\mathbb{Q}\in\mathcal{P}_0(\Xi) : d_k(\mathbb{Q},\hat{\mathbb{P}}_k)\leq \varepsilon_k\bigr\},
\end{align}
where $\hat{\mathbb{P}}_k\in\mathcal{P}_0(\Xi)$ is a reference distribution, $\varepsilon_k\geq 0$ is a radius, and $d_k:\mathcal{P}_0(\Xi)\times\mathcal{P}_0(\Xi)\to[0,+\infty]$ is a discrepancy measure. We require that $d_k(\mathbb{Q},\hat{\mathbb{P}}_k) = 0$ if and only if $\mathbb{Q}=\hat{\mathbb{P}}_k$ and that $d_k(\mathbb{Q},\hat{\mathbb{P}}_k)$ is quasiconvex and weakly lower semicontinuous in~$\mathbb{Q}$. Together these mild conditions ensure that $\hat{\mathbb{P}}_k\in\mathfrak{B}_k$ for every $\varepsilon_k\geq 0$ and that $\mathfrak{B}_k$ is a convex as well as weakly closed ambiguity set. Wasserstein and $\phi$-divergence ambiguity sets are canonical examples satisfying these conditions.

\begin{example}[Wasserstein ambiguity sets]
\label{ex:wasserstein-ambiguity-set}
For $p\in[1,\infty)$ and a norm $\|\cdot\|$ on $\mathbb{R}^m$, the $p$-Wasserstein distance between $\mathbb{Q},\mathbb{Q}'\in\mathcal{P}_0(\Xi)$ is
\begin{align}\label{eq:wasserstein-dist}
    \mathsf{W}_p(\mathbb{Q},\mathbb{Q}') \coloneqq \left(\inf_{\Pi\in\mathcal{C}(\mathbb{Q},\mathbb{Q}')} \mathbb{E}_{\Pi}\bigl[\|\bm{\xi}-\bm{\xi}'\|^p\bigr]\right)^{1/p},
\end{align}
where $\mathcal{C}(\mathbb{Q},\mathbb{Q}')$ denotes the set of couplings of $\mathbb{Q}$ and $\mathbb{Q}'$, that is, the joint distributions on $\Xi\times\Xi$ with marginals $\mathbb{Q}$ and $\mathbb{Q}'$~\citep{villani2008optimal}. The map $\mathsf{W}_p$ is a metric on $\mathcal{P}_0(\Xi)$ and satisfies the conditions on $d_k$ imposed above. Moreover, $\mathsf{W}_p^p(\mathbb{Q},\mathbb{Q}')$ is convex and weakly lower semicontinuous in~$\mathbb{Q}$ \citep[theorem~4.8 and remark~6.12]{villani2008optimal}, and thus $\mathsf{W}_p$ is quasiconvex in~$\mathbb{Q}$. Setting $d_k=\mathsf{W}_p$ in~\eqref{eq:discrepancy-ball} yields the $p$-Wasserstein ambiguity set. Analogously, one can construct weakly closed $\infty$-Wasserstein ambiguity sets \citep[section~2.3.4]{kuhn2019wasserstein}.
\end{example}

\begin{example}[$\phi$-divergence ambiguity sets]
\label{ex:phi-divergence-ambiguity-set}
For a closed convex function $\phi:\mathbb{R}_+\to\mathbb{R}\cup\{+\infty\}$ with $\phi(1)=0$ and $\lim_{t\to\infty} \phi(t)/t=\infty$, the $\phi$-divergence of $\mathbb{Q}\in\mathcal{P}_0(\Xi)$ from $\mathbb{Q}'\in\mathcal{P}_0(\Xi)$ is
\begin{align*}
    \mathsf{D}_\phi(\mathbb{Q}\|\mathbb{Q}') \coloneqq \begin{cases} \displaystyle\int_\Xi \phi\!\left(\frac{\mathrm{d}\mathbb{Q}}{\mathrm{d}\mathbb{Q}'}(\bm{\xi})\right)\mathrm{d}\mathbb{Q}'(\bm{\xi}) & \text{if } \mathbb{Q}\ll\mathbb{Q}', \\ +\infty & \text{otherwise,}\end{cases}
\end{align*}
where $\mathrm{d}\mathbb{Q}/\mathrm{d}\mathbb{Q}'$ denotes the Radon-Nikodym derivative~\citep{csiszar1963phidiv}. Prominent special cases include the Kullback-Leibler divergence ($\phi(t)=t\log t$) or the $\chi^2$-divergence ($\phi(t)=(t-1)^2$). The map $\mathsf{D}_\phi(\mathbb{Q}\|\mathbb{Q}')$ is convex and weakly lower semicontinuous and satisfies $\mathsf{D}_\phi(\mathbb{Q}\|\mathbb{Q}')=0$ if and only if $\mathbb{Q}=\mathbb{Q}'$ under mild conditions on~$\phi$~\citep{csiszar1963phidiv}. Setting $d_k=\mathsf{D}_\phi$ in~\eqref{eq:discrepancy-ball} yields the $\phi$-divergence ambiguity set. We emphasize that the restrictive condition $\lim_{t\to\infty} \phi(t)/t=\infty$ is not needed, that is, one can construct more general weakly closed $\phi$-divergence ambiguity sets at the expense of introducing additional notation \citep[section~2.2]{kuhn2025distributionally}.
\end{example}

We also allow the primitive distribution families $\mathfrak B_k$ in~\eqref{eq:intersection-set} to be moment ambiguity sets.

\begin{example}[Moment ambiguity sets]
\label{ex:moment-ambiguity-sets}
Writing $\mathbb{Q}\sim(\bm{\mu},\bm{\Sigma})$ to mean that a distribution~$\mathbb{Q}$ on $\mathbb R^m$ has mean $\bm{\mu}\in\mathbb R^m$ and covariance matrix $\bm{\Sigma}\in\mathbb S^m_+$, a moment ambiguity set takes the form
\begin{align*}
    \mathfrak{B}_k = \bigl\{\mathbb{Q}\in\mathcal{P}_0(\Xi) : \mathbb{Q}\sim(\bm{\mu},\bm{\Sigma}),\;(\bm{\mu},\bm{\Sigma})\in\mathcal{U}_k\bigr\},
\end{align*}
where $\mathcal{U}_k\subseteq\mathbb{R}^m\times\mathbb{S}_+^m$ is a closed moment uncertainty set. A singleton $\mathcal{U}_k=\{(\hat{\bm{\mu}}_k,\hat{\bm{\Sigma}}_k)\}$ yields the Chebyshev ambiguity set, while $\mathcal{U}_k=\{\hat{\bm{\mu}}_k\}\times\mathbb{S}_+^m$ yields the Markov ambiguity set, see, e.g., \citep[section~2.1]{kuhn2025distributionally}. A popular choice takes $\mathcal{U}_k$ to be a ball in the Gelbrich distance,
\begin{align*}
    \mathcal{U}_k = \bigl\{(\bm{\mu},\bm{\Sigma})\in\mathbb{R}^m\times\mathbb{S}_+^m : \mathsf{G}^2\bigl((\bm{\mu},\bm{\Sigma}),(\hat{\bm{\mu}}_k,\hat{\bm{\Sigma}}_k)\bigr)\leq\varepsilon_k^2\bigr\},
\end{align*}
where the squared Gelbrich distance between two mean-covariance pairs is defined as
\begin{align}\label{eq:gelbrich-dist}
    \mathsf{G}^2\bigl((\bm{\mu},\bm{\Sigma}),(\hat{\bm{\mu}},\hat{\bm{\Sigma}})\bigr)\coloneqq\|\bm{\mu}-\hat{\bm{\mu}}\|_2^2+\mathrm{Tr}\Bigl[\bm{\Sigma}+\hat{\bm{\Sigma}}-2\bigl(\hat{\bm{\Sigma}}^{1/2}\bm{\Sigma}\hat{\bm{\Sigma}}^{1/2}\bigr)^{1/2}\Bigr],
\end{align}
see \citep{gelbrich1990formula}. Unlike $\mathsf{W}_p$ and $\mathsf{D}_\phi$, the Gelbrich distance compares moment pairs rather than distributions and satisfies $\mathsf{G}((\bm{\mu},\bm{\Sigma}),(\hat{\bm{\mu}},\hat{\bm{\Sigma}}))=0$ whenever $\bm{\mu}=\hat{\bm{\mu}}$ and $\bm{\Sigma}=\hat{\bm{\Sigma}}$, regardless of whether the underlying distributions coincide. It therefore does not satisfy the identity of indiscernibles required of $d_k$, and $\mathfrak{B}_k$ is not a special case of~\eqref{eq:discrepancy-ball}.
\end{example}
 
Several of our subsequent results make use of the following regularity conditions. 

\begin{assumption}[Regularity]~\label{assump:general}
    \begin{enumerate}[label=(\roman*)]
     \item \label{assump:general-support} The support set $\Xi\subseteq \bR^m$ is compact.
    \item \label{assump:general-convex-ball} Each ambiguity set $\mathfrak{B}_{k}\subseteq \cP_0(\Xi)$, $k\in[K]$, is convex and weakly closed. 
    \item \label{assump:general-nonempty-intersection} The intersection $\cap_{k=1}^K\mathfrak{B}_k$ is nonempty. 
    \item \label{assump:loss-usc} The loss function $\ell:\Xi\to \mathbb R$ is upper semicontinuous.
\end{enumerate}
\end{assumption}

Together, the conditions in Assumption~\ref{assump:general} ensure that Problem~\eqref{prob:sup-primal} is feasible and that the supremum $v_\mathrm{P}$ is attained. To see this, note first that Assumptions~\ref{assump:general}\,\ref{assump:general-support} and~\ref{assump:general}\,\ref{assump:general-convex-ball} imply that each~$\mathfrak{B}_k$ is weakly compact~\cite[proposition~3.6]{kuhn2025distributionally}. For discrepancy balls~\eqref{eq:discrepancy-ball}, this follows from the conditions imposed on~$d_k$, while for moment ambiguity sets, Assumption~\ref{assump:general}\,\ref{assump:general-support} alone suffices \cite[proposition~3.7]{kuhn2025distributionally}. Assumption~\ref{assump:general}\,\ref{assump:general-nonempty-intersection} then ensures that~$\mathcal{P}$ is nonempty, while Assumptions~\ref{assump:general}\,\ref{assump:general-support} and~\ref{assump:general}\,\ref{assump:loss-usc} together imply that the expected loss $\mathbb{E}_{\mathbb{Q}}[\ell(\bm{\xi})]$ is weakly upper semicontinuous in~$\mathbb Q$ \cite[proposition~3.3]{kuhn2025distributionally}. The supremum in~\eqref{prob:sup-primal} is therefore attained by the extreme value theorem, which applies because $\mathcal{P}$ is weakly compact as the intersection of the weakly compact sets~$\mathfrak{B}_k$, $k\in[K]$. The following proposition shows that Problem~\eqref{prob:sup-primal} admits a strong dual that decouples across the $K$ ambiguity sets, replacing the joint optimization over $\mathcal{P}$ with $K$ independent subproblems over the individual sets~$\mathfrak{B}_k$.

\begin{proposition}[Duality]
\label{prop:reformulation_main}
    Problem~\eqref{prob:sup-primal} admits the following weak dual.
    \begin{align}
    \label{eq:reformulation_dual} \tag{D}
     v_\mathrm{D}:=\left\{\begin{array}{lll}
       \displaystyle {\inf} & \displaystyle \sum_{k=1}^K\underset{\mathbb{Q}_k \in \mathfrak{B}_{k}}{\sup} \mathbb{E}_{\mathbb{Q}_k}[\lambda_k(\bm{\xi})]  \\
         \mathrm{s.t.} & \lambda_k \in  C(\Xi) & \forall k\in[K] \\
         &  \ell( \bm{\xi}) -\sum_{k=1}^K \lambda_k(\bm{\xi})  \leq 0 &  \forall \bm{\xi} \in \Xi
    \end{array}\right.
\end{align}
Thus, we have $v_\mathrm{P}\leq v_\mathrm{D}$. If Assumption~\ref{assump:general} holds, then we have strong duality, that is, $v_\mathrm{P}=v_\mathrm{D}$.
\end{proposition}

Assumption~\ref{assump:general} is sufficient but not necessary for strong duality. One can show that Proposition~\ref{prop:reformulation_main} continues to hold in other settings, for instance when $\Xi=\mathbb{R}^m$ under appropriate growth conditions. The proposition replaces the joint optimization over the fusion ambiguity set~$\mathcal{P}$ with a dual problem whose objective decouples into $K$ worst-case expectation terms, each involving only an elementary ambiguity set $\mathfrak{B}_k$. The remaining difficulty lies in the decision variables $\lambda_k\in C(\Xi)$, which are arbitrary continuous functions coupled through the semi-infinite constraint $\sum_{k=1}^K\lambda_k(\bm{\xi})\geq\ell(\bm{\xi})$ for all $\bm{\xi}\in\Xi$. A natural approach would be to address~\eqref{eq:reformulation_dual} with a decision rule approximation, which imposes a parametric restriction on each~$\lambda_k$. However, even though the objective separates across $k$, the constraint ties all $K$ functions to the same loss value $\ell(\bm{\xi})$ at the same point $\bm{\xi}$, so restricting each $\lambda_k$ independently still constrains how well their sum approximates $\ell$ pointwise. We therefore pursue a different strategy. We introduce auxiliary copies~$\bm{\xi}_1,\ldots,\bm{\xi}_K$ of the uncertainty $\bm{\xi}$ and rewrite the constraint as $\sum_{k=1}^K\lambda_k(\bm{\xi}_k)\geq\ell(\bm{\xi})$ subject to the linking constraints $\bm{\xi}_k=\bm{\xi}$ for all $k\in[K]$. This reformulation is equivalent to~\eqref{eq:reformulation_dual}, but dualizing the linking constraints separates the resulting Lagrangian across~$k$ in the sense that each $\lambda_k$ is evaluated only at its own copy $\bm{\xi}_k$, and the coupling is absorbed into the dual multipliers associated with $\bm{\xi}_k=\bm{\xi}$. This yields an exact reformulation of~\eqref{eq:reformulation_dual} with a fully separable structure, enabling richer and more tractable classes of decision rule approximations.

\begin{lemma}\label{lem:pairwise-symmetric-rules}
    For any weight vector $\bm{t}\in\Delta(K)$, Problem~\eqref{eq:reformulation_dual} is equivalent to
    \begin{equation}
    \begin{aligned}
    \label{eq:prob-with-semi-infinite-constraint}
    \inf \quad & \sum_{k=1}^K \sup_{\mathbb{Q}_k\in\mathfrak{B}_k} \mathbb{E}_{\mathbb{Q}_k}[\lambda_k(\bm{\xi})] \\ \mathrm{s.t.} \quad & \lambda_k\in C(\Xi) \quad \forall k\in[K] \\ 
    & \hspace{-2.5mm} \left. \begin{array}{c@{\;}l}
        \displaystyle \inf_{\{\bm{\alpha}_k\in\mathcal{A}\}_k} & \displaystyle \sup_{\bm{\xi}\in\Xi} \; \sup_{\{\bm{\xi}_k\in\Xi\}_k} \; \sum_{k=1}^K \Bigl(t_k\ell(\bm{\xi}_k) - \lambda_k(\bm{\xi}_k) + \bm{\alpha}_k(\bm{\xi}_k,\bm{\xi})^\top(\bm{\xi}_k-\bm{\xi})\Bigr) \\[2ex]
        \mathrm{s.t.} & \displaystyle \sum_{k=1}^K\bm{\alpha}_k(\bm{\xi}',\bm{\xi}) =\bm{0} \quad \forall \bm \xi,\bm \xi'\in\Xi
        \end{array}\right\} \leq 0 ,
    \end{aligned}
    \end{equation}
    where $\mathcal A$ denotes the space of all functions $\bm\alpha:\Xi\times\Xi\to\mathbb R^m$ that are symmetric in their arguments, that is, $\bm{\alpha}(\bm{\xi}',\bm{\xi}) =\bm{\alpha}(\bm{\xi},\bm{\xi}')$. Consequently, the optimal value of~\eqref{eq:prob-with-semi-infinite-constraint} coincides with~$v_\mathrm{D}$. If additionally Assumption~\ref{assump:general} holds, then this optimal value also coincides with~$v_\mathrm{P}$.
\end{lemma}

Problem~\eqref{eq:prob-with-semi-infinite-constraint} is still intractable due to the infinite-dimensional multiplier functions $\bm{\alpha}_k\in\mathcal{A}$. This places~\eqref{eq:prob-with-semi-infinite-constraint} squarely in the framework of two-stage robust optimization, where each function~$\bm{\alpha}_k$, $k\in[K]$, plays the role of an adjustable decision \citep[section~14]{bental2009robust}. While such problems are generically NP-hard \cite[example~3.1]{bental2004adjustable}, tractable approximations can be obtained by restricting all adjustable decisions to structured decision rules. This prompts us to restrict $\bm{\alpha}_k$ to an affine function of the form
\begin{align}\label{eq:affine-decision-rule}
    \bm{\alpha}_k(\bm{\xi}',\bm{\xi}) = \bm{a}_k + \bm{A}_{k}(\bm{\xi}'+\bm{\xi})
\end{align}
for some coefficient vectors $\bm{a}_k\in\bR^m$ and matrices $\bm{A}_k\in\bR^{m\times m}$ for all~$k\in[K]$. The dependence on $\bm{\xi}'+\bm{\xi}$ makes the decision rule~\eqref{eq:affine-decision-rule} permutation symmetric, as required by the definition of~$\mathcal A$.

Next, we derive a useful reformulation of the resulting decision rule approximation.

\begin{proposition}[Affine Approximation]\label{prop:affine_upper_bound}
Suppose that~$\ell$ is continuous and fix any $\bm{t}\in\Delta(K)$. Restricting the functions $\{\bm{\alpha}_k\}_k$ in Problem~\eqref{eq:prob-with-semi-infinite-constraint} to affine decision rules of the form~\eqref{eq:affine-decision-rule} reduces~\eqref{eq:prob-with-semi-infinite-constraint}~to
\begin{align}\label{prob:affine}\tag{D$_\mathrm{aff}$}
      v_\mathrm{aff} 
        \coloneqq \left\{\begin{array}{lll} \inf & \displaystyle \sum_{k=1}^K \underset{\mathbb{Q}_k \in \mathfrak{B}_{k}}{\sup} \mathbb{E}_{\mathbb{Q}_k}[t_k \ell(\bm{\xi}) + {\bm{a}_k}^\top \bm{\xi} +\langle \bm{A}_k,\bm{\xi} \bm{\xi}^\top \rangle ] \\
         \mathrm{s.t.} & \bm{a}_k\in \mathbb{R}^m, \; \bm{A}_k\in \bS^m \quad \forall k\in[K] \\
         & \sum_{k=1}^K \bm{a}_k =\mathbf{0} \\
         & \sum_{k=1}^K \bm{A}_k =\mathbf{0},
    \end{array}\right.
\end{align}
and therefore $v_\mathrm{P}\leq v_\mathrm{aff}$. Moreover, for any fixed decision rules of the form~\eqref{eq:affine-decision-rule} with $\{\bm{a}_k,\bm{A}_k\}_k$ feasible in~\eqref{prob:affine}, the dual functions $\lambda_k(\bm{\xi}) = t_k\ell(\bm{\xi}) + \bm{a}_k^\top\bm{\xi} + \langle \bm{A}_k, \bm{\xi}\bm{\xi}^\top\rangle$, $k\in[K]$, are optimal in~\eqref{eq:prob-with-semi-infinite-constraint}.\looseness=-1
\end{proposition}

Proposition~\ref{prop:affine_upper_bound} has a key computational consequence. Rather than optimizing over the cumbersome intersection $\cap_{k=1}^K\mathfrak{B}_k$, one can bound~$v_\mathrm{P}$ from above by solving~\eqref{prob:affine}, in which the worst-case expectation over the fusion ambiguity set is replaced by~$K$ separate worst-case expectations, one per ambiguity set~$\mathfrak{B}_k$, each evaluated against a modified loss of the form
\[
    \lambda_k(\bm{\xi}) = t_k\ell(\bm{\xi}) + \bm{a}_k^\top\bm{\xi} + \langle\bm{A}_k,\bm{\xi}\bm{\xi}^\top\rangle.
\]
The~$K$ subproblems interact only through the finite-dimensional zero-sum constraints $\sum_{k=1}^K\bm{a}_k=\bm{0}$ and $\sum_{k=1}^K\bm{A}_k=\bm{0}$. Consequently, \eqref{prob:affine} is tractable whenever each reference-wise worst-case expectation admits a tractable reformulation for losses of the form $\lambda_k$ as described above.

A further simplification arises by restricting the affine decision rule~\eqref{eq:affine-decision-rule} to the constant case by setting~$\bm{A}_k=\bm{0}$ for all $k\in[K]$. This yields a static decision rule of the form
\begin{align}\label{eq:static-decision-rule}
    \bm{\alpha}_k(\bm{\xi}', \bm{\xi}) = \bm{a}_k.
\end{align}
As a direct special case of Proposition~\ref{prop:affine_upper_bound}, the following corollary is immediate.

\begin{corollary}[Static Approximation]\label{cor:static_upper_bound} 
Suppose that~$\ell$ is continuous and fix any $\bm{t}\in\Delta(K)$. Restricting the functions $\{\bm{\alpha}_k\}_k$ in Problem~\eqref{eq:prob-with-semi-infinite-constraint} to static decision rules of the form~\eqref{eq:static-decision-rule} reduces~\eqref{eq:prob-with-semi-infinite-constraint}~to
\begin{align*}\tag{D$_\mathrm{stat}$}\label{prob:static_upper_bound}
        v_\mathrm{stat} \coloneqq \left\{\begin{array}{lll} \inf & \displaystyle \sum_{k=1}^K \underset{\mathbb{Q}_k \in \mathfrak{B}_{k}}{\sup} \big\{ \mathbb{E}_{\mathbb{Q}_k}[t_k \ell( \bm{\xi}) + \bm{a}_k^\top \bm{\xi} ]  \big\} \\
         \mathrm{s.t.} & \bm{a}_k \in \mathbb{R}^m \quad \forall k\in[K] \\
         & \sum_{k=1}^K \bm{a}_k =\bm{0},
    \end{array}\right.
\end{align*}
and therefore $v_\mathrm{P} \leq v_\mathrm{aff} \leq v_\mathrm{stat}$. Moreover, for any fixed decision rules of the form~\eqref{eq:static-decision-rule} with $\{\bm{a}_k\}_k$ feasible in~\eqref{prob:affine}, the dual functions $\lambda_k(\bm{\xi}) = t_k\ell(\bm{\xi}) + \bm{a}_k^\top\bm{\xi}$, $k\in[K]$, are optimal in~\eqref{eq:prob-with-semi-infinite-constraint}.
\end{corollary}

\begin{remark}[Empty Intersections]\label{rem:empty-intersection} If $\cap_{k=1}^K\mathfrak{B}_k=\emptyset$, then $v_\mathrm{P}=-\infty$, and the fusion problem itself is ill-posed, yet $v_\mathrm{aff}$ and $v_\mathrm{stat}$ may remain finite. The bounds $v_\mathrm{P}\leq v_\mathrm{aff}\leq v_\mathrm{stat}$ hold even if Assumption~\ref{assump:general}\,\ref{assump:general-nonempty-intersection} fails, and a finite bound is no evidence that this assumption holds. The approximations detect emptiness only when it is visible in the moments. Indeed, if~$\mathfrak B_k$ is convex for every~$k\in[K]$ and a mild regularity condition holds, then one can show that $v_\mathrm{stat}=-\infty$ if and only if the sets of means of the distributions in $\mathfrak{B}_k$, $k\in[K]$, have no common point, and $v_\mathrm{aff}=-\infty$ if and only if the same holds for their sets of first and second moments. Balls that share moments but no distribution therefore go undetected. This is to be expected, as deciding whether the intersection is empty is NP-hard already for Wasserstein balls \citep[theorem~4]{rychener2024wasserstein}. Nonemptiness is, however, easily enforced by design. For discrepancy balls~\eqref{eq:discrepancy-ball}, choosing radii $\varepsilon_k\geq d_k(\mathbb{Q}_0,\hat{\mathbb{P}}_k)$ for some distribution $\mathbb{Q}_0$ and all $k\in[K]$ ensures that $\mathbb{Q}_0$ lies in the intersection. For Wasserstein balls, natural choices of $\mathbb{Q}_0$ are barycenters of the reference distributions and their cheaper Gaussian proxies (Section~\ref{sec:barycenter}), which is the calibration we adopt in Section~\ref{sec:experiments}. No choice of radii helps, in contrast, for $\phi$-divergence balls around references with disjoint supports or for Chebyshev sets with different moments (Section~\ref{sec:gradient-other}).
\end{remark}

\begin{remark}[Lifted Decision Rules]\label{rem:lifted-rules}
The decision rules~\eqref{eq:affine-decision-rule} and~\eqref{eq:static-decision-rule} can be applied to $\bm{s}(\bm{\xi})$ instead of~$\bm{\xi}$, where $\bm{s}$ is any continuous injective map. Indeed, the linking constraints $\bm{\xi}_k=\bm{\xi}$ discussed after Proposition~\ref{prop:reformulation_main} are equivalent to $\bm{s}(\bm{\xi}_k)=\bm{s}(\bm{\xi})$. Dualizing the latter yields approximations of the form~\eqref{prob:affine} and~\eqref{prob:static_upper_bound}, in which $t_k\ell(\bm{\xi})$ is augmented by the correction terms $\bm{a}_k^\top\bm{s}(\bm{\xi})+\langle\bm{A}_k,\bm{s}(\bm{\xi})\bm{s}(\bm{\xi})^\top\rangle$ and $\bm{a}_k^\top\bm{s}(\bm{\xi})$, respectively. The proofs of Proposition~\ref{prop:affine_upper_bound} and Corollary~\ref{cor:static_upper_bound} carry over with $\bm{\xi}$ replaced by $\bm{s}(\bm{\xi})$ in the correction terms, and the resulting values still bound~$v_\mathrm{P}$ from~above.
\end{remark}

We next identify conditions under which the affine approximation is exact. To this end, we adapt the notion of a positive-semidefinite affine family from~\citet[definition~6]{nguyen2021mean}.

\begin{definition}[Positive-Semidefinite Affine Family] \label{def:psd-affine-family} A set of distributions~$\mathcal{S}\subseteq\mathcal P_0(\mathbb{R}^m)$ with finite second moments is a positive-semidefinite affine family if it satisfies the following properties:
\vspace{-0.5em}
 \begin{enumerate}[label=(\roman*)]
        \item \label{def:psd-affine-closed} $T_{\#}\mathbb{Q}\in\mathcal{S}$ for every $\mathbb{Q}\in\mathcal{S}$ and every affine map $T(\bm{\xi})=\bm{B}\bm{\xi}+\bm{b}$ with $\bm{B}\in\mathbb{S}_+^m$ and $\bm{b}\in\mathbb{R}^m$; \vspace{-0.6em}
        \item \label{def:psd-affine-unique} $\mathcal{S}$ contains exactly one $\mathbb{Q}\sim(\bm{\mu},\bm{\Sigma})$ for every $(\bm{\mu},\bm{\Sigma})\in\mathbb{R}^m\times\mathbb{S}_+^m$.
    \end{enumerate}
\end{definition}

Examples of positive-semidefinite affine families include the family of all Gaussian distributions or the family of all elliptical distributions with a common characteristic generator and finite second moments~\citep[section~2.1]{nguyen2021mean}. The two conditions of Definition~\ref{def:psd-affine-family} ensure that the optimal transport map between any two non-degenerate distributions in~$\mathcal{S}$ is positive-semidefinite affine, so that the $2$-Wasserstein distance~\eqref{eq:wasserstein-dist} between them coincides with the Gelbrich distance~\eqref{eq:gelbrich-dist} between their mean-covariance pairs \citep[theorem~2]{nguyen2021mean}.

\begin{proposition}[Tightness of the Affine Approximation]\label{prop:affine-tightness}
    Suppose the following conditions hold:
    \begin{enumerate}[label=(\roman*)]
    \item $\Xi=\mathbb{R}^m$ and $\ell(\bm{\xi})=\bm{\xi}^\top\bm{P}\bm{\xi}+\bm{q}^\top\bm{\xi}+r$ for some $\bm{P}\in\mathbb{S}^m$, $\bm{q}\in\mathbb{R}^m$, $r\in\mathbb{R}$;
    \item each $\mathfrak{B}_k$ is a $2$-Wasserstein ball~\eqref{eq:discrepancy-ball} with $d_k=\mathsf{W}_2$ induced by the $2$-norm on~$\mathbb R^m$ and~$\varepsilon_k>0$; 
    \item the reference distributions $\hat{\mathbb{P}}_1,\ldots,\hat{\mathbb{P}}_K$ belong to a common positive-semidefinite affine family, and all of them have strictly positive-definite covariance matrices.
\end{enumerate}
Then $v_\mathrm{P}=v_\mathrm{aff}$, and~\eqref{prob:affine} admits a tractable semidefinite programming reformulation.
\end{proposition}

Condition~{\em (iii)} of Proposition~\ref{prop:affine-tightness} is satisfied, for instance, if all reference distributions are non-degenerate Gaussians. Note also that Proposition~\ref{prop:affine-tightness} identifies a setting in which strong duality between~\eqref{prob:sup-primal} and~\eqref{eq:reformulation_dual} holds even though Assumption~\ref{assump:general} is {\em not} satisfied. To see this, recall that weak duality gives $v_\mathrm{P}\leq v_\mathrm{D}$ by Proposition~\ref{prop:reformulation_main}, and that restricting $\bm{\alpha}_k$ to affine decision rules in the minimization problem~\eqref{eq:prob-with-semi-infinite-constraint} can only increase its optimal value, so $v_\mathrm{D}\leq v_\mathrm{aff}$. Under the conditions of Proposition~\ref{prop:affine-tightness}, $v_\mathrm{P}=v_\mathrm{aff}$, so the chain $v_\mathrm{P}\leq v_\mathrm{D}\leq v_\mathrm{aff}$ collapses to $v_\mathrm{P}=v_\mathrm{D}=v_\mathrm{aff}$, establishing strong duality even when $\Xi=\mathbb{R}^m$ is unbounded. We close this section by identifying an exactness result akin to Proposition~\ref{prop:affine-tightness} for the static approximation.

\begin{proposition}[Tightness of Static Approximation] \label{prop:static-tightness} Suppose the following conditions hold:
    \begin{enumerate}[label=(\roman*)]
    \item \label{cond:static-loss} $\Xi=\mathbb{R}^m$ and $\ell(\bm{\xi})=\bm{q}^\top\bm{\xi}+r$ for some $\bm{q}\in\mathbb{R}^m$, $r\in\mathbb{R}$;
    \item \label{cond:static-ball} each $\mathfrak{B}_k$ is a $2$-Wasserstein ball~\eqref{eq:discrepancy-ball} with $d_k=\mathsf{W}_2$ induced by the $2$-norm on~$\mathbb R^m$ and~$\varepsilon_k>0$; 
    \item \label{cond:static-family} the reference distributions $\hat{\mathbb{P}}_1,\ldots,\hat{\mathbb{P}}_K$ belong to a common positive-semidefinite affine family, and share the same covariance matrix.
\end{enumerate}
Then $v_\mathrm{P}=v_\mathrm{stat}$, and~\eqref{prob:static_upper_bound} admits a tractable second-order cone programming reformulation.
\end{proposition}

Restricting the decision rules $\bm{\alpha}_k$ in~\eqref{eq:prob-with-semi-infinite-constraint} to the affine and static classes has two notable consequences. First, Proposition~\ref{prop:affine_upper_bound} shows that the affine decision rule restriction induces dual functions of the form $\lambda_k(\bm{\xi})=t_k\ell(\bm{\xi})+\bm{a}_k^\top \bm{\xi}+\langle\bm{A}_k,\bm{\xi}\bm{\xi}^\top\rangle$. Setting $\bm{A}_k=\bm{0}$  recovers the dual functions of the static decision rule restriction as a special case. This parametrization is not imposed directly on $\lambda_k$ but is instead derived from the restriction on~$\bm{\alpha}_k$. As a result, it is automatically adapted to the loss~$\ell$, incorporating its functional form through the leading term~$t_k\ell(\bm{\xi})$. By contrast, directly restricting~$\lambda_k$ to a prescribed polynomial family, as discussed after Proposition~\ref{prop:reformulation_main}, fixes its form independently of~$\ell$ and may therefore yield a poorer approximation when~$\ell$ has structure not captured by the chosen decision rule family. We return to this comparison in Appendix~\ref{app:dual_lambda_approx}, where we show that the $\bm{\alpha}_k$-based approximations dominate the corresponding $\lambda_k$-based approximations. Second, while higher-order decision rules for~$\bm{\alpha}_k$ could in principle yield tighter upper bounds, they may render the resulting worst-case expectation problems intractable. In general, higher-order polynomial rules introduce terms of degree beyond two into the inner maximizations, which need not admit tractable reformulations. The affine and static classes strike a favorable balance. They are rich enough to achieve exactness under the conditions of Propositions~\ref{prop:affine-tightness} and~\ref{prop:static-tightness}, yet structured enough to preserve the tractability of the individual worst-case expectation problems.

\section{Comparison with Alternative Approximations}\label{sec:comparison}
Rather than approximating the hard uncertainty quantification Problem~\eqref{prob:sup-primal} directly as in Section~\ref{sec:main}, one might instead replace the fusion ambiguity set $\cap_{k=1}^K\mathfrak{B}_k$ by a single, more tractable aggregate set. We study three natural aggregation strategies based on Minkowski sums, probabilistic mixtures and Wasserstein barycenters, respectively. We then show that the first two are dominated by our static decision rule approximation, meaning that~\eqref{prob:static_upper_bound} yields tighter bounds at essentially no additional computational cost. The barycenter construction, by contrast, lacks clean dominance guarantees, is NP-hard to compute, and can be unstable. As a byproduct, we show that when all elementary ambiguity sets $\mathfrak{B}_k$, $k\in[K]$, are $1$-Wasserstein balls sharing a common reference, the Minkowski sum relaxation is reminiscent of classical elastic net regularization, and our static decision rule approximation provably improves upon it.

\subsection{Minkowski Sum Approximation}
\label{sec:minkowski}
A natural relaxation of the fusion ambiguity set~$\mathcal{P}$ is the weighted Minkowski sum of the individual ambiguity sets~$\mathfrak{B}_k$, $k\in[K]$. For a weight vector $\bm{t}\in\Delta(K)$ (the same weight vector that is used in the decision rule approximations), this set consists of all mixture distributions with fixed weights~$t_k$, where each component~$\mathbb{Q}_k$ is selected from~$\mathfrak{B}_k$. We thus approximate~$\mathcal P$ by
\begin{align*}
    \mathcal{P}_\mathrm{mink}\coloneqq\left\{\sum_{k=1}^K t_k\mathbb{Q}_k :\mathbb{Q}_k\in\mathfrak{B}_k\;\forall k\in[K]\right\}.
\end{align*}
This construction has been used in~\citep{ma2024bayesian} to model distributional ambiguity arising from heterogeneous renewable energy patterns. The corresponding worst-case risk is
\begin{align}\tag{P$_\mathrm{mink}$}\label{prob:mink}
    v_\mathrm{mink}\coloneqq\sup_{\mathbb{Q}\in\mathcal{P}_\mathrm{mink}} \mathbb{E}_{\mathbb{Q}}[\ell(\bm{\xi})] = \sum_{k=1}^K t_k \sup_{\mathbb{Q}_k\in\mathfrak{B}_k} \mathbb{E}_{\mathbb{Q}_k}[\ell(\bm{\xi})].
\end{align}
The next proposition shows that~$\mathcal{P}_\mathrm{mink}$ covers~$\mathcal{P}$ and places~$v_\mathrm{mink}$ in the approximation hierarchy.

\begin{proposition}[Minkowski Sum Relaxation]\label{prop:mink}
For any $\bm{t}\in\Delta(K)$, we have $\mathcal{P}\subseteq \mathcal{P}_\mathrm{mink}$ and
\[
    v_\mathrm{P}\leq v_\mathrm{stat}\leq v_\mathrm{mink}.
\]
\end{proposition}

Both~\eqref{prob:mink} and our decision rule approximations decompose into~$K$ independent single-reference DRO subproblems and are therefore tractable whenever each subproblem is. The key structural difference is that~\eqref{prob:mink} leaves these subproblems entirely uncoupled, whereas~\eqref{prob:static_upper_bound} and~\eqref{prob:affine} introduce coupling through the zero-sum constraints on $\{\bm{a}_k\}_k$ and $\{\bm{A}_k\}_k$. Since these constraints are finite-dimensional, they add negligible computational overhead while yielding tighter bounds, as guaranteed by Proposition~\ref{prop:mink}.

\subsection{Mixture Approximation}\label{sec:mixture}
A second approach to bounding $v_\mathrm{P}$ from above is to replace the $K$ reference distributions by the weighted mixture $\mathbb{P}_\mathrm{mix}\coloneqq\sum_{k=1}^K t_k\hat{\mathbb{P}}_k$ and center a single ambiguity ball at this aggregate nominal model, reducing~\eqref{prob:sup-primal} to a single-reference DRO problem. Unlike the Minkowski sum relaxation of Section~\ref{sec:minkowski}, which applies to arbitrary ambiguity sets, this construction requires each~$\mathfrak{B}_k$ to be a discrepancy-based set~\eqref{eq:discrepancy-ball}. Indeed, centering a ball at the mixture requires a notion of discrepancy between distributions. We thus restrict attention to such sets throughout this section and write $\mathfrak{B}_k=\mathfrak{B}^{d_k}_{\varepsilon_k}(\hat{\mathbb{P}}_k)$ to make the discrepancy, radius, and center explicit. Hence, the discrepancy~$d_k$ may change with~$k$. Placing $v_\mathrm{mix}$ in the approximation hierarchy further requires the $K$ discrepancy measures to be compatible with the mixing operation, as formalized below.

\begin{assumption}[Aggregating Discrepancy]\label{assump:aggregating-discrepancy}
There exists a discrepancy measure $d$ in the sense of Section~\ref{sec:main} such that, for every $\bm{t}\in\Delta(K)$ and all distributions $\{\mathbb{Q}_k\}_k,\{\mathbb{P}_k\}_k\subseteq\mathcal{P}_0(\Xi)$, we have
\[
    d\bigg(\sum_{k=1}^K t_k\mathbb{Q}_k,\sum_{k=1}^K t_k\mathbb{P}_k\bigg)\leq\sum_{k=1}^K t_k\,d_k(\mathbb{Q}_k,\mathbb{P}_k).
\]
\end{assumption}

The most important special case arises when all $K$ ambiguity sets share a common discrepancy $d_1=\cdots=d_K=d$, in which case Assumption~\ref{assump:aggregating-discrepancy} holds whenever~$d$ is jointly convex in its two arguments. This is the case if $d=\mathsf{W}_p^p$ is the $p$-th power of the $p$-Wasserstein distance for some~$p\in[1,\infty)$, if $d=\mathsf{D}_\phi$ is the $\phi$-divergence induced by some entropy function~$\phi$, or if $d$ is a maximum mean discrepancy \citep{gretton2012kernel}. Note that a $p$-Wasserstein ball of radius~$\varepsilon_k$ coincides with the ball of radius~$\varepsilon_k^p$ with respect to~$\mathsf{W}_p^p$, so that the mixture ambiguity set is the $p$-Wasserstein ball of radius $(\sum_{k=1}^Kt_k\varepsilon_k^p)^{1/p}$ around the mixture of the reference distributions. When the $d_k$ differ, Assumption~\ref{assump:aggregating-discrepancy} still holds if there exists a jointly convex discrepancy~$d$ with $d\leq d_k$ for all $k\in[K]$. For example, if $d_1=\mathsf{W}_2$, $d_2=\mathsf{TV}$ and $d_3=\mathsf{W}_\infty$ are the $2$-Wasserstein, total variation and $\infty$-Wasserstein distances on a compact set~$\Xi$, then one can use the scaled $1$-Wasserstein distance $d=\mathsf{W}_1/D$ with $D\coloneqq\max\{1,\mathrm{diam}(\Xi)\}$. Indeed, $\mathsf{W}_1$ is jointly convex and satisfies the estimates $\mathsf{W}_1\leq\mathrm{diam}(\Xi)\,\mathsf{TV}$ \citep[theorem~4]{gibbs-metrics:2002} and $\mathsf{W}_1\leq\mathsf{W}_2\leq\mathsf{W}_\infty$.

Under Assumption~\ref{assump:aggregating-discrepancy}, define the mixture radius $\varepsilon\coloneqq\sum_{k=1}^K t_k\varepsilon_k$ and the mixture ambiguity set $\mathcal{P}_\mathrm{mix}\coloneqq\mathfrak{B}^d_\varepsilon(\mathbb{P}_\mathrm{mix})$, with corresponding worst-case risk
\begin{align}\tag{P$_\mathrm{mix}$}\label{prob:mix}
    v_\mathrm{mix}\coloneqq\sup_{\mathbb{Q}\in\mathcal{P}_\mathrm{mix}}\mathbb{E}_{\mathbb{Q}}[\ell(\bm{\xi})].
\end{align}
The following proposition places $v_\mathrm{mix}$ in the approximation hierarchy.

\begin{proposition}[Mixture Relaxation]\label{prop:mix}
Suppose that Assumption~\ref{assump:aggregating-discrepancy} holds. Then, for any weight vector $\bm{t}\in\Delta(K)$, we have $\mathcal{P}_\mathrm{mink}\subseteq\mathcal{P}_\mathrm{mix}$~and
\[
    v_\mathrm{P}\leq v_\mathrm{stat}\leq v_\mathrm{mink}\leq v_\mathrm{mix}.
\]
\end{proposition}

Proposition~\ref{prop:mix} shows that $v_\mathrm{mix}$ is the weakest bound in the hierarchy. Its appeal lies in reducing the $K$-reference problem to a single-reference DRO problem, at the cost of discarding all reference-specific structure by collapsing the $K$ individual ambiguity sets into one.

\subsection{Wasserstein Barycenter Approximation}\label{sec:barycenter}
Throughout this section, we assume that $\Xi$ is convex and that each elementary ambiguity set~$\mathfrak B_k$, $k\in[K]$, is a $p$-Wasserstein ball for some common fixed~$p\in[1,\infty)$. In this special case, a third aggregation strategy is to compute a $p$-Wasserstein barycenter~\citep{agueh2011barycenters} of the reference distributions and construct a single $p$-Wasserstein ball around it. A $p$-Wasserstein barycenter corresponding to a given weight vector~$\bm{w}\in\Delta(K)$ is any distribution
\begin{align}\label{eq:barycenter}
    \mathbb{P}_\mathrm{bary}\in\mathop{\argmin}_{\mathbb{Q}\in\mathcal{P}_0(\mathbb R^m)} ~ \sum_{k=1}^K w_k\mathsf{W}_p^p(\mathbb{Q},\hat{\mathbb{P}}_k).
\end{align}
Problem~\eqref{eq:barycenter} is always solvable \citep[theorem~2]{legouic2017existence} and therefore $\mathbb{P}_\mathrm{bary}$ is guaranteed to exist. As~$\Xi$ is convex, one can further show that $\mathbb{P}_\mathrm{bary}$ is supported on~$\Xi$. For a prescribed radius $\varepsilon_\mathrm{bary}\geq 0$, we introduce the barycenter ambiguity set $\mathcal{P}_\mathrm{bary}\coloneqq\{\mathbb{Q}\in\mathcal{P}_0(\Xi):\mathsf{W}_p(\mathbb{Q},\mathbb{P}_\mathrm{bary})\leq\varepsilon_\mathrm{bary}\}$ and define the corresponding worst-case risk as
\begin{align}\tag{P$_\mathrm{bary}$}\label{prob:barycenter}
    v_\mathrm{bary}\coloneqq\sup_{\mathbb{Q}\in\mathcal{P}_\mathrm{bary}}\mathbb{E}_{\mathbb{Q}}[\ell(\bm{\xi})].
\end{align}

Unlike $\mathcal{P}_\mathrm{mink}$ and $\mathcal{P}_\mathrm{mix}$, the ambiguity set $\mathcal{P}_\mathrm{bary}$ does not automatically contain $\mathcal{P}$. Guaranteeing $\mathcal{P}\subseteq\mathcal{P}_\mathrm{bary}$ requires calibrating $\varepsilon_\mathrm{bary}$ to account for the distances between the reference distributions and the barycenter. For instance, the triangle inequality gives $\mathcal{P}\subseteq\mathcal{P}_\mathrm{bary}$ whenever $\varepsilon_\mathrm{bary}\geq{\min}_{k\in[K]}\{\varepsilon_k+\mathsf{W}_p(\hat{\mathbb{P}}_k,\mathbb{P}_\mathrm{bary})\}$. However, this calibration depends on which barycenter is selected, and the barycenter need not be unique. For example, if $\Xi=[0,1]$ and there are  $K=2$ reference distributions $\hat{\mathbb{P}}_1=\delta_0$ and $\hat{\mathbb{P}}_2=\delta_1$ with equal weights, then one can show that {\em every} distribution in $\mathcal P_0(\Xi)$ constitutes a $1$-Wasserstein barycenter. Thus, infinitely many valid choices of $\mathcal{P}_\mathrm{bary}$ exist with no principled way to select among them, and neither a containment guarantee nor a meaningful comparison with $v_\mathrm{aff}$ can be ensured in general. Beyond these obstacles, the barycenter-based approach has two further shortcomings. First, computing an exact Wasserstein barycenter is NP-hard in general \citep{altschuler2022wasserstein}. The Gaussian case is a tractable exception \citep{knott1994,agueh2011barycenters}. Specifically, if $\hat{\mathbb{P}}_k=\mathcal{N}(\hat{\bm{\mu}}_k,\hat{\bm{\Sigma}}_k)$ is a Gaussian distribution with $\hat{\bm{\Sigma}}_k \succ \bm 0$ for every~$k\in[K]$ and the ground metric is the Euclidean norm on~$\mathbb R^m$, then the unique $2$-Wasserstein barycenter is also a Gaussian distribution, $\mathbb{P}_\mathrm{bary}=\mathcal{N}(\bm{\mu}_\mathrm{bary},\bm{\Sigma}_\mathrm{bary})$, whose mean vector and covariance matrix satisfy 
\begin{align}\label{eq:barycenter-mean-cov}
    \bm{\mu}_\mathrm{bary}\coloneqq\sum_{k=1}^K w_k\hat{\bm{\mu}}_k \quad \text{and} \quad  \bm{\Sigma}_\mathrm{bary}=\sum_{k=1}^K w_k\big(\bm{\Sigma}_\mathrm{bary}^{1/2}\hat{\bm{\Sigma}}_k\bm{\Sigma}_\mathrm{bary}^{1/2}\big)^{1/2}.
\end{align}
In addition, the covariance equation in~\eqref{eq:barycenter-mean-cov} has a unique positive definite solution~$\bm{\Sigma}_\mathrm{bary}$ that can be found via a fixed-point iteration \citep{alvarez2016fixed}. Second, Wasserstein barycenter-based DRO models suffer from instability, that is, their optimal values as well as their optimal solutions can change discontinuously under Wasserstein perturbations of the reference distributions \citep[example~4]{rychener2024wasserstein}. By contrast, fusion DRO is provably stable when the ambiguity set is obtained by intersecting $K$ Wasserstein balls. Specifically, one can show that the optimal value as well as the optimal solution of a fusion DRO problem change continuously under Wasserstein perturbations of the reference distributions \citep[theorem~6 and corollary~3]{rychener2024wasserstein}. In view of these limitations, we do not attempt to include the Wasserstein barycenter approach in our approximation hierarchy. Nevertheless, we use the tractable Gaussian barycenter as a benchmark in the numerical experiments of Section~\ref{sec:experiments}.

\subsection{Elastic Net Regularization}
\label{sec:elastic-net}
We close this section by connecting our approximations to elastic net regularization from supervised machine learning~\citep{zou2005regularization}. Unlike in the previous sections, all elementary ambiguity sets share the same reference distribution $\hat{\mathbb{P}}\in\mathcal{P}_0(\mathbb{R}^m)$ with $\mathbb E_{\hat{\mathbb P}}[\|\bm \xi\|_2]<\infty$. Specifically, for every $k\in[K]$, $\mathfrak{B}_k$ is the $1$-Wasserstein ball of radius $\varepsilon_k>0$ around $\hat{\mathbb{P}}$, with transport cost induced by $\|\cdot\|_{p_k}$. In addition, $\Xi=\mathbb{R}^m$ and $\ell(\bm{\xi})=L(\bm{\beta}^\top\bm{\xi})$ for some model $\bm\beta\in\mathbb R^m$, where  $L:\mathbb{R}\to\mathbb{R}$ is convex and Lipschitz continuous with modulus $\mathrm{Lip}(L)$. The fusion ambiguity set thus intersects $K$ Wasserstein balls, each encoding a different notion of perturbation of the same reference, and captures a richer geometry than any single ball unless all norms $\|\cdot\|_{p_k}$ coincide.

For a fixed weight vector $\bm{t}\in\Delta(K)$, the Minkowski sum approximation of Section~\ref{sec:minkowski} decomposes into $K$ separate worst-case expectations as shown in~\eqref{prob:mink}. Since $\Xi=\mathbb{R}^m$, each worst-case expectation evaluates to $\mathbb{E}_{\hat{\mathbb{P}}}[\ell(\bm{\xi})]+\varepsilon_k\,\mathrm{Lip}(L)\|\bm{\beta}\|_{p_k^*}$, where $p_k^*=p_k/(p_k-1)$ is the H\"older conjugate of~$p_k$ \citep[theorem~6.3]{mohajerin2018data}, giving
\begin{align}\tag{P$_\mathrm{elastic}$}\label{prob:elastic-net}
    v_\mathrm{mink}=v_\mathrm{elastic}\coloneqq\mathbb{E}_{\hat{\mathbb{P}}}[\ell(\bm{\xi})]+\mathrm{Lip}(L)\sum_{k=1}^K t_k\varepsilon_k\|\bm{\beta}\|_{p_k^*}.
\end{align}
This objective augments the expected loss under the reference distribution with a regularization term that penalizes model complexity using a weighted combination of $K$ dual norms, one per ambiguity set, in exact analogy with elastic net regularization. In particular, if~$K=2$, $p_1^*=1$, and $p_2^*=2$, then the regularization term reduces to a weighted sum of $\ell_1$ and $\ell_2$ penalties, mirroring the elastic net of \citet{zou2005regularization} (which penalizes $\|\bm{\beta}\|_2^2$ rather than $\|\bm{\beta}\|_2$).

For the mixture approximation of Section~\ref{sec:mixture}, the $K$ references coincide with~$\hat{\mathbb{P}}$, and therefore the mixture center is $\hat{\mathbb{P}}$ regardless of $\bm{t}$. The tightest admissible transport norm for $d$ in Assumption~\ref{assump:aggregating-discrepancy} satisfies $\|\cdot\|\leq\|\cdot\|_{p_k}$ for all $k$, so the tightest choice is $\|\cdot\|=\min_{k\in[K]}\|\cdot\|_{p_k}$, with dual norm $\max_{k\in[K]}\|\cdot\|_{p_k^*}$. Evaluating the worst-case expectation over this single ball gives
\begin{align*}
    v_\mathrm{mix}=\mathbb{E}_{\hat{\mathbb{P}}}[\ell(\bm{\xi})]+\mathrm{Lip}(L)\Big(\sum_{k=1}^K t_k\varepsilon_k\Big)\max_{k\in[K]}\|\bm{\beta}\|_{p_k^*}.
\end{align*}
The two benchmarks differ only in how radii are paired with norms. Indeed, $v_\mathrm{elastic}$ uses a weighted average of norm-radius products, whereas $v_\mathrm{mix}$ pairs the aggregate radius with the largest dual norm. As $\sum_{k=1}^K t_k\varepsilon_k\|\bm{\beta}\|_{p_k^*}\leq(\sum_{k=1}^K t_k\varepsilon_k)\max_{k\in[K]}\|\bm{\beta}\|_{p_k^*}$, we recover~$v_\mathrm{elastic}\leq v_\mathrm{mix}$, consistent with Proposition~\ref{prop:mix}. In the elastic net case ($K=2$, $p_1^*=1$, $p_2^*=2$), the mixture penalty collapses to the pure $\ell_1$ penalty $(t_1\varepsilon_1+t_2\varepsilon_2)\,\mathrm{Lip}(L)\|\bm{\beta}\|_1$ without $\ell_2$ component.

As $L$ is convex and Lipschitz continuous with modulus~$\mathrm{Lip}(L)$, its conjugate $L^*$ has bounded domain with closure $[\underline{\theta},\overline{\theta}]$, where $\max\{|\underline{\theta}|,|\overline{\theta}|\}=\mathrm{Lip}(L)$ \citep[corollary~13.3.3]{R97:convex_analysis}. Next, we show that the static approximation improves upon the elastic net benchmark.

\begin{proposition}[Elastic Net Approximation]\label{prop:elastic-net}
For any $\bm{t}\in\Delta(K)$, Problem~\eqref{prob:static_upper_bound} reduces to
\begin{align}\label{eq:elastic-static-value}
    v_\mathrm{stat}=\mathbb{E}_{\hat{\mathbb{P}}}[\ell(\bm{\xi})]+\inf_{\substack{\bm{a}_k\in\mathbb{R}^m,\,k\in[K]\\\sum_{k=1}^K\bm{a}_k=\bm{0}}}\sum_{k=1}^K\varepsilon_k\max\Big\{\|t_k\underline \theta\bm{\beta}+\bm{a}_k\|_{p_k^*},\|t_k\overline \theta \bm{\beta}+\bm{a}_k\|_{p_k^*}\Big\}.
\end{align}
Moreover, $\bm{a}_k=\bm{0}$ for all $k\in[K]$ is feasible in~\eqref{eq:elastic-static-value} and recovers~\eqref{prob:elastic-net}, so $v_\mathrm{P}\leq v_\mathrm{stat}\leq v_\mathrm{elastic}$.
\end{proposition}

Proposition~\ref{prop:elastic-net} shows that the static approximation improves on the elastic net benchmark. Optimizing over $\{\bm{a}_k\}_k$ subject to the zero-sum constraint can only tighten the bound. Finally, since $\ell(\bm{\xi})=L(\bm{\beta}^\top\bm{\xi})$ grows at most linearly by Lipschitz continuity of $L$, the optimal $\bm{A}_k$ in~\eqref{prob:affine} must be zero \citep[proposition~2]{gao2023distributionally}, so $v_\mathrm{aff}=v_\mathrm{stat}$ in this setting.

\section{Efficient Computation via Stochastic Gradient Descent}\label{sec:sgd}
So far, the decision $\bm{\beta}$ has been held fixed while we approximated the adversary's problem. We now return to the decision maker's problem. Writing $\ell(\bm{\beta};\bm{\xi})$ for the loss and $v_\mathrm{aff}(\bm{\beta})$ for the optimal value of~\eqref{prob:affine} with $\ell(\bm{\beta};\cdot)$ in place of $\ell$, we replace the inner worst-case expectation in~\eqref{eq:classical_dro} by its affine approximation $v_\mathrm{aff}(\bm{\beta})$. This yields a conservative surrogate of the full DRO problem, which minimizes over $\bm{\beta}$ and the decision rule coefficients jointly. We thus aim to solve
\begin{align}\label{prob:dro-affine}
    \inf_{\bm{\beta}\in\mathcal{B}} \; v_\mathrm{aff}(\bm{\beta})=\left\{\begin{array}{ll}
        \inf & \displaystyle\sum_{k=1}^K\sup_{\mathbb{Q}_k\in\mathfrak{B}_k}\mathbb{E}_{\mathbb{Q}_k}\big[t_k\ell(\bm{\beta};\bm{\xi})+\bm{a}_k^\top\bm{\xi}+\langle\bm{A}_k,\bm{\xi}\bm{\xi}^\top\rangle\big] \\
        \mathrm{s.t.} & \bm{\beta}\in\mathcal{B},\;\bm{a}_k\in\mathbb{R}^m,\;\bm{A}_k\in\mathbb{S}^m\quad\forall k\in[K] \\
        & \sum_{k=1}^K\bm{a}_k=\bm{0},\;\sum_{k=1}^K\bm{A}_k=\bm{0},
    \end{array}\right.
\end{align}
where the feasible set $\mathcal{B}$ is assumed to be convex and closed. In the remainder, we discuss solution schemes for Problem~\eqref{prob:dro-affine} when the elementary ambiguity sets $\mathfrak{B}_k$, $k\in[K]$, are Wasserstein, $\phi$-divergence or moment ambiguity sets. In all these cases, we will argue that~\eqref{prob:dro-affine} is amenable to (stochastic or deterministic) first-order methods that scale to large problem instances.

\subsection{Wasserstein Ambiguity Sets}
\label{sec:gradient-wass}
Throughout this section, each $\mathfrak{B}_k$ is a $p$-Wasserstein ball of radius $\varepsilon_k>0$ around~$\hat{\mathbb{P}}_k\in\mathcal{P}_0(\Xi)$ for some~$p\in[1,\infty)$ and some norm~$\|\cdot\|$ on~$\mathbb{R}^m$ (see Example~\ref{ex:wasserstein-ambiguity-set}), and the loss~$\ell$ is bounded below on~$\mathcal{B}\times\Xi$. Together with Assumption~\ref{assump:general}, boundedness from below makes $\ell(\bm{\beta};\bm{\xi})$ bounded in $\bm{\xi}$ for every $\bm{\beta}\in\mathcal{B}$. By standard duality results for Wasserstein DRO \citep[theorem~4.18]{kuhn2025distributionally}, all embedded maximization problems in the objective of~\eqref{prob:dro-affine} can thus be replaced by equivalent minimization problems over Lagrange multipliers~$\gamma_k\geq 0$, $k\in[K]$. Hence, \eqref{prob:dro-affine} can be recast as
\begin{align}\label{prob:wass-dual}
    \min_{\bm{\theta}\in \Theta} \left\{ F(\bm{\theta})\coloneqq\sum_{k=1}^K\mathbb{E}_{\hat{\mathbb{P}}_k}\Big[\max_{\bm{\xi}\in\Xi}f_k(\bm{\theta};\bm{\xi},\hat{\bm{\xi}}_k)\Big] \right\} ,
\end{align}
where $\bm{\theta}\coloneqq(\bm{\beta},\{\bm{a}_k,\bm{A}_k,\gamma_k\}_k)$ collects the original decisions, the decision rule coefficients and the Lagrange multipliers. To keep notation simple, we also introduce the effective loss functions
\begin{align*}
    f_k(\bm{\theta};\bm{\xi},\hat{\bm{\xi}})\coloneqq t_k\ell(\bm{\beta};\bm{\xi})+\bm{a}_k^\top\bm{\xi}+\langle\bm{A}_k,\bm{\xi}\bm{\xi}^\top\rangle-\gamma_k\|\bm{\xi}-\hat{\bm{\xi}}\|^p+\gamma_k\varepsilon_k^p, \quad k\in[K],
\end{align*}
set $\Theta_0= \mathbb R^n \times (\mathbb R^m \times \mathbb S^m\times \mathbb R)^K$, and define the effective feasible set 
\[
	\Theta \coloneqq \left\{ \bm{\theta} \in \Theta_0 : \bm \beta\in\mathcal B,\; \sum_{k=1}^K\bm{a}_k=\bm{0},\;\sum_{k=1}^K\bm{A}_k=\bm{0} ,\; \gamma_k\geq 0\;\forall k\in[K]\right\}.
\]
The maximum in~\eqref{prob:wass-dual} is attained because $f_k$ is upper semicontinuous in $\bm{\xi}$ on the compact set $\Xi$.

When~$\hat{\mathbb{P}}_k$ is a discrete distribution supported on~$N_k$ points for every~$k\in[K]$, the loss~$\ell(\bm{\beta};\bm{\xi})$ is convex in~$\bm{\beta}$ and piecewise concave in~$\bm{\xi}$, the support set     $\Xi$ is convex, and we focus on the static approximation by setting $\bm A_k=\bm 0$ for every $k\in[K]$, then Problem~\eqref{prob:wass-dual} can be reformulated as a finite convex program by appealing to standard robust optimization techniques~\citep[theorem~7.20]{kuhn2025distributionally}. The resulting program has $\mathcal{O}(\sum_{k=1}^K N_k)$ variables and constraints, in stark contrast with exact reformulations of the DRO~\eqref{eq:classical_dro} over the original fusion ambiguity set, which require $\mathcal{O}(\prod_{k=1}^K N_k)$ variables and constraints (\textit{cf}.~Section~\ref{sec:intro}). The affine approximation does not inherit this reduction, because the term $\langle\bm{A}_k,\bm{\xi}\bm{\xi}^\top\rangle$ is generically an indefinite quadratic in $\bm{\xi}$.

Outside these cases, that is, when the reference distributions are not discrete, when the loss lacks the required structure, or when the affine approximation is used, we solve~\eqref{prob:wass-dual} by a stochastic subgradient method that requires only samples from the $\hat{\mathbb{P}}_k$ and a routine for the inner maximization. To cover convex and smooth non-convex losses alike, we use the following notion of generalized convexity. A function $h(\bm \theta)$ is \emph{$\rho$-weakly convex} for some $\rho\geq 0$ if $h(\bm \theta)+\frac{\rho}{2}\|\bm \theta\|_2^2$ is convex. Thus, convex functions are $0$-weakly convex, and continuously differentiable functions with $\rho$-Lipschitz gradient are $\rho$-weakly convex. For such $h$, the subdifferential is defined as
\[
	\partial h(\bm{\theta})\coloneqq \left\{\bm{g}:h(\bm{\theta}')\geq h(\bm{\theta})+\langle\bm{g},\bm{\theta}'-\bm{\theta}\rangle-\frac{\rho}{2}\|\bm{\theta}'-\bm{\theta}\|_2^2 \quad \forall\bm{\theta}' \in \Theta_0 \right\}, 
\]
which reduces to the convex subdifferential when $\rho=0$ and to the gradient when $h$ is differentiable~\citep[section~2]{davis2019stochastic}. To define the stochastic subgradient oracle, let $\bm{\xi}_k^\star (\bm{\theta},\hat{\bm{\xi}}_k) \in\argmax_{\bm{\xi}\in\Xi}f_k(\bm{\theta};\bm{\xi},\hat{\bm{\xi}}_k)$ be a maximizer chosen as a Borel measurable function, let $\bm g_k(\bm{\theta};\bm{\xi},\hat{\bm{\xi}})\in\partial_{\bm{\theta}}f_k(\bm{\theta};\bm{\xi},\hat{\bm{\xi}})$ be a Borel measurable selection of subgradients, and set
 \begin{align}\label{eq:stoch-subgrad}
    \bm G(\bm{\theta};\hat{\bm{\xi}})\coloneqq\sum_{k=1}^K \bm g_k \left(\bm{\theta};\bm{\xi}_k^\star (\bm{\theta},\hat{\bm{\xi}}_k) ,\hat{\bm{\xi}}_k \right), \quad\hat{\bm{\xi}}=(\hat{\bm{\xi}}_1,\ldots,\hat{\bm{\xi}}_K)\sim\hat{\mathbb P} \coloneqq \hat{\mathbb{P}}_1\otimes\cdots\otimes\hat{\mathbb{P}}_K.
\end{align}
The existence of these measurable selections is part of the following result.

\begin{proposition}[Stochastic Subgradient Oracle]
\label{prop:sgd_wass}
Suppose that Assumption~\ref{assump:general} holds, that the loss function~$\ell$ is bounded below on $\mathcal{B}\times\Xi$, and that, for every $\bm{\xi}\in\Xi$, $\ell(\bm{\beta};\bm{\xi})$ is $\rho$-weakly convex and $L$-Lipschitz continuous in $\bm{\beta}$ on $\mathbb R^n$. Then $F$ is finite and $\rho$-weakly convex on $\Theta_0$, the oracle~\eqref{eq:stoch-subgrad} is well-defined and satisfies $\mathbb{E}_{\hat{\mathbb P}} [\bm G(\bm{\theta};\hat{\bm{\xi}})]\in\partial F(\bm{\theta})$ for every $\bm{\theta}\in \Theta$, and $\sup_{\bm{\theta}\in \Theta}\mathbb{E}_{\hat{\mathbb P}} \|\bm G(\bm{\theta};\hat{\bm{\xi}})\|_2^2<\infty$.
\end{proposition}

Proposition~\ref{prop:sgd_wass} supplies exactly the ingredients required by the projected stochastic subgradient method $\bm{\theta}_{t+1}=\Pi_{\Theta}(\bm{\theta}_t-\alpha_t \bm G(\bm{\theta}_t;\hat{\bm{\xi}}_t))$ with step sizes $\alpha_t>0$ and independent samples $\hat{\bm{\xi}}_t$, $t\in\mathbb Z_+$, where $\Pi_{\Theta}$ is the Euclidean projection onto the convex closed set $\Theta$, which projects $\bm{\beta}$ onto $\mathcal{B}$, clips each $\gamma_k$ at zero, and subtracts the averages $\frac{1}{K}\sum_{k=1}^K\bm{a}_k$ and $\frac{1}{K}\sum_{k=1}^K\bm{A}_k$ from the decision rule coefficients. As $F$ is shorthand for the objective of~\eqref{prob:wass-dual} and thus of~\eqref{prob:dro-affine}, we have
\[
	\inf_{\bm \theta\in\Theta}F(\bm{\theta}) \geq \inf_{\bm \beta\in\mathcal B} \mathbb{E}_{\mathbb{Q}}[\ell(\bm{\beta};\bm{\xi})]\geq\inf_{\bm\beta\in\mathcal{B}} \inf_{\bm\xi\in\Xi}\ell(\bm \beta;\bm \xi) >-\infty 
\]
for any $\mathbb{Q}\in\cap_{k=1}^K\mathfrak{B}_k$, which exists by Assumption~\ref{assump:general}\,\ref{assump:general-nonempty-intersection}. Here, the first inequality follows from the zero-sum constraints on the decision rule coefficients and the assumption that~$\bm{t}\in\Delta(K)$. This uniform lower bound on $F$ lets us invoke the convergence guarantees of \citet{davis2019stochastic} for the stochastic subgradient algorithm. Progress of this algorithm is measured through the Moreau envelope $F_\lambda(\bm{\theta})\coloneqq\min_{\bm \theta'\in \Theta}\{F(\bm \theta')+\frac{1}{2\lambda}\|\bm \theta'-\bm{\theta}\|_2^2\}$ with $\lambda\in(0,1/\rho)$, whose minimizer $\mathrm{prox}_{\lambda F}(\bm{\theta})$ is unique by strong convexity. One can show that $F_\lambda(\bm\theta)$ is continuously differentiable and that any point with $\|\nabla F_\lambda(\bm{\theta})\|_2\leq\delta$ lies within $\lambda\delta$ of a $\delta$-stationary point of~\eqref{prob:wass-dual} \citep[section~2]{davis2019stochastic}. With $\alpha_t\propto 1/\sqrt{T}$, $\lambda=1/(2\rho)$, and an output drawn uniformly at random from the first $T$ iterates, the algorithm returns $\bm{\theta}$ with $\mathbb{E}\|\nabla F_\lambda(\bm{\theta})\|_2\leq\delta$ after $T=\mathcal{O}(\delta^{-4})$ stochastic subgradient evaluations~\citep[theorem~3.1]{davis2019stochastic}. If $\ell(\bm{\beta};\bm{\xi})$ is convex in $\bm{\beta}$ for every $\bm{\xi}$, then $\rho=0$, $F$ is convex, and the averaged iterate is $\delta$-optimal in expectation after~$\mathcal{O}(\delta^{-2})$ evaluations~\citep[theorem~4.1]{davis2019stochastic}.

\begin{remark}[Inner Maximization]
Each stochastic subgradient evaluation requires the global maximization of $f_k(\bm{\theta};\bm{\xi},\hat{\bm{\xi}}_k)$ over $\bm{\xi}\in\Xi$ for every $k\in[K]$. These problems are non-convex in general because $\langle\bm{A}_k,\bm{\xi}\bm{\xi}^\top\rangle$ is an indefinite quadratic and $\ell(\bm{\beta};\bm{\xi})$ is merely upper semicontinuous in~$\bm{\xi}$. They are tractable when $\Xi$ is finite. In addition, they are also tractable when $p=2$, $\|\cdot\|=\|\cdot\|_2$ is the Euclidean norm, $\ell(\bm{\beta};\bm{\xi})$ is quadratic in $\bm{\xi}$, and $\Xi$ is an ellipsoid, in which case the maximization over $\bm \xi$ is solved exactly by a semidefinite program via the S-lemma. For the static approximation ($\bm{A}_k=\bm{0}$) the quadratic term disappears, and the inner problem is concave whenever $\ell(\bm{\beta};\bm{\xi})$ is concave in~$\bm{\xi}$. An inexact inner maximization introduces a bias into $\bm{G}$ that is controlled by the suboptimality of the maximizer. We do not pursue this refinement here.
\end{remark} 

\subsection{Other Ambiguity Sets}
\label{sec:gradient-other}
The stochastic subgradient descent methods of Section~\ref{sec:gradient-wass} readily extend from intersections of Wasserstein balls to other fusion ambiguity sets. The separation of the~$K$ inner suprema in the objective function of Problem~\eqref{prob:dro-affine} ensures that the stochastic subgradient in~$\bm{\theta}$ is always representable as a sum of~$K$ per-ball contributions, which can be evaluated in parallel. The computation of each contribution depends on the structure of the underlying ambiguity set~$\mathfrak{B}_k$. For Wasserstein and $\phi$-divergence ambiguity sets, it is obtained from stochastic oracles, and for moment ambiguity sets it can be computed deterministically and in closed form. Crucially, however, the per-ball terms only interact through the zero-sum constraints on the decision rule coefficients. Intersections of ambiguity sets of different types therefore require no additional machinery beyond combining the respective routines. What changes from one type to another is the nature of the per-ball supremum problem, the regularity needed from the loss, and whether discrete reference distributions make the intersection $\cap_{k=1}^K\mathfrak B_k$ meaningful at all.

For $\phi$-divergence balls (see Example~\ref{ex:phi-divergence-ambiguity-set}) the stochastic subgradient oracles are cheaper than for Wasserstein balls. The dual of each inner supremum minimizes the expectation under the corresponding reference distribution of a composition of the conjugate entropy function~$\phi^*$ with the loss, and it does {\em not} involve a global maximization over~$\Xi$ \citep[theorem~4.15]{kuhn2025distributionally}. If the loss is convex in~$\bm{\beta}$, then the dual objectives are convex in the decision, the decision rule coefficients and two scalar multipliers per ball, and differentiating under the expectation yields an unbiased estimate of a subgradient. 
Two caveats apply. First, each entropy function must grow superlinearly, so that $\phi^*$ is finite everywhere. Otherwise the dual problems retain a semi-infinite constraint involving the loss function (which may be nonconcave in~$\xi$). Second, convergence of the stochastic subgradient method is guaranteed only if the iterates remain in a compact region on which the multipliers of the divergence constraints are bounded away from zero \citep{namkoong2016stochastic}, because each of these multipliers enters a denominator in the oracle and the variance of the subgradient estimate explodes as it approaches zero. Unlike in Section~\ref{sec:gradient-wass}, where compactness of~$\Xi$ already yields a uniform second-moment bound, the decision, the decision rule coefficients and the multipliers must therefore be confined to a bounded set. Such a restriction shrinks the feasible set of the dual minimization problem and hence preserves the validity of the bound $v_\mathrm{P}\leq v_\mathrm{aff}$, at the price of loosening it if the region is chosen too small. When the reference distribution is discrete, then the expectation in the dual objective reduces to a finite sum, in which case one can even compute deterministic subgradients. However, under our standing assumption that $\phi$ grows superlinearly, the intersection of multiple $\phi$-divergence balls is empty unless the supports of the reference distributions overlap. Thus, the intersection is empty almost surely when the references are empirical distributions of independent samples from continuous populations. Intersections of $\phi$-divergence balls are therefore meaningful mainly for continuous references with a common support, such as fitted parametric models or kernel density estimators, which is precisely the regime that calls for the stochastic methods described above.

An intersection of moment ambiguity sets (see Example~\ref{ex:moment-ambiguity-sets}) is again a moment ambiguity set that jointly enforces all moment conditions of the elementary ambiguity sets. While intersections of Chebyshev ambiguity sets are generically empty, intersections of sufficiently large Gelbrich balls are nondegenerate, and worst-case expectations over such intersections admit semidefinite reformulations reminiscent of those for single Gelbrich balls \citep{nguyen2021mean}. For this reason, deterministic subgradients of the dual objective can be computed efficiently. As the expectation of the correction terms $\bm{a}_k^\top\bm{\xi}+\langle\bm{A}_k,\bm{\xi}\bm{\xi}^\top\rangle$ depends on a distribution only through its mean and second-moment matrix, one can show that the affine decision rule approximation is exact if~$\mathfrak B_k$ is a moment ambiguity set and the corresponding uncertainty set $\mathcal{U}_k$ is convex in the mean and the second-moment matrix for every $k\in[K]$. This is the case for Gelbrich balls, for example.

\section{Numerical Experiments}
\label{sec:experiments}
When several sources of information bear on a decision, should a decision maker intersect the elementary ambiguity sets they induce or aggregate the sources into a single ambiguity set? Our theory shows that the decision-rule approximations of the intersection yield tighter worst-case bounds than the Minkowski and mixture relaxations of Section~\ref{sec:comparison}. It leaves open, however, whether tighter bounds translate into better decisions out of sample, how much the affine decision rules improve over the static ones, and what this costs computationally. We answer these questions for a classical portfolio optimization problem, a multi-item newsvendor problem and a neural network classification problem, by comparing the affine approximation of Proposition~\ref{prop:affine_upper_bound} (\textsc{affine}) and the static approximation of Corollary~\ref{cor:static_upper_bound} (\textsc{static}) against standard benchmarks.

\paragraph{Data.} In the experiments of Sections~\ref{sec:portfolio} and~\ref{sec:newsvendor}, the decision maker observes data from $K$ sources. Source $k\in[K]$ provides $N_k$ samples, and the reference distribution $\hat{\mathbb{P}}_k$ is the empirical distribution of these samples. Pooling the samples of all sources yields the pooled empirical distribution, which is the mixture of the reference distributions in the sense of Section~\ref{sec:mixture} with weights $t_k=N_k/\sum_{j=1}^KN_j$. In the experiment of Section~\ref{sec:SGD-experiment}, by contrast, the decision maker observes only noisy summary statistics of a hidden dataset.

\paragraph{Methods.} Table~\ref{tab:baselines} summarizes the methods and the experiments in which they are used. Our \textsc{affine} and \textsc{static} approximations intersect the~$K$ elementary ambiguity sets and restrict the Lagrange multipliers to decision rules as in Section~\ref{sec:main}. They apply to any elementary ambiguity sets. The benchmarks instead aggregate the information sources into a single ambiguity set, as discussed in Section~\ref{sec:comparison}. Most of them aggregate the reference distributions and thus require samples. Hence, they are used only in the experiments of Sections~\ref{sec:portfolio} and~\ref{sec:newsvendor}. \textsc{pool} centers a single ambiguity set of the same type as the elementary ones at the pooled empirical distribution, and \textsc{erm} is its special case with radius zero. \textsc{mixture} centers such a set at the mixture of the reference distributions with the weights described below (Section~\ref{sec:mixture}). Since exact Wasserstein barycenters are hard to compute (Section~\ref{sec:barycenter}), \textsc{barycenter} uses a Gaussian proxy. We evaluate the mean vector and covariance matrix of each $\hat{\mathbb{P}}_k$, aggregate them via~\eqref{eq:barycenter-mean-cov} into a Gaussian Gelbrich barycenter $\mathbb{P}_\mathrm{bary}\coloneqq\mathcal{N}(\bm{\mu}_\mathrm{bary},\bm{\Sigma}_\mathrm{bary})$, and center a single ambiguity set at the empirical distribution $\hat{\mathbb{P}}_\mathrm{bary}$ of $N_\mathrm{bary}$ independent samples from~$\mathbb{P}_\mathrm{bary}$. \textsc{rkhs}~\citep{lixianyu2026robust} represents each source by the kernel mean embedding of its reference distribution, learns a support measure clustering ball around these embeddings~\citep{muandet2013one} with conservatism controlled by a regularization parameter~$\nu\in(0,1)$, and solves the robust problem over the induced ambiguity set. In contrast, \textsc{minkowski} aggregates the elementary ambiguity sets themselves through their Minkowski sum (Section~\ref{sec:minkowski}) and thus applies in all three experiments. Further references that are specific to individual experiments are introduced where they are used.

\paragraph{Ambiguity sets and solution methods.} In Sections~\ref{sec:portfolio} and~\ref{sec:newsvendor}, each elementary ambiguity set $\mathfrak{B}_k$ is a $2$-Wasserstein ball around $\hat{\mathbb{P}}_k$ with the Euclidean norm (Example~\ref{ex:wasserstein-ambiguity-set}). This geometry suits our quadratic decision rules because it underpins the exactness result of Proposition~\ref{prop:affine-tightness} and aligns with the Gelbrich geometry underlying \textsc{barycenter}. In both experiments, \textsc{affine} reduces to a semidefinite program, \textsc{erm} to a linear program and all other methods except \textsc{rkhs} to second-order cone programs, which we solve with MOSEK, while \textsc{rkhs} is solved with the procedure of~\citet{lixianyu2026robust}. Section~\ref{sec:SGD-experiment}, by contrast, intersects a moment ambiguity set and a histogram ambiguity set in a classification problem involving a neural network. As the loss is nonconvex in the decisions, we address this problem with the subgradient method of Section~\ref{sec:sgd}.

\begin{table}[t]
\centering
\small
\setlength{\tabcolsep}{5pt}
\caption{Overview of methods. Check marks indicate the sections in which each method is used.}
\label{tab:baselines}
\begin{tabular}{@{}llccc@{}}
\toprule
Method & Ambiguity set & \S\ref{sec:portfolio} & \S\ref{sec:newsvendor} & \S\ref{sec:SGD-experiment} \\
\midrule
\textsc{affine} & intersection, affine rules & \checkmark & \checkmark & \checkmark \\
\textsc{static} & intersection, static rules & \checkmark & \checkmark & \checkmark \\
\textsc{minkowski} & Minkowski sum & \checkmark & \checkmark & \checkmark \\
\textsc{mixture} & set around mixture & \checkmark & \checkmark & \\
\textsc{pool} & set around pooled distribution & & \checkmark & \\
\textsc{erm} & pooled distribution & & \checkmark & \\
\textsc{barycenter} & set around Gaussian barycenter & & \checkmark & \\
\textsc{rkhs} & ball of kernel mean embeddings & & \checkmark & \\
\midrule
\textsc{exact} & intersection, exact & \checkmark & & \\
\textsc{moment}, \textsc{tv} & single elementary set & & & \checkmark \\
\textsc{oracle} & hidden empirical distribution & & & \checkmark \\
\bottomrule
\end{tabular}
\end{table}

\paragraph{Calibration and implementation.} In Sections~\ref{sec:portfolio} and~\ref{sec:newsvendor}, the radii of the elementary ambiguity sets of \textsc{affine} and \textsc{static} are set to $\varepsilon_k=(1+\delta)\,\mathsf{W}_2(\hat{\mathbb{P}}_k,\hat{\mathbb{P}}_\mathrm{bary})$ for a parameter $\delta\geq 0$. This parametrization ensures that $\cap_{k=1}^K\mathfrak{B}_k$ contains $\hat{\mathbb{P}}_\mathrm{bary}$ and is thus nonempty. It also anchors the fusion ambiguity set at the same proxy as \textsc{barycenter}, so that differences in performance reflect how the sources are combined rather than where the ambiguity sets are centered. The parameter~$\delta$, the hyperparameters of the benchmarks and the weights~$t_k$ are chosen differently in each experiment, as described in the respective sections. In Section~\ref{sec:SGD-experiment}, all radii are instead calibrated from the known noise distributions of the summary statistics. The experiments of Sections~\ref{sec:portfolio} and~\ref{sec:newsvendor} were run on two 18-core Intel Xeon Gold 6240 CPUs at 2.60\,GHz with 376\,GiB of RAM, and those of Section~\ref{sec:SGD-experiment} on NVIDIA GH200 GPUs with 96\,GB of memory. Code and data are available at 
\url{https://github.com/busesen/fusion-DRO}.

\subsection{Portfolio Optimization}\label{sec:portfolio}
The first experiment assesses the tightness and computational cost of different approximations of the exact fusion DRO problem. In a portfolio problem small enough for the exact reformulation to be solved, we find that \textsc{affine} stays significantly closer to the exact value than the aggregation-based bounds, that this advantage grows with the number of sources, and that \textsc{affine} remains tractable long after the exact reformulation has run out of memory.

\paragraph{Model.} An investor allocates one unit of capital across $m=10$ assets with uncertain returns~$\bm{\xi}\in\mathbb{R}^m$. As short sales are not allowed, the portfolio weights $\bm{x}$ range over the probability simplex $\Delta(m)\subseteq\mathbb{R}^m$. The investor seeks to minimize the mean of the portfolio loss $-\bm{x}^\top\bm{\xi}$ plus~$\kappa=10$ times the conditional value-at-risk (CVaR) of this loss at level $\eta = 0.2$, in the worst case over all distributions in the fusion ambiguity set. Consequently, the investor aims to solve
\[
    \inf_{\bm{x}\in\Delta(m)}\;\sup_{\mathbb{Q}\in\cap_{k=1}^K\mathfrak{B}_k}\Big\{\mathbb{E}_{\mathbb{Q}}[-\bm{x}^\top\bm{\xi}]+\kappa\cdot\mathrm{CVaR}_\eta^{\mathbb{Q}}(-\bm{x}^\top\bm{\xi})\Big\}.
\]
By using the variational formula for the CVaR \citep{rockafellar2000optimization}, one can express the objective of this minimax problem as an infimum over an auxiliary variable~$\tau\in\mathbb{R}$. This infimum can be interchanged with the supremum over $\mathbb{Q}$ by the lopsided minimax theorem for optimized certainty equivalents~\citep[theorem~5.18]{kuhn2025distributionally}, which applies because $\kappa>0$, $\eta<1$, and all distributions in the fusion ambiguity set have finite first moments. Proceeding as in \citep[\S~7.1]{mohajerin2018data}, one can then recast the investor's problem as an instance of~\eqref{eq:classical_dro}, where $\bm{\beta}\coloneqq(\bm{x},\tau)\in\mathbb R^n$ for $n=m+1$, $\mathcal{B}\coloneqq\Delta(m)\times\mathbb{R}$, $\mathcal P\coloneqq\cap_{k=1}^K \mathfrak B_k$, and 
\begin{equation}
	\label{eq:pw-loss}
	\ell(\bm{\beta};\bm{\xi})\coloneqq\max_{i\in\{1,2\}}\{b_i\bm{x}^\top\bm{\xi}+c_i\tau\}
\end{equation}
with $b_1=-1$, $b_2=-1-\kappa/\eta$, $c_1=\kappa$ and $c_2=\kappa(1-1/\eta)$. As $\mathfrak B_k$ is a $2$-Wasserstein ball around a discrete distribution with $N_k$ atoms for every $k\in[K]$, this instance of~\eqref{eq:classical_dro} has an exact second-order cone reformulation with $\mathcal O(\prod_{k=1}^KN_k)$ variables and constraints \citep[corollary~2]{rychener2024wasserstein}, which we refer to as \textsc{exact}. It serves as the benchmark for the approximations. 

\paragraph{Data.} We assume that each source $k$ provides $N_k$ independent samples from $\mathbb{P}_k=\mathcal{N}(\bm{\mu}_k,\bm{\Sigma}_k)$ with covariance matrix $\bm{\Sigma}_k=\sigma_\mathrm{sys}^2\bm{u}_k\bm{u}_k^\top+\sigma_\mathrm{id}^2\bm{I}_m$, where $\bm{u}_k$ is drawn uniformly from the unit sphere in $\mathbb R^m$, while $\sigma_\mathrm{sys}=0.02$ and $\sigma_\mathrm{id}=0.015$ capture source-specific co-movement and idiosyncratic noise, respectively. Following~\citet[\S~6.1]{rychener2024wasserstein}, we further assume that source $k$ assigns a high mean return to the assets in a block $G_k\subseteq[m]$ of four assets and a low mean return to all other assets, that is, $(\bm{\mu}_k)_j=\mu_\mathrm{high}+(k-1)\gamma$ for $j\in G_k$ and $(\bm{\mu}_k)_j=\mu_\mathrm{low}-(k-1)\gamma$ for $j\notin G_k$, where $\mu_\mathrm{high}=0.015$, $\mu_\mathrm{low}=-0.005$ and $\gamma=0.001$. The blocks are constructed recursively. We set $G_1=\{1,2,3,4\}$, and for $k\geq 2$, we choose $G_k$ as a four-element subset of $[m]$  that differs from $G_1,\ldots,G_{k-1}$ and shares exactly two assets with $G_{k-1}$, using a deterministic rule detailed in our code, so that some assets are favored by more sources than others, and the gap between high and low mean returns widens with~$k$. Consecutive sources thus agree on two of their four high-return assets and disagree on the other two.

\paragraph{Experimental setup.} We compare \textsc{affine}, \textsc{static}, \textsc{minkowski} and \textsc{mixture} with \textsc{exact}. The other methods of Table~\ref{tab:baselines} are omitted because they are not guaranteed to upper-bound \textsc{exact}. To isolate tightness and cost from calibration, we fix weights $t_k=1/K$ and radii $\varepsilon_k=(1+\delta)\,\mathsf{W}_2(\hat{\mathbb{P}}_k,\hat{\mathbb{P}}_\mathrm{bary})$ with $\delta=0.01$ for all $k\in[K]$. \textsc{exact}, \textsc{affine}, \textsc{static} and \textsc{minkowski} use these radii for the individual balls, and \textsc{mixture} uses the radius $(\sum_{k=1}^Kt_k\varepsilon_k^2)^{1/2}$, for which Proposition~\ref{prop:mix} applies to $2$-Wasserstein balls. We vary $K\in\{2,3,\ldots,20\}$ with fixed sample size $N_k=10$, and $N_k\in\{5,10,20,\ldots,100\}$ with fixed number of sources $K=3$. In the sweep over~$K$, $\hat{\mathbb{P}}_\mathrm{bary}$ consists of $N_\mathrm{bary}=10$ samples, is computed once from all $20$ sources and is held fixed. Adding a source then adds a ball to the intersection. Hence, the fusion ambiguity sets of a replication are nested and shrink as $K$ grows. All results are averaged over $10$ replications. 

\begin{figure}
    \centering
    \includegraphics[width=0.9\linewidth]{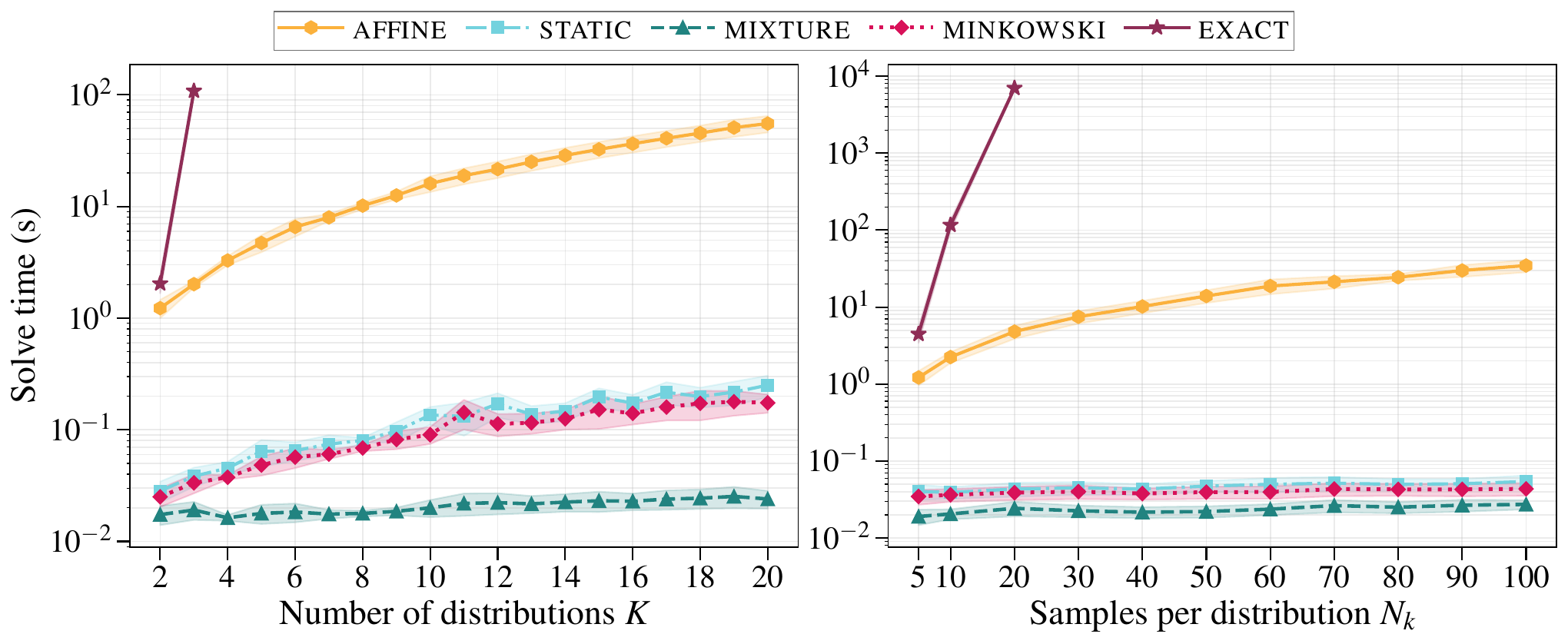}
    \caption{Mean solve time (logarithmic scale) with 95\% confidence intervals as a function of the number of sources $K$ for $N_k=10$ (left) and of the sample size $N_k$ for $K=3$ (right).}
    \label{fig:scaling}
\end{figure}

\begin{figure}
    \centering
    \includegraphics[width=0.9\linewidth]{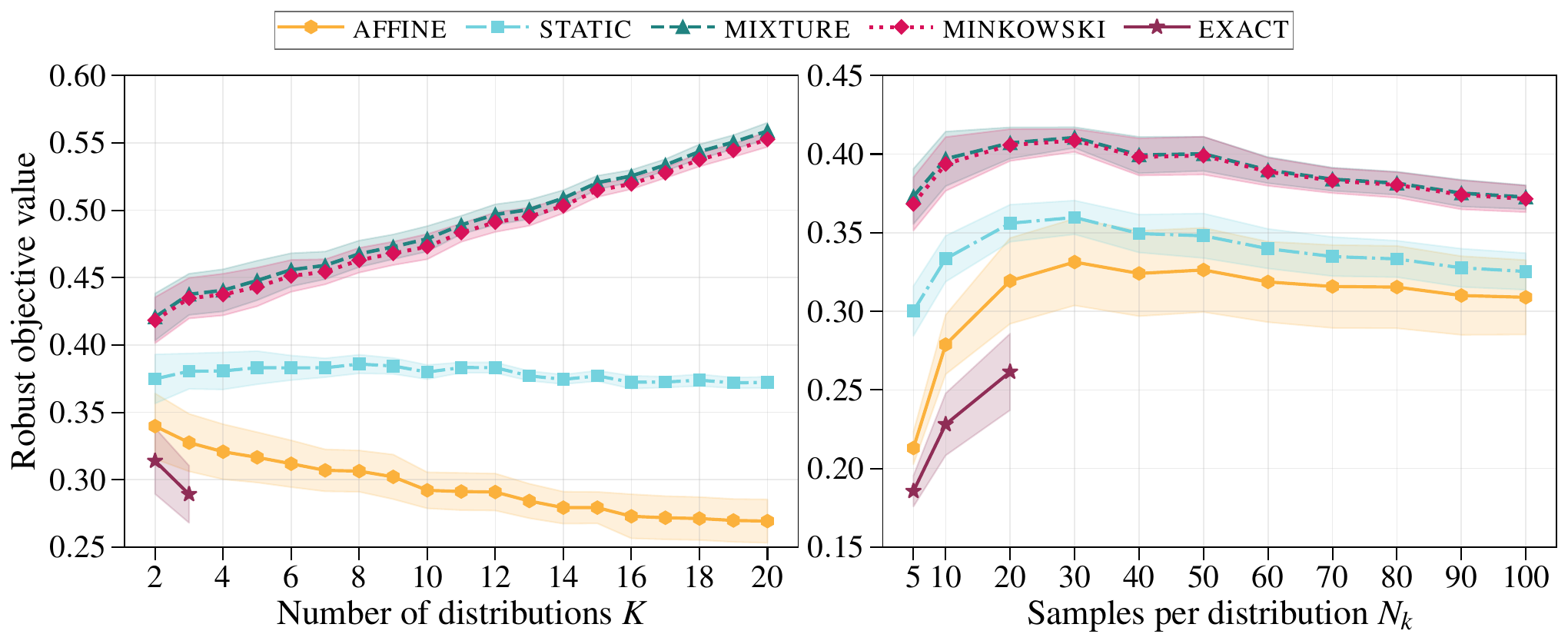}
    \caption{Mean in-sample robust objective value with 95\% confidence intervals as a function of the number of sources $K$ for $N_k=10$ (left) and of the sample size $N_k$ for $K=3$ (right).}
    \label{fig:objective}
\end{figure}

\paragraph{Results.} Figure~\ref{fig:scaling} confirms that the exact reformulation is impractical beyond toy sizes. Its solve time explodes with both $K$ and $N_k$, and larger instances exhaust the available memory, as anticipated by its $\mathcal O(\prod_{k=1}^KN_k)$ scaling. \textsc{mixture}, \textsc{minkowski} and \textsc{static} remain fast throughout, with \textsc{mixture} essentially insensitive to $K$ because it always solves a single-ball problem. \textsc{affine} gives rise to a semidefinite program and is thus slower than these second-order cone surrogates but scales more gracefully than \textsc{exact} over the entire range. Figure~\ref{fig:objective} shows the in-sample objective values. They are ordered exactly as the theory predicts, \textsc{affine}$\,\leq\,$\textsc{static}$\,\leq\,$\textsc{minkowski}$\,\leq\,$\textsc{mixture}. Where \textsc{exact} can be computed,\footnote{For \textsc{exact} with $N_k=20$, MOSEK terminated without certifying optimality in $5$ of the $10$ replications. In these cases, we use the cost of the last iterate instead of the exact minimum.} \textsc{affine} stays within $14\%$ of it in the sweep over $K$ and within $23\%$ in the sweep over $N_k$, whereas \textsc{minkowski} and \textsc{mixture} overshoot it by $33$--$51\%$ and by up to $101\%$ for small $N_k$, respectively. In addition, the gaps widen with~$K$. At $K=20$, \textsc{affine} lies $28\%$ below \textsc{static} and about $52\%$ below the two aggregation-based bounds. This reflects a structural difference. Every additional source shrinks the fusion ambiguity set, and the decision rule approximations, which retain the coupling between the balls, translate this into less conservative bounds. The aggregation-based bounds average over the individual balls instead. Consequently, additional sources cannot make them less conservative in the same way. \textsc{minkowski} and \textsc{mixture} are nearly indistinguishable here, presumably because the radii $\varepsilon_k$ are of similar size, so that the mixture radius $(\sum_{k=1}^Kt_k\varepsilon_k^2)^{1/2}$ is close to each of them.

\begin{figure}[t]
    \centering
    \includegraphics[width=\textwidth]{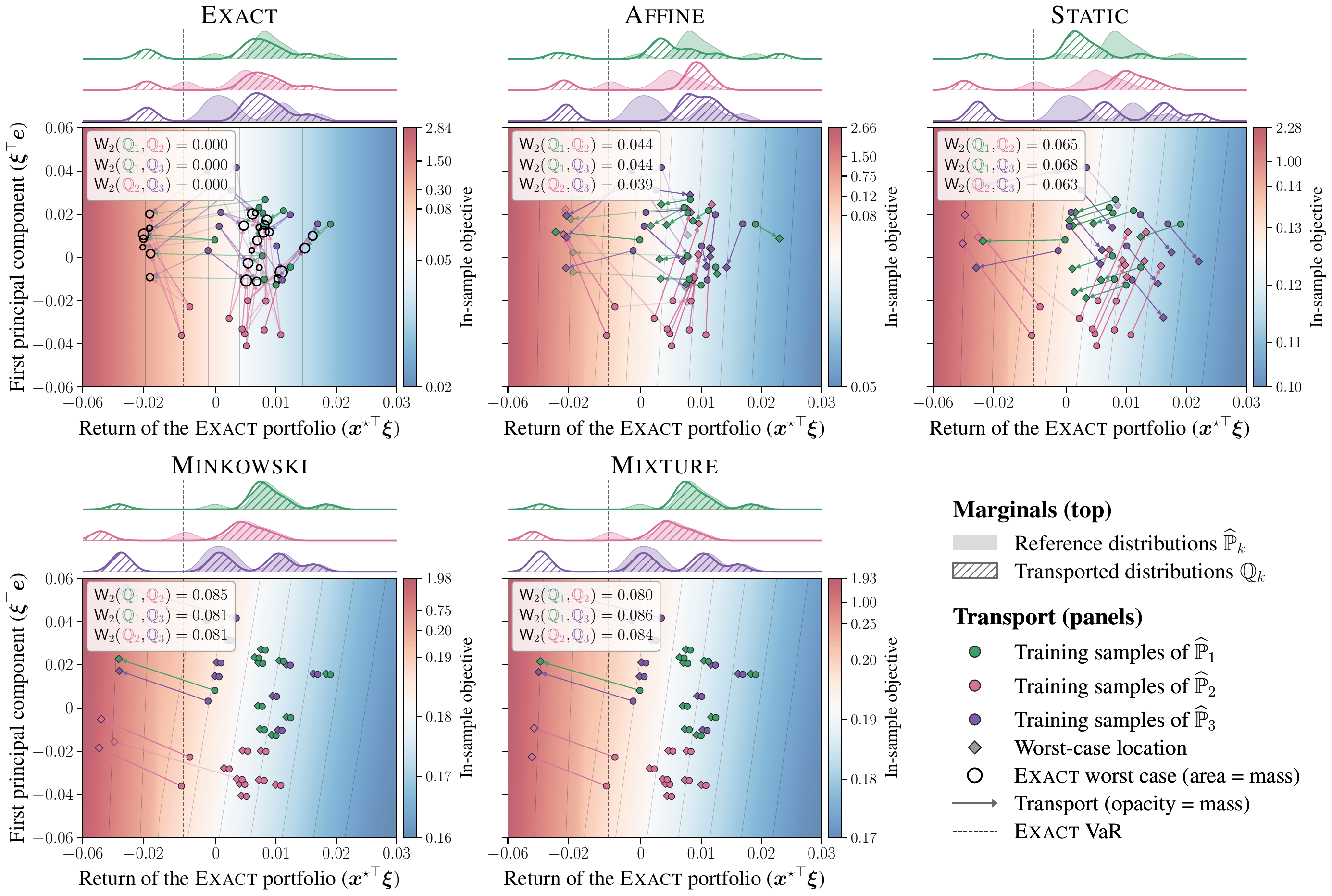}
    \caption{Worst-case transport of the training samples for one instance with $K=3$ and $N_k=10$. The background shows the loss~\eqref{eq:pw-loss} at the decision $\bm{x}$ of the method in each panel, with a separate color scale per panel. The dashed line marks the value-at-risk of \textsc{exact}, and the curves above each panel show the marginal distributions along the horizontal axis.}
    \label{fig:worstcase}
\end{figure}

\paragraph{Anatomy of the worst-case distribution.} Figure~\ref{fig:worstcase} explains why the intersection-based bounds are tighter. For one instance with $K=3$ sources, it shows how the worst-case distribution of each method transports the training samples of each source. The horizontal axis is the return $\bm{\xi}^\top\bm{x}^\star$ of the optimal portfolio $\bm{x}^\star$ of \textsc{exact}, so that points further left incur higher losses. The vertical axis is the coordinate of $\bm{\xi}$ along the direction in which the source means differ most, that is, the leading principal component of $\bm{\mu}_1,\ldots,\bm{\mu}_K$, after removing its component along $\bm{x}^\star$. \textsc{exact} must find a single distribution that lies in every ball, so the samples of all three sources are transported to the same distribution, here concentrated on a few common points. The decision rule approximations relax this requirement. They allow a separate transported distribution per source, but the static rules penalize differences in their means and the affine rules differences in their first and second moments. Accordingly, \textsc{static} and \textsc{affine} move the samples of different sources towards similar locations. \textsc{minkowski} and \textsc{mixture} impose no such link. Each source's samples are transported independently, and a few of them are moved far into the loss region. These uncoordinated moves inflate the aggregation-based bounds.

\subsection{Multi-Item Newsvendor Problem}\label{sec:newsvendor}
The second experiment asks whether the tighter bounds translate into better decisions out of sample. To test this, we consider a coffee shop that must decide how many small, regular and large cups of a coffee type to stock, but has only a few weeks of reliable sales data for it. Long sales histories of other coffee types are available, but they are partly corrupted by recording errors. Fusion DRO combines these sources by intersecting one ambiguity set around the scarce reliable data with one around each corrupted history. When each of five coffee types plays the role of the scarce product in turn, the affine approximation achieves the lowest out-of-sample cost in all five cases, whereas pooling all data, the natural default, incurs losses in every case.

\paragraph{Model.} The coffee shop chooses the numbers $\bm{\beta}\in\mathbb{R}^m_+$ of small, regular and large cups to stock ({\em i.e.}, $m=3$), before observing the respective demands $\bm{\xi}\in\mathbb{R}^m$. It incurs the newsvendor cost
\[
    \ell(\bm{\beta};\bm{\xi})=\bm{c}^\top\bm{\beta}-\bm{p}^\top\min\{\bm{\beta},\bm{\xi}\}+\bm{h}^\top\max\{\bm{\beta}-\bm{\xi},\bm{0}\}+\bm{b}^\top\max\{\bm{\xi}-\bm{\beta},\bm{0}\},
\]
where $\bm{c}$, $\bm{p}$, $\bm{h}$ and $\bm{b}$ are the unit ordering costs, selling prices, holding costs and backorder costs, and the minimum and maximum are taken componentwise. Negative costs correspond to profits.

\paragraph{Data.} Our case study builds on that of~\citet[\S~4.1.2]{lixianyu2026robust}, who use the Coffee Shop Sales dataset of Maven Roasters,\footnote{\url{https://mavenanalytics.io/data-playground/coffee-shop-sales}} with transactions recorded for 27 weeks from three New York City stores. Like them, we consider five coffee products sold in small, regular and large cups (Brazilian, Colombian, Ethiopian, Jamaican and Old Time Diner Blend), and we define the demand $\bm{\xi}$ for a product as its average daily sales of the three cup sizes over a week. We sum the sales over the three stores, so that each product provides $27$ weekly observations of $\bm{\xi}$. Following~\citet{lixianyu2026robust}, we use the modal Brazilian unit prices $\bm{p}=(2.2,3.0,3.5)$ for all products, and set $\bm{c}=(0.4,1.9,2.0)$, $\bm{h}=0.05\,\bm{p}$ and $\bm{b}=1.5\,\bm{p}$. Each coffee type in turn serves as the scarce target product. Five of the $27$ target observations, drawn at random, define the target reference distribution, and the remaining $22$ observations are used for testing. The full observation histories of the four other products define one reference distribution each, so that $K=5$. To mimic unreliable records, such as bulk orders or duplicated entries booked as regular demand, we multiply each observation of each of the other products by $4$ with a probability that differs across these products ($10\%$, $20\%$, $30\%$ and $40\%$ in random order). The target observations are never corrupted. The histories of the other products therefore overstate demand by different amounts, while the few target observations are accurate.

\paragraph{Experimental setup.} Methods that weight the sources (\textsc{affine}, \textsc{static}, \textsc{minkowski}, \textsc{mixture} and \textsc{barycenter}) assign weight $0.5$ to the target reference distribution and $0.125$ to each of the other four, so that the accurate target observations count as much as all corrupted histories combined. The empirical distribution~$\hat{\mathbb{P}}_\mathrm{bary}$, which defines \textsc{barycenter} as well as the radii of \textsc{affine} and \textsc{static}, consists of $N_\mathrm{bary}=128$ independent samples from~$\mathbb{P}_\mathrm{bary}$. In addition, each method except \textsc{erm} has a single hyperparameter that controls its conservatism, namely the parameter $\delta$ for \textsc{affine} and \textsc{static}, the radius~$\varepsilon$ for \textsc{pool}, \textsc{mixture} and \textsc{barycenter}, the common radius $\varepsilon_k=\varepsilon$ of the $K$ elementary ambiguity sets for \textsc{minkowski}, and the regularization parameter~$\nu$ for \textsc{rkhs}. We tune $\delta$ and $\varepsilon$ over the grid
\[
    \mathcal{E}\coloneqq\{0,0.005,0.01,0.05,0.1,0.2,0.5,1,2,5,10,20,50\}
\]
and $\nu$ over the coarser grid $\{0.1,0.3,0.5,0.7,0.9\}$. The hyperparameter of each method is selected by leave-one-out cross-validation on the five target observations. In each round, four observations define the target reference distribution and the fifth is used for validation, while the reference distributions of the other products remain unchanged. The method is then re-solved with all five target observations, using the hyperparameter with the lowest average validation cost, and evaluated on the $22$ test observations. For \textsc{rkhs}, we stop each solve after one hour and select~$\nu$ among the values for which all five rounds are solved within this limit. Each experiment is repeated $10$ times, and the target observations, the corrupted observations, the samples defining~$\hat{\mathbb{P}}_\mathrm{bary}$ and the random samples of \textsc{rkhs} are redrawn in every repetition.

\paragraph{Results.} Table~\ref{tab:coffee} reports the average test cost and the solve time for each target product. Three findings stand out. First, \textsc{affine} achieves the lowest cost for every target, earning profits of $31$ to $40$ per day. The closest competitors are \textsc{mixture} and \textsc{minkowski}, which nearly coincide, and \textsc{affine} improves on them by $2$ to $7$ per day. Second, how the sources are combined matters more than robustness itself. \textsc{erm} and \textsc{pool} incur positive costs ({\em i.e.}, losses) for all five targets, and \textsc{pool} essentially reproduces \textsc{erm}, since cross-validation selects a radius of nearly zero. Pooling weights each source by its sample size, so the five accurate target weeks are swamped by more than a hundred historical weeks, many of them inflated, and the resulting orders are far too large. Intersecting the ambiguity sets avoids this. Every candidate distribution must lie in the target's ball as well, which prevents the inflated histories from dictating the worst case. Third, \textsc{affine} outperforms \textsc{static} by a wide margin. This gap is not a consequence of Corollary~\ref{cor:static_upper_bound}, which compares bounds for the same radii, whereas here each method selects its own hyperparameters; \textsc{static} even falls behind \textsc{mixture} and \textsc{minkowski} for the same reason. A plausible explanation is that the corruption inflates second moments much more than means, and only the affine decision rules couple the second moments of the transported distributions across sources. The price of \textsc{affine} is computational. Solving the underlying semidefinite program takes about $19$ seconds, compared with fractions of a second for the second-order cone programs. This remains far below the $1,150$ to $1,900$ seconds of \textsc{rkhs}, which also performs worst. As observed by~\citet{han2026robust}, proximity in the RKHS does not control how far the admissible distributions spread their mass, and here the worst case concentrates on the largest observed demands. 

\begin{table}[t]
\centering
\fontsize{8}{9.5}\selectfont
\setlength{\tabcolsep}{3.5pt}
\renewcommand{\arraystretch}{0.95}
\begin{tabular}{@{}llrrrrr@{}}
\toprule
Method & Metric & Brazilian & Colombian & Ethiopian & Jamaican & Old Time \\
\midrule
\multirow{2}{*}{\textsc{affine}} & Test cost & \textbf{-40.40} & \textbf{-33.92} & \textbf{-31.14} & \textbf{-33.41} & \textbf{-36.70} \\
 & Time & 18.78 & 18.59 & 18.78 & 18.50 & 18.76 \\ \midrule
\multirow{2}{*}{\textsc{static}} & Test cost & -19.67 & -6.40 & -4.31 & -4.29 & -14.67 \\
 & Time & 0.17 & 0.16 & 0.17 & 0.16 & 0.16 \\ \midrule
\multirow{2}{*}{\textsc{minkowski}} & Test cost & -37.01 & -26.83 & -29.04 & -31.12 & -35.17 \\
 & Time & 0.14 & 0.15 & 0.14 & 0.14 & 0.14 \\ \midrule
\multirow{2}{*}{\textsc{mixture}} & Test cost & -37.01 & -26.82 & -29.04 & -31.12 & -35.19 \\
 & Time & \textbf{0.04} & \textbf{0.04} & \textbf{0.04} & \textbf{0.04} & \textbf{0.04} \\ \midrule
\multirow{2}{*}{\textsc{barycenter}} & Test cost & -10.07 & 1.75 & 5.53 & 7.40 & -1.25 \\
 & Time & \textbf{0.04} & \textbf{0.04} & \textbf{0.04} & \textbf{0.04} & \textbf{0.04} \\ \midrule
\multirow{2}{*}{\textsc{pool}} & Test cost & 0.74 & 29.96 & 32.04 & 32.30 & 15.38 \\
 & Time & \textbf{0.04} & \textbf{0.04} & \textbf{0.04} & \textbf{0.04} & \textbf{0.04} \\ \midrule
\multirow{2}{*}{\textsc{erm}} & Test cost & 0.74 & 29.96 & 32.04 & 32.30 & 14.33 \\
 & Time & \textbf{0.04} & \textbf{0.04} & \textbf{0.04} & \textbf{0.04} & \textbf{0.04} \\ \midrule
\multirow{2}{*}{\textsc{rkhs}} & Test cost & 187.54 & 142.64 & 142.50 & 141.10 & 118.06 \\
 & Time & 1903.04 & 1560.54 & 1437.84 & 1793.16 & 1150.52 \\
\bottomrule
\end{tabular}
\caption{Average daily cost on the test observations and time in seconds to solve the final problem with the selected hyperparameter, excluding cross-validation, for each method and target product, averaged over $10$ repetitions. Bold entries mark the best value per column.}
\label{tab:coffee}
\end{table}

\begin{remark}[Intermediate Decision Rules]
In this experiment, \textsc{affine} performs best but requires solving semidefinite programs, whereas \textsc{static} is cheap but performs considerably worse. Intermediate decision rules can combine the strengths of both. For instance, in separate experiments on the newsvendor instances above, we restricted the matrices~$\bm{A}_k$ in~\eqref{eq:affine-decision-rule} to be diagonal and used the radii selected for \textsc{affine}. The resulting test costs were within about $5\%$ of those of \textsc{affine}, and the solve times were below half a second, as for \textsc{static}. We do not report these results in detail, because which restriction works well depends on the problem at hand, and a systematic study of intermediate decision rules is beyond the scope of this paper.
\end{remark}

\subsection{Neural Network Training for Classification}\label{sec:SGD-experiment}
The third experiment leaves the comfort zone of convex losses and considers training a neural network for a binary classification task (credit default). 
In contrast with the earlier experiments, we consider ambiguity sets that are not Wasserstein balls around the empirical distribution but construct two ambiguity sets based on moment information (for the continuous features) and the total variation distance (for the categorical features). As the resulting problem is nonconvex and has about $60{,}000$ decision variables, we solve it with the subgradient method of Section~\ref{sec:sgd}.~The affine approximation outperforms every applicable benchmark at every sample size, and for~$6{,}000$ records it roughly halves the suboptimality of training over either elementary ambiguity set alone.

\paragraph{Model.} We train a neural network to predict a binary label $y\in\{0,1\}$ from numerical features $\bm{x}_\mathrm{num}$ that range over the unit $\ell_\infty$-ball $\mathcal{N}\subseteq\mathbb{R}^{m_\mathrm{num}}$ and categorical features $\bm{x}_\mathrm{cat}$ that range over a finite set $\mathcal{C}\subseteq\mathbb{Z}^{m_\mathrm{cat}}$. The uncertain parameter $\bm{\xi}=(\bm{x}_\mathrm{num},\bm{x}_\mathrm{cat},y)$ thus ranges over the set $\Xi\coloneqq\mathcal{N}\times\mathcal{C}\times\{0,1\}$, which is compact as required by Assumption~\ref{assump:general}\,\ref{assump:general-support}. The decision $\bm{\beta}\in\mathbb{R}^n$ collects the weights of a fully connected neural network that maps the features $\bm{x}\coloneqq(\bm{x}_\mathrm{num},\bm{x}_\mathrm{cat})$ to a logit $f_{\bm{\beta}}(\bm{x})$, \textit{i.e.}, to the predicted log-odds that $y=1$. The network first one-hot encodes the categorical features and then passes all features through three hidden layers of widths $256$, $128$ and $64$ with SiLU activations~\citep{elfwing2018sigmoid} and a linear output layer. The decision maker incurs the cross-entropy loss $\ell(\bm{\beta};\bm{\xi})=-y\log\sigma(f_{\bm{\beta}}(\bm{x}))-(1-y)\log(1-\sigma(f_{\bm{\beta}}(\bm{x})))$, where $\sigma(t)\coloneqq1/(1+e^{-t})$. This loss is nonconvex in $\bm{\beta}$, but the smoothness of the activation functions makes it weakly convex and Lipschitz continuous in~$\bm{\beta}$ on every bounded convex set.

The fusion ambiguity set $\mathcal{P}=\mathfrak{B}_1\cap\mathfrak{B}_2$ is the intersection of two elementary ambiguity sets. The set $\mathfrak{B}_1$ contains all distributions on~$\Xi$ for which the vector of distinct first and second moments of~$(\bm{x}_\mathrm{num},y)$ lies within Euclidean distance~$\varepsilon_1$ of a noisy estimate. As this restricts the mean and covariance matrix of~$(\bm{x}_\mathrm{num},y)$, $\mathfrak{B}_1$ constitutes a moment ambiguity set in the sense of Example~\ref{ex:moment-ambiguity-sets}, and it leaves the categorical features free. The set~$\mathfrak{B}_2$ is defined through a partition of the values of~$(\bm{x}_\mathrm{cat},y)$ into finitely many cells. It contains all distributions on~$\Xi$ under which the probabilities of the cells lie within total variation distance~$\varepsilon_2$ of a noisy histogram. Thus, $\mathfrak{B}_2$ leaves the numerical features free. The radii~$\varepsilon_1$ and~$\varepsilon_2$ are known and ensure that each set contains the unknown distribution of interest with a prescribed confidence. Such radii arise naturally when a data holder adds Gaussian noise to the moments \citep{balle2018improving,biswas2020coinpress} or Laplace noise to the cell counts \citep{hay2010boosting} to protect privacy, or when the statistics are estimated from a random sample \citep{hotelling1931generalization,delage2010distributionally,weissman2003inequalities}; see Appendix~\ref{app:sgd-experiment}. Both decision rule approximations of~\eqref{eq:classical_dro} use the lifted decision rules of Remark~\ref{rem:lifted-rules}, where $\bm{s}(\bm{\xi})$ stacks the numerical features, the one-hot encoded categorical features, the label and the cell indicators, and where the quadratic terms omit these indicators.

\paragraph{Data.} We use the Default of Credit Card Clients dataset,\footnote{\url{https://archive.ics.uci.edu/dataset/350/default+of+credit+card+clients}} which records whether each of $30{,}000$ credit card clients defaults in the following month. The $m_\mathrm{num}=14$ numerical features are the credit limit, the age, six monthly bill amounts and six monthly payment amounts, where monetary amounts are log-transformed and all numerical features are scaled to $[-1,1]$ using data-independent bounds. The $m_\mathrm{cat}=9$ categorical features are the sex, the education level, the marital status and six monthly repayment statuses, which add $45$ one-hot inputs to the network. The histogram has $432$ cells, each fixing the label and an aggregated version of the categorical features (Appendix~\ref{app:sgd-experiment}). In each simulation, we sample $N$ records, which form the hidden empirical distribution $\mathbb{P}_N$. We compute the $134$ distinct moments and the histogram of~$\mathbb{P}_N$, perturb them with Gaussian and Laplace noise, respectively, and choose the radii $\varepsilon_1$ and $\varepsilon_2$ so that the two sets contain the moments and the histogram of $\mathbb{P}_N$ with $95\%$ confidence each (Appendix~\ref{app:sgd-experiment}).\footnote{In every simulation, we certify that $\mathfrak{B}_1\cap\mathfrak{B}_2$ is nonempty, as required by Assumption~\ref{assump:general}\,\ref{assump:general-nonempty-intersection}, by constructing a finitely supported distribution on~$\Xi$ that belongs to both sets; cf.~Remark~\ref{rem:empty-intersection}.} We vary $N$ from $2{,}000$ to $18{,}000$ and reserve $12{,}000$ further records for testing.

\paragraph{Experimental setup.} We compare \textsc{affine} and \textsc{static} over $\mathfrak{B}_1\cap\mathfrak{B}_2$ with \textsc{minkowski} (Section~\ref{sec:minkowski}) and with \textsc{moment} and \textsc{tv}, which train over $\mathfrak{B}_1$ and $\mathfrak{B}_2$ alone, respectively. As both elementary ambiguity sets describe the same hidden distribution, we use the uniform weights $t_1=t_2=1/2$. The sample-based benchmarks of Table~\ref{tab:baselines} do not apply. Indeed, \textsc{pool}, \textsc{mixture} and \textsc{barycenter} require reference distributions on~$\Xi$, and \textsc{erm} and \textsc{rkhs} require samples, which the decision maker does not have. For reference, \textsc{oracle} minimizes the empirical risk under the hidden distribution $\mathbb{P}_N$, using stochastic gradient descent and stopping at the update with the lowest test loss. It uses both the hidden records and the test data and is thus unattainable. All robust models are trained with the subgradient method of Section~\ref{sec:sgd} for $75{,}000$ updates. Each update evaluates the worst-case expectations of the loss augmented by its correction terms over~$\mathfrak{B}_1$ and~$\mathfrak{B}_2$, and Danskin's theorem yields their gradients from the corresponding maximizers. For the moment set~$\mathfrak{B}_1$, we replace the worst-case expectation by its dual \citep[theorem~4.5]{kuhn2025distributionally} and update the dual multipliers together with the network weights (Section~\ref{sec:gradient-other}), which requires one maximization over~$\Xi$. For the histogram set~$\mathfrak{B}_2$, the worst-case expectation is available in closed form. It places the probability of each cell on a maximizer of the augmented loss within that cell and shifts probability~$\varepsilon_2$ from the cells with the smallest maxima to the cell with the largest maximum \citep[proposition~6.13]{kuhn2025distributionally}, which requires one maximization per cell. These maximizations are nonconvex and, because of the categorical features, combinatorial. We approximate them by multi-start projected gradient ascent with rounding to valid one-hot encodings (Appendix~\ref{app:sgd-experiment}). A monolithic approach is out of reach, as \textsc{affine} would require a nonconvex program over $56{,}577$ network weights and $4{,}644$ decision rule coefficients with semi-infinite constraints over the mixed-integer support~$\Xi$. The subgradient method, by contrast, runs ten simulations in parallel on a single NVIDIA GH200 GPU in $2.1$ hours for \textsc{affine} and $2.0$ hours for \textsc{static}. The quadratic terms of the affine rules thus come almost for free, whereas the semidefinite program made \textsc{affine} about $110$ times slower than \textsc{static} in Section~\ref{sec:newsvendor}.

\paragraph{Results.} Table~\ref{tab:SGD_result} reports the mean suboptimality of each method over $50$ simulations, defined as the relative excess of its out-of-sample cross-entropy over that of \textsc{oracle}. The out-of-sample cross-entropy of \textsc{oracle} itself ranges from $0.430$ to $0.440$ across all~$N$. Each simulation redraws the hidden records, the test records, the noise and the initial network weights. \textsc{affine} attains the lowest out-of-sample cross-entropy at every sample size and improves on \textsc{static} throughout. The elementary ambiguity sets fail for complementary reasons. The worst cases of \textsc{moment} can rearrange the categorical features freely, including the repayment statuses that carry most of the predictive information. As a consequence, its suboptimality stalls at about $25\%$, however large $N$ becomes. \textsc{tv} constrains the categorical features but leaves the numerical ones free. It improves with $N$ but remains above \textsc{affine}. Which elementary set is better depends on $N$. \textsc{moment} wins for $N\leq 6{,}000$ and \textsc{tv} for $N\geq 10{,}000$. \textsc{affine} beats both at every $N$, so the decision maker need not know which statistic is more informative, and for $N=6{,}000$ it roughly halves the suboptimality of both. Aggregation, in contrast, combines the weaknesses of both sets. \textsc{minkowski} averages two worst cases, each of which exploits the features that its own set leaves free, and for $N\geq 6{,}000$ it performs worse than either set alone. Finally, the theory guarantees less than what we observe. With uniform weights, Proposition~\ref{prop:mink} bounds the (in-sample) worst-case values of \textsc{static} and \textsc{affine} only by the average of the two single-set values, and indeed \textsc{static} falls behind \textsc{moment} for $N=2{,}000$. \textsc{affine}, however, beats both elementary sets at every~$N$.\looseness=-1

\begin{table}[t]
\centering
\caption{Average suboptimality (in \%) of each method relative to \textsc{oracle} across $50$ simulations, with standard errors in parentheses. Bold entries mark the best method for each sample size~$N$.}
\label{tab:SGD_result}
\begin{tabular}{r*{5}{r@{\,}l}}
\toprule
$N$ & \multicolumn{2}{c}{\textsc{affine}} & \multicolumn{2}{c}{\textsc{static}} & \multicolumn{2}{c}{\textsc{minkowski}} & \multicolumn{2}{c}{\textsc{moment}} & \multicolumn{2}{c}{\textsc{tv}} \\
\midrule
$2{,}000$  & \textbf{37.83} & (0.59) & 50.31 & (0.40) & 56.32 & (0.24) & 44.66 & (0.50) & 57.08 & (0.16) \\
$6{,}000$  & \textbf{18.41} & (0.21) & 22.08 & (0.19) & 35.76 & (0.25) & 31.88 & (0.27) & 35.18 & (0.20) \\
$10{,}000$ & \textbf{10.13} & (0.13) & 12.20 & (0.11) & 29.51 & (0.14) & 27.41 & (0.16) & 19.65 & (1.02) \\
$14{,}000$ & \textbf{6.99} & (0.10) & 8.93 & (0.07) & 26.83 & (0.12) & 25.17 & (0.14) & 9.53 & (0.08) \\
$18{,}000$ & \textbf{5.59} & (0.08) & 7.59 & (0.07) & 25.61 & (0.09) & 24.06 & (0.11) & 8.05 & (0.07) \\
\bottomrule
\end{tabular}
\end{table}

\paragraph{Acknowledgments}
This work was supported as a part of NCCR Automation, a National Centre of Competence in Research, funded by the Swiss National Science Foundation (grant number 51NF40\_225155). The authors acknowledge the use of resources provided by the Isambard-AI National AI Research Resource (AIRR). Isambard-AI is operated by the University of Bristol and is funded by the UK Government’s Department for Business, Innovation, Science and Trade (DBIST) via UK Research and Innovation; and the Science and Technology Facilities Council [ST/AIRR/I-A-I/1023].

\bibliographystyle{informs2014} 
\bibliography{refinf.bib}

\newpage
\appendix
\section*{Appendix}
\section{Proofs}\label{sec:appendix-proofs}
Below we provide the proofs for all theoretical results in the main text. 

\begin{proof}[Proof of Proposition~\ref{prop:reformulation_main}]
We first prove weak duality ($v_\mathrm{P}\leq v_\mathrm{D}$). To this end, let $\mathbb{Q}$ be feasible in~\eqref{prob:sup-primal} and let $\{\lambda_k\}_k$ be feasible in~\eqref{eq:reformulation_dual}. Then, we have
\begin{align*}
    \mathbb{E}_{\mathbb{Q}}[\ell(\bm{\xi})]\leq\mathbb{E}_{\mathbb{Q}}\bigg[\sum_{k=1}^K\lambda_k(\bm{\xi})\bigg]\leq\sum_{k=1}^K\mathbb{E}_{\mathbb{Q}}[\lambda_k(\bm{\xi})]\leq\sum_{k=1}^K\sup_{\mathbb{Q}_k\in\mathfrak{B}_k}\mathbb{E}_{\mathbb{Q}_k}[\lambda_k(\bm{\xi})],
\end{align*}
where the first inequality follows from the feasibility of $\{\lambda_k\}_k$ in~\eqref{eq:reformulation_dual} and the monotonicity of the expectation. The second inequality uses the subadditivity of the expectation. The right-hand side of this inequality equals~$\infty$ unless the positive part of every~$\lambda_k$ is $\mathbb{Q}$-integrable, in which case the inequality holds as an equality. The third inequality holds because $\mathbb{Q}$ is feasible in~\eqref{prob:sup-primal}, which implies that $\mathbb{Q}\in\mathfrak{B}_k$ for every $k\in[K]$. Taking the supremum over all feasible~$\mathbb{Q}$ on the left-hand side and the infimum over all feasible $\{\lambda_k\}_k$ on the right-hand side yields $v_\mathrm{P}\leq v_\mathrm{D}$.

Next, we prove strong duality ($v_\mathrm{P} = v_\mathrm{D}$) under Assumption~\ref{assump:general}. By introducing auxiliary variables~$\mathbb{Q}_k\in\mathfrak{B}_k$, $k\in[K]$, and enforcing $\mathbb{Q}_k=\mathbb{Q}$ via equality constraints, we can reformulate~\eqref{prob:sup-primal} as \looseness=-1
\begin{align}\label{prob:sup-primal_alternative}
    v_\mathrm{P} = \left\{\begin{array}{cll}
        \sup & \mathbb{E}_{\mathbb{Q}}[\ell(\bm{\xi})] \\
        \mathrm{s.t.} & \mathbb{Q}\in\cP_0(\Xi),\;\mathbb{Q}_k\in\mathfrak{B}_k & \forall k\in[K] \\
        & \mathbb{Q}_k=\mathbb{Q} & \forall k\in[K].
    \end{array}\right.
\end{align}
Indeed, $\mathbb{Q}\in\cP_0(\Xi)$ satisfies $\mathbb{Q}_k=\mathbb{Q}\in\mathfrak{B}_k$ for all~$k\in[K]$ if and only if $\mathbb{Q}\in\cap_{k=1}^K\mathfrak{B}_k=\cP$. Dualizing the equality constraints with Lagrange multipliers~$\lambda_k\in C(\Xi)$ then yields
\begin{align}
\label{prob:decoupled_supinf_formulation}
    v_\mathrm{P} = \sup_{\mathbb{Q}\in\cP_0(\Xi)}\sup_{\{\mathbb{Q}_k\in\mathfrak{B}_k\}_k} \inf_{\{\lambda_k\in C(\Xi)\}_k} \mathbb{E}_{\mathbb{Q}}\!\left[\ell(\bm{\xi})-\sum_{k=1}^K\lambda_k(\bm{\xi})\right] + \sum_{k=1}^K\mathbb{E}_{\mathbb{Q}_k}[\lambda_k(\bm{\xi})].
\end{align}
To see this, recall that $\Xi$ is compact such that every $\lambda_k\in C(\Xi)$ is bounded, and that bounded continuous functions separate probability measures \citep[lemma~9.3.2]{dudley2002real}. If $\mathbb{Q}_k\neq\mathbb{Q}$ for some $k\in[K]$, there thus exists $\lambda_k\in C(\Xi)$ with $\mathbb{E}_{\mathbb{Q}_k}[\lambda_k(\bm{\xi})]\neq\mathbb{E}_{\mathbb{Q}}[\lambda_k(\bm{\xi})]$, and scaling~$\lambda_k$ drives the inner infimum in~\eqref{prob:decoupled_supinf_formulation} to~$-\infty$. If $\mathbb{Q}_k=\mathbb{Q}$ for all $k\in[K]$, in contrast, all terms involving the multipliers cancel, and the inner infimum equals $\mathbb{E}_{\mathbb{Q}}[\ell(\bm{\xi})]$.

Interchanging the order of maximization and minimization in~\eqref{prob:decoupled_supinf_formulation} leads to the problem
\begin{align}\label{eq:lagrange_UB}
    v_\mathrm{D} = \inf_{\{\lambda_k\in C(\Xi)\}_k} \sup_{\mathbb{Q}\in\cP_0(\Xi)} \mathbb{E}_{\mathbb{Q}}\!\left[\ell(\bm{\xi})-\sum_{k=1}^K\lambda_k(\bm{\xi})\right] + \sum_{k=1}^K\sup_{\mathbb{Q}_k\in\mathfrak{B}_k}\mathbb{E}_{\mathbb{Q}_k}[\lambda_k(\bm{\xi})].
\end{align}
To prove that the right-hand side of~\eqref{eq:lagrange_UB} is indeed equivalent to Problem~\eqref{eq:reformulation_dual}, we first note that 
\begin{align*}
    \sup_{\mathbb{Q}\in\cP_0(\Xi)}\mathbb{E}_{\mathbb{Q}}\!\left[\ell(\bm{\xi})-\sum_{k=1}^K\lambda_k(\bm{\xi})\right] = \sup_{\bm{\xi}\in\Xi}\left\{\ell(\bm{\xi})-\sum_{k=1}^K\lambda_k(\bm{\xi})\right\}.
\end{align*}
Clearly, the supremum over~$\mathbb{Q}\in\cP_0(\Xi)$ is smaller or equal to the supremum over~$\bm{\xi}\in\Xi$, and the converse inequality holds because $\cP_0(\Xi)$ contains the Dirac measure $\delta_{\bm{\xi}}$ for every $\bm{\xi}\in\Xi$. Introducing an epigraphical variable $\tau$ then allows us to recast the right-hand side of~\eqref{eq:lagrange_UB} as
\begin{align*}
\begin{array}{lll}
    \displaystyle \underset{}{\inf} & \displaystyle \tau
    + \sum_{k=1}^K \sup_{\mathbb{Q}_k \in \mathfrak{B}_{k}} \mathbb{E}_{\mathbb{Q}_k}[\lambda_k(\bm{\xi})]
    \\
    \mathrm{s.t.} & \tau\in\mathbb R, ~ \lambda_k \in C(\Xi) \quad&\forall k\in[K] \\
    & \ell(\bm{\xi}) - \sum_{k=1}^K\lambda_k(\bm{\xi})\leq \tau \quad &\forall \bm{\xi} \in \Xi.
\end{array}
\end{align*}
Without loss of generality, we may fix~$\tau=0$. To see this, consider any feasible solution $(\{\lambda_k\}_k,\tau)$, and define $(\{\lambda'_k\}_k, \tau')$ through $\lambda'_k=\lambda_k+\tau/K$ for every~$k\in[K]$ and $\tau'=0$. By construction, this new solution remains feasible and has the same objective function value as $(\{\lambda_k\}_k,\tau)$. This reasoning shows that the right-hand side of~\eqref{eq:lagrange_UB} is thus indeed equivalent to Problem~\eqref{eq:reformulation_dual}. 

Strong duality now follows from Sion's minimax theorem for extended real-valued functions \citep[theorem~2.11]{bendavid2023new}. The feasible set $\cP_0(\Xi)\times\prod_{k=1}^K\mathfrak{B}_k$ of the supremum is convex by Assumption~\ref{assump:general}\,\ref{assump:general-convex-ball} and weakly compact. Indeed, $\cP_0(\Xi)$ is weakly compact because $\Xi$ is compact \citep[proposition~3.6]{kuhn2025distributionally}, each~$\mathfrak{B}_k$ is weakly compact as a weakly closed subset of~$\cP_0(\Xi)$ by Assumption~\ref{assump:general}\,\ref{assump:general-convex-ball}, and finite products of compact sets are compact. The feasible set $C(\Xi)^K$ of the infimum is a vector space and hence convex. For fixed $(\lambda_1,\ldots,\lambda_K)\in C(\Xi)^K$, the Lagrangian is affine (hence concave) and weakly upper semicontinuous in $(\mathbb{Q},\mathbb{Q}_1,\ldots,\mathbb{Q}_K)$. Indeed, $\mathbb{E}_{\mathbb{Q}}[\ell(\bm{\xi})]$ is weakly upper semicontinuous in~$\mathbb Q$ by Assumptions~\ref{assump:general}\,\ref{assump:general-support} and~\ref{assump:general}\,\ref{assump:loss-usc} together with \citet[proposition~3.3]{kuhn2025distributionally}, and each expectation $\mathbb{E}_{\mathbb{Q}_k}[\lambda_k(\bm{\xi})]$ is weakly continuous in~$\mathbb Q_k$ because $\lambda_k\in C(\Xi)$ is continuous as well as bounded. For fixed $(\mathbb{Q},\mathbb{Q}_1,\ldots,\mathbb{Q}_K)$, the Lagrangian is identically $-\infty$ if $\mathbb{E}_{\bQ}[\ell(\bm{\xi})]=-\infty$. Otherwise, it is affine in $(\lambda_1,\ldots,\lambda_K)$ and continuous with respect to the supremum norm. In both cases, it is convex and lower semicontinuous in $(\lambda_1,\ldots,\lambda_K)$. Hence, Sion's minimax theorem applies, and we conclude that $v_\mathrm{P}=v_\mathrm{D}$.
\end{proof}

In order to prove Lemma~\ref{lem:pairwise-symmetric-rules}, we first derive an intermediate reformulation of Problem~\eqref{eq:reformulation_dual}.

\begin{lemma}\label{lem:adjustable_opt_trick}
For any weight vector $\bm{t}\in\Delta(K)$, Problem~\eqref{eq:reformulation_dual} is equivalent to
\begin{equation}
\label{eq:adjustable_reform}
\begin{aligned}
    \inf \quad & \sum_{k=1}^K \sup_{\mathbb{Q}_k\in\mathfrak{B}_k} \mathbb{E}_{\mathbb{Q}_k}[\lambda_k(\bm{\xi})] \\ \mathrm{s.t.} \quad & \lambda_k\in C(\Xi) \quad \forall k\in[K] \\ & \inf_{\{\bm{\alpha}'_k\in\mathcal{A}'\}_k} \sup_{\bm{\xi}\in\Xi} \sup_{\{\bm{\xi}_k\in\Xi\}_k} \sum_{k=1}^K \Bigl(t_k\ell(\bm{\xi}_k) - \lambda_k(\bm{\xi}_k) + \bm{\alpha}'_k(\bm{\xi}_1,\ldots,\bm{\xi}_K,\bm{\xi})^\top(\bm{\xi}_k-\bm{\xi})\Bigr) \leq 0, 
\end{aligned}
\end{equation}
where $\mathcal{A}'$ denotes the space of all functions $\bm{\alpha}':\Xi^{K+1}\to\mathbb{R}^m$. Thus, the optimal value of~\eqref{eq:adjustable_reform} coincides with~$v_\mathrm{D}$. If Assumption~\ref{assump:general} holds, then this optimal value also coincides with~$v_\mathrm{P}$.
\end{lemma}

\begin{proof}
    Problems~\eqref{eq:reformulation_dual} and~\eqref{eq:adjustable_reform} share the same objective function as well as the same domain constraints~$\lambda_k\in C(\Xi)$, $k\in[K]$. Thus, it suffices to show that, for any fixed $\{\lambda_k\}_k$, the inequality constraint of~\eqref{eq:reformulation_dual} holds if and only if that of~\eqref{eq:adjustable_reform} holds. We establish this through a chain of equivalences. Specifically, the inequality constraint of Problem~\eqref{eq:reformulation_dual} can be recast as
    \begin{align*}
        & \sup_{\bm{\xi} \in \Xi} \: \Big\{ \ell(\bm{\xi}) - \sum_{k=1}^K\lambda_k(\bm{\xi}) \Big\} \leq 0  \notag \\
        \iff ~ & \sup_{\bm{\xi} \in \Xi} \sup_{\{\bm{\xi}_k\in\Xi\}_k} \: \left\{ \sum_{k=1}^K \big( t_k \ell(\bm{\xi}_k)  - \lambda_k(\bm{\xi}_k) \big) : \bm{\xi}_k = \bm{\xi} \; \forall k\in[K]\right\} \leq 0,
    \end{align*}
    where the equivalence follows by introducing auxiliary variables~$\bm{\xi}_k$ constrained to equal~$\bm{\xi}$ and recalling that~$\bm{t}\in\Delta(K)$. Next, we dualize the equality constraints $\bm{\xi}_k = \bm{\xi}$ using Lagrange multipliers $\bm{\alpha}'_k \in \mathbb{R}^m$ to show that the inequality constraint of Problem~\eqref{eq:reformulation_dual} is equivalent to
     \begin{align*}
        \sup_{\bm{\xi} \in \Xi}  \sup_{\{\bm{\xi}_k\in\Xi\}_k} \: \inf_{\{\bm{\alpha}'_k \in \mathbb{R}^m\}_{k}} \: \left\{ \sum_{k=1}^K \big( t_k \ell(\bm{\xi}_k)  - \lambda_k(\bm{\xi}_k)   + {\bm{\alpha}'_k}^\top (\bm{\xi}_k - \bm{\xi}) \big) \right\} \leq 0.
    \end{align*}
    Finally, we replace the inner infimum over static Lagrange multipliers $\bm{\alpha}'_k \in \mathbb{R}^m$ with an outer infimum over functional Lagrange multipliers $\bm{\alpha}_k' \in \mathcal{A}'$, which adapt to $(\bm{\xi}_1,\ldots,\bm{\xi}_K,\bm{\xi})\in\Xi^{K+1}$. This yields exactly the inequality constraint in~\eqref{eq:adjustable_reform}, and thus the claim follows.
\end{proof}

Armed with Lemma~\ref{lem:adjustable_opt_trick}, we are now ready to prove Lemma~\ref{lem:pairwise-symmetric-rules}.
\begin{proof}[Proof of Lemma~\ref{lem:pairwise-symmetric-rules}]
    
    Fix any weight vector~$\bm{t}\in\Delta(K)$ and any continuous functions~$\lambda_k\in C(\Xi)$, $k\in[K]$. The claim follows if we can prove that the infimum over~$\{\bm\alpha_k\}_k$ in~\eqref{eq:prob-with-semi-infinite-constraint} matches the infimum over~$\{\bm\alpha'_k\}_k$ in~\eqref{eq:adjustable_reform}. We first show that the infimum in~\eqref{eq:adjustable_reform} is no larger than that in~\eqref{eq:prob-with-semi-infinite-constraint}. To see this, fix~$\bm\alpha_k\in\mathcal A$ feasible in~\eqref{eq:prob-with-semi-infinite-constraint} and define $\bm\alpha'_k\in\mathcal A'$ through $\bm{\alpha}'_k(\bm{\xi}_1,\ldots,\bm{\xi}_K,\bm{\xi}) =\bm{\alpha}_k(\bm{\xi}_k,\bm{\xi})$ for every $k\in[K]$. It is easy to see that the functions $\{\bm\alpha'_k\}_k$ constructed in this way belong to~$\mathcal A'$ and attain the same objective function value in~\eqref{eq:adjustable_reform} as~$\{\bm\alpha_k\}_k$ in~\eqref{eq:prob-with-semi-infinite-constraint}. Hence, the infimum in~\eqref{eq:adjustable_reform} is no larger than that in~\eqref{eq:prob-with-semi-infinite-constraint}. To prove the reverse inequality, fix~$\bm{\alpha}'_k\in\mathcal A'$ and define an auxiliary function $\lambda_k'\in C(\Xi)$ through $\lambda'_k(\bm{\xi})=\lambda_k(\bm{\xi})-t_k\sum_{j=1}^K\lambda_j(\bm{\xi})$ for every~$k\in[K]$. Note that~$\sum_k\lambda'_k(\bm{\xi})=0$ for all~$\bm{\xi}\in\Xi$ because~$\bm{t}\in\Delta(K)$. Next, define~$\bm{\alpha}_k\in\mathcal A$ through 
    \begin{align*}
     \bm{\alpha}_k(\bm{\xi}',\bm{\xi})
    =\begin{cases}
        \displaystyle\frac{\lambda'_k(\bm{\xi}')-\lambda'_k(\bm{\xi})}{\|\bm{\xi}'-\bm{\xi}\|_2^2}(\bm{\xi}'-\bm{\xi}) 
        &\text{if } \bm{\xi}'\neq\bm{\xi},\\[1ex]
        \bm{0} &\text{otherwise},
    \end{cases}
    \end{align*}
    for all $k\in[K]$. Note that $\bm{\alpha}_k(\bm{\xi}',\bm{\xi})$ is ostensibly invariant under permutations of~$\bm\xi$ and~$\bm \xi'$, which guarantees that~$\bm\alpha_k$ is indeed an element of~$\mathcal A$. Moreover, we have
\begin{align*}
    \sum_{k=1}^K\bm{\alpha}_k(\bm{\xi}',\bm{\xi})
    =\frac{\sum_k\lambda'_k(\bm{\xi}')-\sum_k\lambda'_k(\bm{\xi})}{\|\bm{\xi}'-\bm{\xi}\|_2^2}(\bm{\xi}'-\bm{\xi})
    =\bm{0} \quad \forall \bm\xi,\bm\xi'\in \Xi,
\end{align*}
which shows that $\{\bm\alpha_k\}_k$ is feasible in~\eqref{eq:prob-with-semi-infinite-constraint}. The definition of $\bm\alpha_k$ also implies that
\begin{align*}
    \bm{\alpha}_k(\bm{\xi}_k,\bm{\xi})^\top(\bm{\xi}_k-\bm{\xi})
    =\lambda'_k(\bm{\xi}_k)-\lambda'_k(\bm{\xi}) \quad \forall \bm \xi,\bm\xi_k\in\Xi.
\end{align*}
Thus, the objective function value of $\{\bm\alpha_k\}_k$ in~\eqref{eq:prob-with-semi-infinite-constraint} evaluates to
\begin{align*}
    &\sup_{\bm{\xi}\in\Xi}\sup_{\{\bm{\xi}_k\in\Xi\}_k} \sum_{k=1}^K\Bigl(t_k\ell(\bm{\xi}_k)-\lambda_k(\bm{\xi}_k)
        +\bm{\alpha}_k(\bm{\xi}_k,\bm{\xi})^\top(\bm{\xi}_k-\bm{\xi})\Bigr) \\ 
        & = ~ \sup_{\bm{\xi}\in\Xi}\sup_{\{\bm{\xi}_k\in\Xi\}_k} \sum_{k=1}^K \Bigl(t_k\ell(\bm{\xi}_k)-\lambda_k(\bm{\xi}_k)  +\lambda'_k(\bm{\xi}_k)-\lambda'_k(\bm{\xi})\Bigr)\\
    &= ~ \sup_{\bm{\xi}\in\Xi}\sup_{\{\bm{\xi}_k\in\Xi\}_k} \sum_{k=1}^Kt_k\bigg(\ell(\bm{\xi}_k)-\sum_{j=1}^K\lambda_j(\bm{\xi}_k)\bigg) \leq \sup_{\bm{\xi}'\in\Xi}\bigg\{\ell(\bm{\xi}')-
    \sum_{k=1}^K\lambda_k(\bm{\xi}')\bigg\} ,
\end{align*}
where the second equality holds because~$\sum_k\lambda'_k(\bm{\xi})=0$. Next, we compare the resulting upper bound with the objective function value of $\{\bm\alpha'_k\}_k$ in~\eqref{eq:adjustable_reform}, which evaluates to 
\begin{align*}
    \sup_{\bm{\xi}\in\Xi}\sup_{\{\bm{\xi}_k\in\Xi\}_k}
    \sum_{k=1}^K\Bigl(t_k\ell(\bm{\xi}_k)-\lambda_k(\bm{\xi}_k)
    +\bm{\alpha}'_k(\bm{\xi}_1,\ldots,\bm{\xi}_K,\bm{\xi})^\top(\bm{\xi}_k-\bm{\xi})\Bigr) \geq 
    \sup_{\bm{\xi}'\in\Xi}\bigg\{\ell(\bm{\xi}')-
    \sum_{k=1}^K\lambda_k(\bm{\xi}')\bigg\},
\end{align*}
where the inequality follows from the restriction $\bm\xi_1=\cdots=\bm\xi_K=\bm\xi$. Combining the two above estimates reveals that the objective function value of $\{\bm\alpha_k\}_k$ in~\eqref{eq:prob-with-semi-infinite-constraint} is no larger than that of~$\{\bm\alpha_k'\}_k$ in~\eqref{eq:adjustable_reform}. Since this construction applies to every~$\{\bm\alpha'_k\}_k$ feasible in~\eqref{eq:adjustable_reform}, the infimum in~\eqref{eq:prob-with-semi-infinite-constraint} is no larger than that in~\eqref{eq:adjustable_reform}. In summary, we have shown that the inequality constraints in~\eqref{eq:prob-with-semi-infinite-constraint} and~\eqref{eq:adjustable_reform} are equivalent, and therefore the claim follows.
\end{proof}

The proof of Proposition~\ref{prop:affine_upper_bound} relies on the following lemma, which shows that, for a specific class of optimization problems with translation-invariant objective and constraint functions, the inequality constraints can be aggregated without sacrificing optimality.

\begin{lemma}[Lossless Constraint Aggregation]\label{lemma:decouple}
    Let $\Lambda$ be a nonempty family of real-valued functions on $\Xi$ that is closed under constant shifts, i.e., if $\lambda\in\Lambda$ then $ \lambda+\alpha\in\Lambda$ for every $\alpha\in\mathbb{R}$. Let $S_k:\Lambda\to(-\infty,\infty]$ and $T_k:\Lambda\to[-\infty,\infty]$, $k\in[K]$, be translation invariant functionals with 
    \[
    S_k(\lambda+\alpha)=S_k(\lambda)-\alpha \quad \text{and} \quad T_k(\lambda+\alpha)=T_k(\lambda)+\alpha
        \quad \forall\lambda\in\Lambda, \alpha\in\mathbb{R}. 
    \]
    Then, the following identity holds
    \begin{align*}
    \mspace{-25mu}
        \inf\bigg\{\sum_{k=1}^K T_k(\lambda_k) :
            \lambda_k\in\Lambda \; \forall k\in[K],\;
            \sum_{k=1}^K S_k(\lambda_k)\leq 0\bigg\}\\
        =
        \inf\bigg\{\sum_{k=1}^K T_k(\lambda_k) :
            \lambda_k\in\Lambda,\;
            S_k(\lambda_k)\leq 0\;\forall k \in [K] \bigg\}.
    \end{align*}
\end{lemma}
\begin{proof}
    Let $(\lambda_1, \dots, \lambda_K)$ be a feasible solution of the disaggregated problem with individual constraints $S_k(\lambda_k) \leq 0$ for all $k \in [K]$. By adding these $K$ inequalities we obtain $\sum_{k=1}^K S_k(\lambda_k) \leq 0$, which implies that $(\lambda_1, \dots, \lambda_K)$ satisfies the single constraint of the aggregated problem. Hence, the aggregated problem constitutes a relaxation of the disaggregated problem, and the infimum of the aggregated problem must be smaller than (or equal to) that of the disaggregated problem.
    
    To establish the reverse inequality, let $(\lambda_1, \dots, \lambda_K)$ be any feasible solution of the aggregated problem, and construct a family of shifted functions $\lambda'_k \coloneqq \lambda_k + \alpha_k$ with constant shifts
    \[
        \alpha_k \coloneqq S_k(\lambda_k) - \frac{1}{K}\sum_{j=1}^K S_j(\lambda_j) \quad \forall k \in [K].
    \]
    These shifts are real numbers because $\sum_{k=1}^K S_k(\lambda_k)\leq 0$ and no $S_k(\lambda_k)$ equals~$-\infty$. As $\Lambda$ is closed under constant shifts, we have $\lambda'_k \in \Lambda$. In addition, we find
    \[
        S_k(\lambda'_k) = S_k(\lambda_k) - \alpha_k = \frac{1}{K}\sum_{j=1}^K S_j(\lambda_j) \leq 0 \quad \forall k \in [K],
    \]
    where the first equality exploits the translation invariance of $S_k$, the second equality follows from the definition of $\alpha_k$, and the inequality holds because $(\lambda_1, \dots, \lambda_K)$ is feasible in the aggregated problem. Thus, $(\lambda'_1, \dots, \lambda'_K)$ is feasible in the disaggregated problem. Similarly, we find
    \[
        \sum_{k=1}^K T_k(\lambda'_k) = \sum_{k=1}^K \big( T_k(\lambda_k) + \alpha_k \big) = \sum_{k=1}^K T_k(\lambda_k),
    \]
    where the two equalities follow from the translation invariance of~$T_k$ and 
    the observation that the shifts $\alpha_k$, $k\in[K]$, sum to~$0$. In summary, we have thus shown that for any feasible solution of the aggregated problem there exists a feasible solution of the disaggregated problem with the same objective function value. Hence, the optimal value of the aggregated problem is non-inferior to that of the disaggregated problem. This establishes the equality and completes the proof.
\end{proof} 

With this lemma in hand, we are now ready to prove Proposition~\ref{prop:affine_upper_bound}.

\begin{proof}[Proof of Proposition~\ref{prop:affine_upper_bound}] Substituting~\eqref{eq:affine-decision-rule} into the inequality constraint of Problem~\eqref{eq:prob-with-semi-infinite-constraint} yields 
\begin{align}\label{eq:inf_constraint_explicit}
    \left.\begin{array}{cl}
        \inf & \displaystyle\sup_{\bm{\xi}\in\Xi}\sup_{\{\bm{\xi}_k\in\Xi\}_k}\sum_{k=1}^K\big(t_k\ell(\bm{\xi}_k)-\lambda_k(\bm{\xi}_k)+(\bm{a}_k+\bm{A}_k(\bm{\xi}_k+\bm{\xi}))^\top(\bm{\xi}_k-\bm{\xi})\big) \\[3ex]
        \mathrm{s.t.} & \bm a_k\in\mathbb R^m,\; \bm A_k\in\mathbb R^{m\times m}\quad \forall k\in[K] \\[1ex]
        & \displaystyle \sum_{k=1}^K\bm{a}_k+\bigg(\sum_{k=1}^K\bm{A}_k\bigg)(\bm{\xi}'+\bm{\xi})
    =\bm{0}\quad\forall\bm{\xi}',\bm{\xi}\in\Xi
    \end{array}\right\}\leq 0.
\end{align}
The semi-infinite equality constraint can be restricted to the zero-sum conditions~$\sum_{k=1}^K\bm{a}_k=\bm{0}$ and~$\sum_{k=1}^K \bm{A}_k=\bm{0}$ without increasing the optimal value of the minimization problem in~\eqref{eq:inf_constraint_explicit}. To see this, suppose that $\{\bm a_k\}_k$ and $\{\bm A_k\}_k$ are feasible in~\eqref{eq:inf_constraint_explicit}, and set
\begin{align*}
    \bm{a}'_k=\bm{a}_k-t_k\sum_{j=1}^K\bm{a}_j
    \quad \text{and} \quad
    \bm{A}'_k=\bm{A}_k-t_k\sum_{j=1}^K\bm{A}_j
    \quad\forall k\in[K].
\end{align*}
Since~$\bm{t}\in\Delta(K)$, these updated coefficients satisfy~$\sum_{k=1}^K\bm{a}'_k=\bm{0}$ and~$\sum_{k=1}^K\bm{A}'_k=\bm{0}$. In addition, the semi-infinite equality constraint in~\eqref{eq:inf_constraint_explicit} ensures that $\bm{a}'_k+\bm{A}'_k(\bm{\xi}'+\bm{\xi}) = \bm{a}_k+\bm{A}_k(\bm{\xi}'+\bm{\xi})$ for all~$\bm{\xi}',\bm{\xi}\in\Xi$ and~$k\in[K]$. Thus, replacing~$\bm{a}_k,\bm{A}_k$ with~$\bm{a}'_k,\bm{A}'_k$ has no impact on the objective function in~\eqref{eq:inf_constraint_explicit}. This confirms that the semi-infinite equality constraint in~\eqref{eq:inf_constraint_explicit} may indeed be replaced with the zero-sum conditions~$\sum_{k=1}^K\bm{a}_k=\bm{0}$ and~$\sum_{k=1}^K \bm{A}_k=\bm{0}$.By suitably adjusting~$\bm a_k$, we can further restrict each~$\bm A_k$ to~$\mathbb S^m$ without sacrificing optimality. Hence, replacing the set of feasible decision rule coefficients with
\[
    \mathcal{Z} \coloneqq \Big\{\big(\{\bm{a}_k\}_k,\{\bm{A}_k\}_k\big)\in(\mathbb{R}^m)^K\times(\mathbb{S}^m)^K : \textstyle\sum_{k=1}^K\bm{a}_k=\bm{0},\;\sum_{k=1}^K\bm{A}_k=\bm{0}\Big\}
\]
does not change the infimum in~\eqref{eq:inf_constraint_explicit}. In summary, the above reasoning shows that the restriction of Problem~\eqref{eq:prob-with-semi-infinite-constraint} to affine decision rules of the form~\eqref{eq:affine-decision-rule} can be expressed equivalently as
    \begin{equation}
    \begin{aligned}
    \label{eq:inf_constraint}
    \inf \quad & \sum_{k=1}^K \sup_{\mathbb{Q}_k\in\mathfrak{B}_k} \mathbb{E}_{\mathbb{Q}_k}[\lambda_k(\bm{\xi})] \\ \mathrm{s.t.} \quad & \lambda_k\in C(\Xi) \quad \forall k\in[K] \\[-0.5ex]
    & \hspace{-2.5mm} \left.\begin{array}{cl}
        \inf & \displaystyle\sup_{\bm{\xi}\in\Xi}\sup_{\{\bm{\xi}_k\in\Xi\}_k}\sum_{k=1}^K\big(t_k\ell(\bm{\xi}_k)-\lambda_k(\bm{\xi}_k)+(\bm{a}_k+\bm{A}_k(\bm{\xi}_k+\bm{\xi}))^\top(\bm{\xi}_k-\bm{\xi})\big) \\
        \mathrm{s.t.} & \big(\{\bm{a}_k\}_k,\{\bm{A}_k\}_k\big)\in\mathcal{Z}
    \end{array}\right\}\leq 0.
    \end{aligned}
\end{equation}
Going forward, we denote the optimal value of Problem~\eqref{eq:inf_constraint} by~$v_\mathrm{inf}$. The terms depending on~$\bm{\xi}$ in the inequality constraint of~\eqref{eq:inf_constraint} are given by $-\sum_{k=1}^K\bm{a}_k^\top\bm{\xi}-\sum_{k=1}^K\langle\bm{A}_k,\bm{\xi}\bm{\xi}^\top\rangle$ and vanish thanks to the zero-sum conditions. Thus, the variable~$\bm{\xi}$ drops out. Next, we apply a standard trick from robust optimization and approximate the inequality constraint in~\eqref{eq:inf_constraint} by the condition
\begin{align}\label{eq:there-exists-constraint}
    \exists\,\big(\{\bm{a}_k,\bm{A}_k\}_k\big)\in\mathcal{Z}: ~ \sup_{\{\bm{\xi}_k\in\Xi\}_k}\sum_{k=1}^K\big(t_k\ell(\bm{\xi}_k)-\lambda_k(\bm{\xi}_k)+\bm{a}_k^\top\bm{\xi}_k+\langle\bm{A}_k,\bm{\xi}_k\bm{\xi}_k^\top\rangle\big)\leq 0.
\end{align}
Let~$v_\mathrm{exist}$ denote the optimal value of the optimization problem obtained by replacing the original inequality constraint in~\eqref{eq:inf_constraint} with~\eqref{eq:there-exists-constraint}. We claim that~$v_\mathrm{inf}=v_\mathrm{exist}$. Since~\eqref{eq:there-exists-constraint} is a restriction of the inequality constraint in~\eqref{eq:inf_constraint} (the infimum is bounded above by any single feasible value), we have $v_\mathrm{inf}\leq v_\mathrm{exist}$. For the reverse inequality, note that the inequality in~\eqref{eq:inf_constraint} holds if and only if for every $\delta>0$ there exist coefficients $(\{\bm{a}_k, \bm{A}_k\}_k) \in \mathcal{Z}$ such that the supremum in~\eqref{eq:there-exists-constraint} is at most~$\delta$. Given any such $\{\lambda_k\}_k$ and $\delta>0$, shifting each $\lambda_k$ to $\lambda_k+\delta/K$ simultaneously brings the supremum to zero (certifying~\eqref{eq:there-exists-constraint}) and raises the objective by exactly~$\delta$. Hence, we have $v_\mathrm{exist}\leq v_\mathrm{inf}+\delta$, and since $\delta>0$ was arbitrary, $v_\mathrm{exist}\leq v_\mathrm{inf}$. This reasoning shows that replacing the inequality constraint in~\eqref{eq:inf_constraint} with~\eqref{eq:there-exists-constraint} has no impact on the problem's optimal value.

For any fixed $\lambda_k \in C(\Xi)$, $k\in[K]$, and $(\{\bm{a}_k,\bm{A}_k\}_k) \in \mathcal Z$, the objective function in~\eqref{eq:there-exists-constraint} is additively separable across $\bm{\xi}_1,\ldots,\bm{\xi}_K\in\Xi$. Thus, the joint supremum decomposes as follows.
\begin{align*}
    &\sup_{\{\bm{\xi}_k\in\Xi\}_k } \:  \Big\{ \sum_{k=1}^K \big(  t_k \ell(\bm{\xi}_k)  - \lambda_k(\bm{\xi}_k)  + \bm{a}_k^\top \bm{\xi}_k +\langle \bm{A}_k,\bm{\xi}_k \bm{\xi}_k^\top \rangle \big) \Big\} \\
    & \quad = \;
    \sum_{k=1}^K \sup_{\bm{\xi}_k \in \Xi} \: \Big\{ t_k \ell(\bm{\xi}_k) - \lambda_k(\bm{\xi}_k) + \bm{a}_k^\top \bm{\xi}_k + \langle \bm{A}_k,\bm{\xi}_k \bm{\xi}_k^\top \rangle \Big\}.
\end{align*}
Substituting the last expression into the restriction of Problem~\eqref{eq:inf_constraint} with constraint~\eqref{eq:there-exists-constraint} yields
\begin{align*} 
      \begin{array}{cll}  \inf & \displaystyle \sum_{k=1}^K \sup_{\mathbb{Q}_k \in \mathfrak{B}_{k}} \mathbb{E}_{\mathbb{Q}_k}[\lambda_k(\bm{\xi})]  & \\
         \mathrm{s.t.} & \big(\{\bm{a}_k, \bm{A}_k\}_k\big)\in\cZ, \; \lambda_k \in C(\Xi) \; \quad \forall k\in[K]\\
         &  \displaystyle \sum_{k=1}^K \sup_{\bm{\xi}_k \in \Xi} \: \left\{ t_k \ell(\bm{\xi}_k) - \lambda_k(\bm{\xi}_k)  +  {\bm{a}_k}^\top \bm{\xi}_k +\langle \bm{A}_k, \bm{\xi}_k\bm{\xi}_k^\top\rangle\right\}  \leq 0 .
    \end{array}
\end{align*}
Next, we use Lemma~\ref{lemma:decouple} with the translation invariant functionals
\begin{align*}
     S_k(\lambda_k) = \sup_{\bm{\xi}_k \in \Xi} \: \left\{ t_k \ell(\bm{\xi}_k) - \lambda_k(\bm{\xi}_k)  +  \bm{a}_k^\top \bm{\xi}_k +\langle \bm{A}_{k}, \bm{\xi}_k\bm{\xi}_k^\top \rangle\right\} \quad \text{and} \quad T_k(\lambda_k) = \sup_{\mathbb{Q}_k \in \mathfrak{B}_{k}} \mathbb{E}_{\mathbb{Q}_k}[\lambda_k(\bm{\xi})]
\end{align*}
to show that the constraint $\sum_{k=1}^KS_k(\lambda_k)\leq 0$ can be disaggregated into $K$ separate constraints $S_k(\lambda_k)\leq 0$, $k\in[K]$, without changing the optimal value. Lemma~\ref{lemma:decouple} applies because each~$S_k$ is the supremum of a real-valued function over the nonempty set~$\Xi$ and thus never equals~$-\infty$. The above problem thus reduces to
\begin{align*}
    \begin{array}{cl}
        \inf & \displaystyle\sum_{k=1}^K\sup_{\mathbb{Q}_k\in\mathfrak{B}_k}\mathbb{E}_{\mathbb{Q}_k}[\lambda_k(\bm{\xi})] \\
        \mathrm{s.t.} & \big(\{\bm{a}_k,\bm{A}_k\}_k\big)\in\mathcal{Z},\;\lambda_k\in C(\Xi)\quad\forall k\in[K] \\
        & t_k\ell(\bm{\xi}_k)+\bm{a}_k^\top\bm{\xi}_k+\langle\bm{A}_k,\bm{\xi}_k\bm{\xi}_k^\top\rangle\leq\lambda_k(\bm{\xi}_k)\quad\forall\bm{\xi}_k\in\Xi,\;k\in[K].
    \end{array}
\end{align*}
For any fixed $(\{\bm{a}_k,\bm{A}_k\}_k)\in\mathcal{Z}$, the constraint bounds each $\lambda_k$ pointwise from below, and the objective is nondecreasing in each $\lambda_k$ under the pointwise order. The optimal choice is therefore the pointwise smallest admissible function, $\lambda_k(\bm{\xi})=t_k\ell(\bm{\xi})+\bm{a}_k^\top\bm{\xi}+\langle\bm{A}_k,\bm{\xi}\bm{\xi}^\top\rangle$, which belongs to~$C(\Xi)$ by the continuity of~$\ell$. Substituting these optimal functions~$\lambda_k$, $k\in[K]$, into the objective yields~\eqref{prob:affine}. In summary, we have thus shown that the restriction of Problem~\eqref{eq:prob-with-semi-infinite-constraint} to affine decision rules of the form~\eqref{eq:affine-decision-rule} is equivalent to~\eqref{prob:affine}. By Proposition~\ref{prop:reformulation_main} and Lemma~\ref{lem:pairwise-symmetric-rules}, the optimal value of~\eqref{eq:prob-with-semi-infinite-constraint} upper bounds~$v_\mathrm{P}$, and as~\eqref{prob:affine} is a restriction of~\eqref{eq:prob-with-semi-infinite-constraint}, we have $v_\mathrm{P}\leq v_\mathrm{aff}$.
\end{proof}

The proof of Proposition~\ref{prop:affine-tightness} relies on three lemmas. The first of these lemmas establishes that~$v_\mathrm{aff}$ coincides with the optimal value of a tractable semidefinite program.

\begin{lemma}\label{lem:affine-sdp-equivalence}
Suppose the conditions of Proposition~\ref{prop:affine-tightness} hold. Then, 
$v_\mathrm{aff}$ equals
\begin{align}\label{eq:sdp-inf}\tag{SDP-P}
    v_{\mathrm{SDP\text{-}P}}\coloneqq \left\{\begin{array}{cll}
        \inf & r + \tau + \bm{1}^\top\bm{s} \\
        \mathrm{s.t.} & \tau\in\mathbb{R},\;\gamma_k\in\mathbb{R}_+,\;s_k\in\mathbb{R},\;\bm{Q}_k,\bm{Z}_k\in\mathbb{S}^m & \forall k\in[K] \\
        & \sum_{k=1}^K\bm{Q}_k = \Big(\sum_{k=1}^K\gamma_k\Big)\bm{I}_m - \bm{P} \\[0.4em]
        & \begin{bmatrix} \sum_{k=1}^K\bm{Q}_k & \dfrac{1}{2}\Big(\bm{q}+ \sum_{k=1}^K 2\gamma_k\hat{\bm{\mu}}_k\Big) \\[0.6em] \dfrac{1}{2}\Big(\bm{q}+ \sum_{k=1}^K 2\gamma_k\hat{\bm{\mu}}_k\Big)^\top & \tau \end{bmatrix} \succeq 0 \\[1.5em]
        & \begin{bmatrix} \bm{Z}_k & \gamma_k\bm{I}_m \\ \gamma_k\bm{I}_m & \bm{Q}_k \end{bmatrix} \succeq 0 & \forall k\in[K] \\[1.0em]
        & \gamma_k\varepsilon_k^2 + \mathrm{Tr}(\hat{\bm{\Sigma}}_k\bm{Z}_k) - \gamma_k\mathrm{Tr}(\hat{\bm{\Sigma}}_k) - \gamma_k\|\hat{\bm{\mu}}_k\|_2^2 \leq s_k & \forall k\in[K].
    \end{array} \right.
\end{align}
\end{lemma}

\begin{proof}
Substituting the quadratic loss into~\eqref{prob:affine} gives
\begin{align}\label{prob:affine_quad_loss}
    \begin{array}{rlll}
        & \inf & \displaystyle\sum_{k=1}^K\sup_{\mathbb{Q}_k\in\mathfrak{B}_k}
            \mathbb{E}_{\mathbb{Q}_k}\big[t_k(\bm{\xi}^\top\bm{P}\bm{\xi}
            +\bm{q}^\top\bm{\xi}+r)+\bm{a}_k^\top\bm{\xi}
            +\langle\bm{A}_k,\bm{\xi}\bm{\xi}^\top\rangle\big] \\
        & \mathrm{s.t.} & \bm{a}_k\in\mathbb{R}^m,\;\bm{A}_k\in\mathbb{S}^m
            \quad\forall k\in[K] \\
        & & \sum_{k=1}^K\bm{A}_k=\bm{0},\quad\sum_{k=1}^K\bm{a}_k=\bm{0}.
    \end{array}
\end{align}
Since $\Xi=\mathbb{R}^m$, the loss is quadratic, $\hat{\mathbb P}_k$ has finite second moments and~$\varepsilon_k>0$, Wasserstein strong duality applies to each inner supremum \citep[theorem~4.18]{kuhn2025distributionally}. Thus, we find
\begin{align*}
    &\sup_{\mathbb{Q}_k\in\mathfrak{B}_k}\mathbb{E}_{\mathbb{Q}_k}
        \big[t_k(\bm{\xi}^\top\bm{P}\bm{\xi}+\bm{q}^\top\bm{\xi}+r)
        +\bm{a}_k^\top\bm{\xi}+\langle\bm{A}_k,\bm{\xi}\bm{\xi}^\top\rangle\big] \\
    =\;&\inf_{\gamma_k\geq 0}\gamma_k\varepsilon_k^2+t_kr
        +\mathbb{E}_{\hat{\mathbb{P}}_k}\Big[\sup_{\bm{\xi}\in\mathbb{R}^m}
        \big\{\bm{\xi}^\top(t_k\bm{P}+\bm{A}_k)\bm{\xi}
        +(t_k\bm{q}+\bm{a}_k)^\top\bm{\xi}
        -\gamma_k\|\bm{\xi}-\hat{\bm{\xi}}\|_2^2\big\}\Big].
\end{align*}
Setting $\bm{Q}_k\coloneqq\gamma_k\bm{I}_m-t_k\bm{P}-\bm{A}_k$ and $\bm{b}_k\coloneqq t_k\bm{q}+\bm{a}_k$, the inner supremum over $\bm{\xi}$ equals $+\infty$ unless $\bm{Q}_k\succeq \bm{0}$, in which case it equals $\frac{1}{4}(\bm{b}_k+2\gamma_k\hat{\bm{\xi}})^\top\bm{Q}_k^{-1} (\bm{b}_k+2\gamma_k\hat{\bm{\xi}}) -\gamma_k \|\hat{\bm{\xi}}\|_2^2$. Here, $\bm{Q}_k^{-1}$ denotes the Moore-Penrose pseudoinverse, and we adopt the convention that $\bm{c}^\top\bm{Q}_k^{-1} \bm{c}=+\infty$ whenever $\bm{c}\notin\mathrm{range}(\bm{Q}_k^{-1})$. If $\gamma_k>0$ and $\bm{Q}_k$ is singular, then $\bm{b}_k+2\gamma_k\hat{\bm{\xi}}\notin\mathrm{range}(\bm{Q}_k)$ with positive probability under~$\hat{\mathbb{P}}_k$ because $\hat{\bm{\Sigma}}_k\succ\bm{0}$, and thus the expectation of the inner supremum is $+\infty$. Hence, we may impose $\bm{Q}_k\succ\bm{0}$ whenever $\gamma_k>0$ without changing the infimum in~\eqref{prob:affine_quad_loss}. Under this condition, taking the expectation with respect to the reference distribution $\hat{\mathbb{P}}_k$ and expressing it in terms of the mean vector $\hat{\bm{\mu}}_k$ and covariance matrix $\hat{\bm{\Sigma}}_k$ of $\hat{\mathbb{P}}_k$ yields
\begin{align*}
    &\mathbb{E}_{\hat{\mathbb{P}}_k}\Big[\frac{1}{4}
        (\bm{b}_k+2\gamma_k\hat{\bm{\xi}})^\top\bm{Q}_k^{-1}
        (\bm{b}_k+2\gamma_k\hat{\bm{\xi}})\Big]
        -\gamma_k\mathbb{E}_{\hat{\mathbb{P}}_k}[\|\hat{\bm{\xi}}\|_2^2] \\
    =\;&\frac{1}{4}\bm{b}_k^\top\bm{Q}_k^{-1}\bm{b}_k
        +\gamma_k\bm{b}_k^\top\bm{Q}_k^{-1}\hat{\bm{\mu}}_k
        +\gamma_k\langle\hat{\bm{\Sigma}}_k
        +\hat{\bm{\mu}}_k\hat{\bm{\mu}}_k^\top,
        \gamma_k\bm{Q}_k^{-1}-\bm{I}_m\rangle.
\end{align*}
Substituting back into~\eqref{prob:affine_quad_loss} and applying the change of variable $\bm{b}_k\leftarrow\bm{b}_k+2\gamma_k\hat{\bm{\mu}}_k$ gives
\begin{align}\label{prob:affine_quad_loss_reformulated}
    \begin{array}{cll}
        \inf & \displaystyle r+\sum_{k=1}^K\Big(\gamma_k\varepsilon_k^2
            +\frac{1}{4}\bm{b}_k^\top\bm{Q}_k^{-1}\bm{b}_k
            +\gamma_k\langle\hat{\bm{\Sigma}}_k,
            \gamma_k\bm{Q}_k^{-1}-\bm{I}_m\rangle
            -\gamma_k\|\hat{\bm{\mu}}_k\|_2^2\Big) \\
        \mathrm{s.t.} & \bm{a}_k,\bm{b}_k\in\mathbb{R}^m,\;
            \bm{A}_k\in\mathbb{S}^m,\;(\gamma_k,\bm{Q}_k)\in\big(\{0\}\times\mathbb{S}_+^m\big)\cup\big(\mathbb{R}_{++}\times\mathbb{S}_{++}^m\big) & \forall k\in[K] \\
        & \bm{Q}_k=\gamma_k\bm{I}_m-t_k\bm{P}-\bm{A}_k,\quad
            \bm{b}_k=t_k\bm{q}+\bm{a}_k+2\gamma_k\hat{\bm{\mu}}_k
            & \forall k\in[K] \\
        & \sum_{k=1}^K\bm{A}_k=\bm{0},\quad\sum_{k=1}^K\bm{a}_k=\bm{0}.
    \end{array}
\end{align}
The domain of $(\gamma_k,\bm{Q}_k)$ in~\eqref{prob:affine_quad_loss_reformulated} captures the requirement $\bm{Q}_k\succ\bm{0}$ whenever $\gamma_k>0$ imposed above. Fix now any $\{\gamma_k,\bm{Q}_k\}_k$ feasible in~\eqref{prob:affine_quad_loss_reformulated}. The minimization over $\{\bm{a}_k,\bm{b}_k\}_k$ reduces to
\begin{align*}
    \inf_{\{\bm{b}_k\}_k}\bigg\{\frac{1}{4}\sum_{k=1}^K\bm{b}_k^\top
        \bm{Q}_k^{-1}\bm{b}_k :
        \sum_{k=1}^K\bm{b}_k=\bm{q}+\sum_{k=1}^K 2\gamma_k\hat{\bm{\mu}}_k\bigg\},
\end{align*}
where we expressed $\bm a_k$ in terms of~$\bm b_k$ and used $\sum_{k=1}^K t_k=1$ to simplify the constraint $\sum_{k=1}^K\bm{a}_k=\bm{0}$. The objective function of the resulting optimization problem is convex in $\{\bm{b}_k\}_k$ and the constraint is affine. By our conventions for the pseudoinverse, the objective is finite only if $\bm{b}_k\in\mathrm{range}(\bm{Q}_k)$ for all $k\in[K]$. Assume for now that $\bm{q}+\sum_{k=1}^K 2\gamma_k\hat{\bm{\mu}}_k=(\sum_{k=1}^K\bm{Q}_k)\bm{z}$ for some~$\bm{z}\in\mathbb{R}^m$, which we later show to be without loss of generality. Then $\bm{b}_k=\bm{Q}_k\bm{z}$, $k\in[K]$, is a Slater point, and strong duality holds. Introducing a multiplier $\bm{x}\in\mathbb{R}^m$ for the equality constraint, the minimization over $\{\bm{b}_k\}_k$ admits the closed-form solution $\bm{b}_k=2\bm{Q}_k\bm{x}$, giving
\begin{align*}
    & \max_{\bm{x}\in\mathbb{R}^m}\left\{\bm{x}^\top\left(\bm{q}
        +\sum_{k=1}^K 2\gamma_k\hat{\bm{\mu}}_k\right)
        -\bm{x}^\top\Big(\sum_{k=1}^K\bm{Q}_k\Big)\bm{x}\right\} \\
    & \quad =\frac{1}{4}\left(\bm{q}+\sum_{k=1}^K 2\gamma_k\hat{\bm{\mu}}_k\right)^\top
        \left(\sum_{k=1}^K\bm{Q}_k \right)^{-1}
        \left(\bm{q}+\sum_{k=1}^K 2\gamma_k\hat{\bm{\mu}}_k\right),
\end{align*}
where $(\sum_{k=1}^K\bm{Q}_k)^{-1}$ denotes again the Moore-Penrose pseudoinverse with the same conventions as above. The assumption that $\bm{q}+\sum_{k=1}^K 2\gamma_k\hat{\bm{\mu}}_k=(\sum_{k=1}^K\bm{Q}_k)\bm{z}$ for some $\bm{z}\in\mathbb{R}^m$ is without loss of generality. Indeed, otherwise this vector lies outside $\sum_{k=1}^K\mathrm{range}(\bm{Q}_k)=\mathrm{range}(\sum_{k=1}^K\bm{Q}_k)$, so every feasible $\{\bm{b}_k\}_k$ has an infinite objective value, and both sides of the identity equal $+\infty$ as well. Summing $\bm{Q}_k=\gamma_k\bm{I}_m-t_k\bm{P}-\bm{A}_k$ over $k$ and using $\sum_{k=1}^K t_k=1$ and $\sum_{k=1}^K \bm{A}_k=\bm{0}$ gives $\sum_{k=1}^K\bm{Q}_k=(\sum_{k=1}^K\gamma_k)\bm{I}_m-\bm{P}$. This shows that~$\bm{A}_k$ is irrelevant and can be eliminated.

After eliminating $\bm{a}_k$, $\bm b_k$ and $\bm{A}_k$, Problem~\eqref{prob:affine_quad_loss_reformulated} involves the two nonlinear terms $(\sum_{k=1}^K \bm{Q}_k)^{-1}$ and $\gamma_k^2\bm{Q}_k^{-1}$. Introducing an epigraphical variable $\tau$ for the former and auxiliary matrix variables~$\bm{Z}_k$ for the latter, both inverse relations can be linearized using Schur complements. We also use epigraphical variables $s_k$ to absorb the scalar terms $\gamma_k\varepsilon_k^2+\langle\hat{\bm{\Sigma}}_k,\bm{Z}_k\rangle-\gamma_k\langle\hat{\bm{\Sigma}}_k,\bm{I}_m\rangle-\gamma_k\|\hat{\bm{\mu}}_k\|_2^2\leq s_k$. Combining these substitutions yields~\eqref{eq:sdp-inf}, and thus the claim follows.
\end{proof}

The second lemma establishes a semidefinite lower bound on~$v_\mathrm{P}$. 

\begin{lemma}\label{lem:rho-star-sdp-lowerbnd} Suppose the conditions of Proposition~\ref{prop:affine-tightness} hold. Then, 
$v_\mathrm{P}$ is bounded below by
\begin{equation}\label{eq:sdp-sup}\tag{SDP-D}
        v_{\mathrm{SDP\text{-}D}}\coloneqq \left\{\begin{array}{cll} \sup &   \langle \bm{P}, \bm{M}\rangle + \bm{q}^\top \bm{\mu}+r    \\[1ex]
        \mathrm{s.t.} & \bm{\mu}\in\mathbb{R}^m, \;\bm{M},\bm{H}\in\mathbb{S}^m, \; \bm{C}_k \in \mathbb{R}^{m \times m} & \forall k\in[K]\\[1ex]
        & \|\hat{\bm{\mu}}_k\|_2^2-2\hat{\bm{\mu}}_k^\top \bm{\mu} + \mathrm{Tr}(\bm{M}+\widehat{\bm{\Sigma}}_k -2\bm{C}_k)\leq \varepsilon_k^2 & \forall k\in[K]\\[1ex]
        & \begin{bmatrix}
            \bm{H} & \bm{C}_k\\ \bm{C}_k^\top & \widehat{\bm{\Sigma}}_k
        \end{bmatrix} \succeq 0 & \forall k\in[K]\\[3ex]
        & \begin{bmatrix}
            \bm{M}-\bm{H} & \bm{\mu}\\ \bm{\mu}^\top & 1
        \end{bmatrix} \succeq 0 .
    \end{array} \right.
    \end{equation}
\end{lemma}

\begin{proof} If~\eqref{eq:sdp-sup} is infeasible, then $v_{\mathrm{SDP\text{-}D}}=-\infty\leq v_\mathrm{P}$ holds trivially. Otherwise, fix any feasible solution $(\bm{\mu},\bm{M},\bm{H},\{\bm{C}_k\}_k)$ of Problem~\eqref{eq:sdp-sup}, and set $\bm{\Sigma}\coloneqq\bm{M}-\bm{\mu}\bm{\mu}^\top$. The last block constraint gives $\bm{M}-\bm{H}\succeq\bm{\mu}\bm{\mu}^\top$ by the Schur complement lemma, and the first block constraint gives $\bm{H}\succeq\bm{0}$. We may thus conclude that $\bm{\Sigma}=\bm{M}-\bm{\mu}\bm{\mu}^\top\succeq\bm{H}\succeq\bm{0}$. Consequently, $(\bm{\mu},\bm{\Sigma})$ is a valid mean-covariance pair satisfying $\bm{\Sigma}\succeq\bm{H}$.

For each $k\in[K]$, define $\phi_k(\bm{H})\coloneqq\mathrm{Tr}[(\hat{\bm{\Sigma}}_k^{1/2}\bm{H}\hat{\bm{\Sigma}}_k^{1/2})^{1/2}]$. By the variational representation of the square-root matrix trace~\citep[proposition~2.2]{kuhn2025distributionally}, we have
\[
    \phi_k(\bm{H})=\max_{\bm{C} \in\mathbb R^{m\times m}}\left\{\mathrm{Tr}(\bm{C}):\begin{bmatrix}\bm{H}&\bm{C}\\\bm{C}^\top&\hat{\bm{\Sigma}}_k\end{bmatrix}\succeq\bm{0}\right\}.
\]
Since $\bm{C}_k$ is feasible in this maximization problem, we have $\mathrm{Tr}(\bm{C}_k)\leq\phi_k(\bm{H})$. In addition, since $\phi_k$ is monotone with respect to the positive-semidefinite order and $\bm{\Sigma}\succeq\bm{H}$, we further have $\mathrm{Tr}(\bm{C}_k)\leq\phi_k(\bm{H})\leq\phi_k(\bm{\Sigma})$. Expanding $\|\bm{\mu}-\hat{\bm{\mu}}_k\|_2^2=\|\hat{\bm{\mu}}_k\|_2^2-2\hat{\bm{\mu}}_k^\top\bm{\mu}+\|\bm{\mu}\|_2^2$ and using $\mathrm{Tr}(\bm{\Sigma})=\mathrm{Tr}(\bm{M})-\|\bm{\mu}\|_2^2$, the scalar constraint of~\eqref{eq:sdp-sup} and the inequality above together give
\begin{align*}
    \mathsf{G}^2\big((\bm{\mu},\bm{\Sigma}),(\hat{\bm{\mu}}_k,\hat{\bm{\Sigma}}_k)\big) &=\|\bm{\mu}-\hat{\bm{\mu}}_k\|_2^2+\mathrm{Tr}(\bm{\Sigma}+\hat{\bm{\Sigma}}_k)-2\phi_k(\bm{\Sigma}) \\
    &\leq\|\hat{\bm{\mu}}_k\|_2^2-2\hat{\bm{\mu}}_k^\top\bm{\mu}+\mathrm{Tr}(\bm{M}+\hat{\bm{\Sigma}}_k-2\bm{C}_k)\leq\varepsilon_k^2 \qquad\forall k\in[K].
\end{align*}
Let $\mathcal{S}$ be the positive-semidefinite affine family containing $\hat{\mathbb{P}}_1,\ldots,\hat{\mathbb{P}}_K$. Since $\hat{\bm{\Sigma}}_k\succ\bm{0}$, the matrix
\[
	\bm{B}_k\coloneqq\hat{\bm{\Sigma}}_k^{-1/2}(\hat{\bm{\Sigma}}_k^{1/2}\bm{\Sigma}\hat{\bm{\Sigma}}_k^{1/2})^{1/2}\hat{\bm{\Sigma}}_k^{-1/2}\in\mathbb{S}_+^m
\]
satisfies $\bm{B}_k\hat{\bm{\Sigma}}_k\bm{B}_k=\bm{\Sigma}$. Therefore, the affine map $T_k(\bm{\xi})\coloneqq\bm{B}_k\bm{\xi}+\bm{\mu}-\bm{B}_k\hat{\bm{\mu}}_k$ pushes $\hat{\mathbb{P}}_k$ forward to a distribution with moments $(\bm{\mu},\bm{\Sigma})$, which lies in $\mathcal{S}$ by Condition~\ref{def:psd-affine-closed} of Definition~\ref{def:psd-affine-family}. By Condition~\ref{def:psd-affine-unique}, $\mathcal{S}$ contains exactly one distribution with moments $(\bm{\mu},\bm{\Sigma})$, and thus all $K$ pushforwards coincide. From now on we denote this common distribution by $\mathbb{Q}$. As the Gelbrich bound is exact for positive-semidefinite affine pushforwards~\citep[theorem~2]{nguyen2021mean}, we have
\[
    \mathsf{W}_2^2(\mathbb{Q},\hat{\mathbb{P}}_k)=\mathsf{G}^2\big((\bm{\mu},\bm{\Sigma}),(\hat{\bm{\mu}}_k,\hat{\bm{\Sigma}}_k)\big)\leq\varepsilon_k^2 \qquad\forall k\in[K],
\]
which shows that $\mathbb{Q}$ resides within the fusion ambiguity set~$\mathcal{P}$. Recall now that $\bm{\mu}=\mathbb{E}_{\mathbb{Q}}[\bm{\xi}]$ and $\bm{M}=\mathbb{E}_{\mathbb{Q}} [\bm{\xi}\bm{\xi}^\top]$ by construction, and that $\ell(\bm{\xi})=\bm{\xi}^\top\bm{P}\bm{\xi} +\bm{q}^\top\bm{\xi}+r$ is quadratic. Therefore, we find
\[
    \mathbb{E}_{\mathbb{Q}}[\ell(\bm{\xi})]=\langle\bm{P},\bm{M}\rangle
    +\bm{q}^\top\bm{\mu}+r,
\]
which equals the objective of~\eqref{eq:sdp-sup} at the given feasible solution. Since $\mathbb{Q}\in\mathcal{P}$ and~\eqref{prob:sup-primal} is a maximization problem, this value is at most $v_\mathrm{P}$. As the feasible solution of~\eqref{eq:sdp-sup} was arbitrary, we have $v_{\mathrm{SDP\text{-}D}}\leq v_\mathrm{P}$, which concludes the proof.
\end{proof}

The third lemma shows that the semidefinite programs from Lemmas~\ref{lem:affine-sdp-equivalence} and~\ref{lem:rho-star-sdp-lowerbnd} are strong~duals.

\begin{lemma}\label{lem:sdp-duality}
    Problems~\eqref{eq:sdp-inf} and~\eqref{eq:sdp-sup} are strong duals, and their optimal values coincide.
\end{lemma}

\begin{proof}
We associate Lagrange multipliers with all constraints of~\eqref{eq:sdp-inf}. The equality constraint is assigned a matrix $\bm{M}\in\mathbb{S}^m$, and the $K$ scalar inequalities are assigned scalars $\alpha_k\geq 0$. In addition, the linear matrix inequalities are assigned positive-semidefinite block matrices
\[
    \begin{bmatrix}\bm{U}&-\bm{\mu}\\-\bm{\mu}^\top&\lambda\end{bmatrix}\succeq\bm{0} \quad \text{and} \quad \begin{bmatrix}\bm{A}_k&-\bm{C}_k^\top\\-\bm{C}_k&\bm{H}_k\end{bmatrix}\succeq\bm{0} \quad\forall k\in[K],
\]
respectively. Letting $\bm{X}\coloneqq(\tau,\bm{s},\bm{\gamma},\{\bm{Q}_k,\bm{Z}_k\}_k)$ and $\bm{Y}\coloneqq(\bm{M},\bm{U},\bm{\mu},\lambda,\{\bm{A}_k,\bm{C}_k,\bm{H}_k,\alpha_k\}_k)$ collect the primal and dual variables, the Lagrangian can be represented as
\begin{align*}
        L(\bm{X},\bm{Y}) &= r + \langle \bm{M}, \bm{P} \rangle + \bm{q}^\top \bm{\mu} + \tau(1 - \lambda) + \sum_{k=1}^K s_k(1 - \alpha_k) + \sum_{k=1}^K \langle \bm{Z}_k, \alpha_k \hat{\bm{\Sigma}}_k - \bm{A}_k \rangle \\
        &\quad + \sum_{k=1}^K \langle \bm{Q}_k, \bm{M} - \bm{U} - \bm{H}_k \rangle \\
        &\quad + \sum_{k=1}^K \gamma_k \Big( \alpha_k\varepsilon_k^2 - \alpha_k\|\hat{\bm{\mu}}_k\|_2^2 + 2\hat{\bm{\mu}}_k^\top \bm{\mu} - \mathrm{Tr}(\bm{M}) - \alpha_k\mathrm{Tr}(\hat{\bm{\Sigma}}_k) + 2\mathrm{Tr}(\bm{C}_k) \Big).
     \end{align*}
The dual objective $\inf_{\bm{X}}L(\bm{X},\bm{Y})$ evaluates to~$-\infty$ unless the coefficient of every free primal variable vanishes and the coefficient of every nonnegative primal variable is nonnegative. Stationarity with respect to $\tau$ requires $\lambda=1$, and stationarity with respect to $s_k$ requires $\alpha_k=1$. Stationarity with respect to $\bm{Z}_k\in\bS^m$ requires $\bm{A}_k=\hat{\bm{\Sigma}}_k$. Stationarity with respect to $\gamma_k\geq 0$ requires $\varepsilon_k^2-\|\hat{\bm{\mu}}_k\|_2^2+2\hat{\bm{\mu}}_k^\top\bm{\mu}-\mathrm{Tr}(\bm{M})-\mathrm{Tr}(\hat{\bm{\Sigma}}_k)+2\mathrm{Tr}(\bm{C}_k)\geq 0$, and stationarity with respect to $\bm{Q}_k\in\mathbb{S}^m$ requires $\bm{H}_k=\bm{M}-\bm{U}$. As the right-hand side is independent of $k$, we write $\bm{H}\coloneqq\bm{H}_k$ for all $k \in[K]$. The dual objective then reduces to $\langle\bm{P},\bm{M}\rangle+\bm{q}^\top\bm{\mu}+r$. Maximizing over the remaining dual variables and using the congruence identity
\[
    \begin{bmatrix}\bm{A}&-\bm{B}\\-\bm{B}^\top&\bm{C}\end{bmatrix}\succeq\bm{0} \quad \iff \quad \begin{bmatrix}\bm{A}&\bm{B}\\\bm{B}^\top&\bm{C}\end{bmatrix}\succeq\bm{0}
\]
for any matrices $\bm{A}$, $\bm{B}$ and $\bm{C}$ of compatible dimensions yields~\eqref{eq:sdp-sup}. Hence, \eqref{eq:sdp-sup} is the Lagrangian dual of~\eqref{eq:sdp-inf}. A Slater point for~\eqref{eq:sdp-inf} is obtained by choosing any $\bm Q_k\succ \bm 0$ and $\gamma_k>0$ satisfying the equality constraint and selecting $\tau>0$, $\bm{Z}_k\succ\bm{0}$ and $s_k$ sufficiently large to strictly satisfy all inequalities. Thus, strong duality holds, and $v_{\mathrm{SDP-P}} = v_{\mathrm{SDP-D}}$.
\end{proof}

With the three lemmas in place, we can now prove Proposition~\ref{prop:affine-tightness}. 

\begin{proof}[Proof of Proposition~\ref{prop:affine-tightness}]
The results obtained so far imply that
\[
    v_\mathrm{P} \leq v_\mathrm{aff} = v_{\mathrm{SDP\text{-}P}} = v_{{\mathrm{SDP\text{-}D}}} \leq v_\mathrm{P}
\]
where the two inequalities follow from Proposition~\ref{prop:affine_upper_bound} and from Lemma~\ref{lem:rho-star-sdp-lowerbnd}, while the two equalities hold thanks to Lemma~\ref{lem:affine-sdp-equivalence} and Lemma~\ref{lem:sdp-duality}. Thus, both inequalities collapse to equalities, which confirms that the affine decision rule approximation is exact and that~$v_\mathrm{P}$ can be computed by solving the semidefinite program~\eqref{eq:sdp-inf} or its dual~\eqref{eq:sdp-sup}.
\end{proof}

\begin{proof}[Proof of Proposition~\ref{prop:static-tightness}]
We proceed in two steps. We first reduce \eqref{prob:sup-primal} and \eqref{prob:static_upper_bound} to second-order cone programs, and then we show that they are strong duals. Throughout, $\hat{\bm{\Sigma}}$ denotes the common covariance matrix of $\hat{\mathbb{P}}_1,\ldots,\hat{\mathbb{P}}_K$. Since the loss is affine, $\mathbb{E}_{\mathbb{Q}}[\ell(\bm{\xi})]=\bm{q}^\top\bm{\mu}+r$ for every $\mathbb{Q}\in\mathcal{P}_0(\mathbb{R}^m)$ with mean $\bm{\mu}\coloneqq\mathbb{E}_{\mathbb{Q}}[\bm{\xi}]$. Thus, the objective depends on $\mathbb{Q}$ only through its mean.

In order to reduce~\eqref{prob:sup-primal} and \eqref{prob:static_upper_bound} to second-order cone programs, we first show that a vector $\bm{\mu}\in\mathbb{R}^m$ is the mean of some $\mathbb{Q}\in\mathfrak{B}_k$ if and only if $\|\bm{\mu}-\hat{\bm{\mu}}_k\|_2\leq\varepsilon_k$. For the necessity, consider any~$\mathbb{Q}\in\mathfrak{B}_k$ with mean $\bm{\mu}$ and let $\Pi$ be any coupling of $\mathbb{Q}$ and $\hat{\mathbb{P}}_k$. Then, we have
\begin{align*}
    \|\bm{\mu}-\hat{\bm{\mu}}_k\|_2=\big\|\mathbb{E}_\Pi[\bm{\xi}-\bm{\xi}']\big\|_2\leq\mathbb{E}_\Pi\big[\|\bm{\xi}-\bm{\xi}'\|_2\big]\leq\Big(\mathbb{E}_\Pi\big[\|\bm{\xi}-\bm{\xi}'\|_2^2\big]\Big)^{1/2},
\end{align*}
where both inequalities follow from Jensen's inequality. Minimizing the right-hand side over all couplings gives $\|\bm{\mu}-\hat{\bm{\mu}}_k\|_2\leq\mathsf{W}_2(\mathbb{Q},\hat{\mathbb{P}}_k)\leq\varepsilon_k$. For the sufficiency, suppose that $\|\bm{\mu}-\hat{\bm{\mu}}_k\|_2\leq\varepsilon_k$, and define $T_k(\bm{\xi})\coloneqq\bm{\xi}+\bm{\mu}-\hat{\bm{\mu}}_k$. As the Jacobian of $T_k$ is the identity matrix, the pushforward distribution $(T_k)_\#\hat{\mathbb{P}}_k$ belongs to the positive-semidefinite affine family~$\mathcal S$ from Condition~\ref{cond:static-family} and has the same moments $(\bm{\mu},\hat{\bm{\Sigma}})$ for every~$k\in[K]$. By Definition~\ref{def:psd-affine-family}\ref{def:psd-affine-unique}, $\mathcal S$ contains exactly one distribution, say~$\mathbb{Q}$, with this moment pair. This implies that $\mathbb{Q}=(T_k)_\#\hat{\mathbb{P}}_k$ for every $k\in[K]$. As~$T_k$ displaces every point by the constant vector $\bm{\mu}-\hat{\bm{\mu}}_k$, we find $\mathsf{W}_2(\mathbb{Q},\hat{\mathbb{P}}_k)\leq\|\bm{\mu}-\hat{\bm{\mu}}_k\|_2\leq\varepsilon_k$, and hence $\mathbb{Q}\in\mathfrak{B}_k$. Since the same $\mathbb{Q}$ satisfies this for every $k\in[K]$, we have $\mathbb{Q}\in\mathcal{P}$. 

Using the above reasoning to express both the expected loss with respect to~$\mathbb Q$ as well as the conditions~$\mathbb Q\in\mathfrak B_k$, $k\in[K]$, in terms of the mean vector~$\bm \mu$ of~$\mathbb Q$, we can recast Problem~\eqref{prob:sup-primal} as
\begin{align}
    \label{eq:v_P-reformulation}
    v_\mathrm{P}=r+\sup_{\bm{\mu}\in\mathbb{R}^m}\Big\{\bm{q}^\top\bm{\mu}:\|\bm{\mu}-\hat{\bm{\mu}}_k\|_2\leq\varepsilon_k\quad\forall k\in[K]\Big\}.
\end{align}
Similarly, setting $\bm{z}_k\coloneqq t_k\bm{q}+\bm{a}_k$, the $k$-th worst-case expectation in~\eqref{prob:static_upper_bound} reduces to
\begin{align*}
    \sup_{\mathbb{Q}_k\in\mathfrak{B}_k}\mathbb{E}_{\mathbb{Q}_k}\big[t_k\ell(\bm{\xi})+\bm{a}_k^\top\bm{\xi}\big] =t_kr+\sup_{\bm{\mu}\in\mathbb{R}^m}\Big\{\bm{z}_k^\top\bm{\mu}:\|\bm{\mu}-\hat{\bm{\mu}}_k\|_2\leq\varepsilon_k\Big\} =t_kr+\hat{\bm{\mu}}_k^\top\bm{z}_k+\varepsilon_k\|\bm{z}_k\|_2,
\end{align*}
where the second equality is obtained by solving the maximization problem over~$\bm \mu$ analytically. Summing both sides of the above equation over $k\in[K]$, recalling that $\sum_{k=1}^K t_k=1$, and noting that the constraint $\sum_{k=1}^K\bm{a}_k=\bm{0}$ holds if and only if $\sum_{k=1}^K\bm{z}_k=\bm{q}$, \eqref{prob:static_upper_bound} then simplifies to
\begin{align}
    \label{eq:v_stat-reformulation}
    v_\mathrm{stat}=r+\inf_{\{\bm{z}_k\}_k} \left\{\sum_{k=1}^K\big(\hat{\bm{\mu}}_k^\top\bm{z}_k+\varepsilon_k\|\bm{z}_k\|_2\big):\sum_{k=1}^K\bm{z}_k=\bm{q}\right\}.
\end{align}

It remains to be shown that~\eqref{eq:v_P-reformulation} and~\eqref{eq:v_stat-reformulation} are strong duals. Denoting the Lagrange multiplier of the equality constraint by~$\bm{\mu}\in\mathbb{R}^m$, we can represent the problem dual to~\eqref{eq:v_stat-reformulation} as follows.
\[
    \sup_{\bm{\mu}\in\mathbb{R}^m} \left\{ r+\bm{q}^\top\bm{\mu}+ \sum_{k=1}^K \inf_{\bm{z}_k \in\mathbb R^m} \left\{(\hat{\bm{\mu}}_k-\bm{\mu})^\top\bm{z}_k+\varepsilon_k\|\bm{z}_k\|_2\right\} \right\}.
\]
Here, the inner infimum over $\bm{z}_k$ equals $0$ if $\|\bm{\mu}-\hat{\bm{\mu}}_k\|_2\leq\varepsilon_k$ (attained at $\bm{z}_k=\bm{0}$) and $-\infty$ otherwise. Hence, the objective function of the dual problem is finite if and only if $\|\bm{\mu}-\hat{\bm{\mu}}_k\|_2\leq\varepsilon_k$ for all $k\in[K]$, in which case it equals $r+\bm{q}^\top\bm{\mu}$. This confirms that~\eqref{eq:v_stat-reformulation} is indeed dual to~\eqref{eq:v_P-reformulation}. A Slater point for~\eqref{eq:v_stat-reformulation} is given by $\bm{z}_k=\bm{q}/K$ for every $k\in[K]$, and thus strong duality holds.

Taken together, the above reasoning implies that~$v_\mathrm{P}=v_\mathrm{stat}$ and that~\eqref{prob:static_upper_bound} admits a tractable second-order cone reformulation. Thus, the claim follows.
\end{proof}


\begin{proof}[Proof of Proposition~\ref{prop:mink}]
We first show that~$\mathcal{P}\subseteq\mathcal{P}_\mathrm{mink}$ and $v_\mathrm{P}\leq v_\mathrm{mink}$. For any $\mathbb{Q}\in\mathcal{P}=\cap_{k=1}^K\mathfrak{B}_k$, setting $\mathbb{Q}_k=\mathbb{Q}$ for every~$k\in[K]$ gives $\sum_k t_k\mathbb{Q}_k=\mathbb{Q}$ with $\mathbb{Q}_k\in\mathfrak{B}_k$ for every~$k\in[K]$, which implies~$\mathbb{Q}\in\mathcal{P}_\mathrm{mink}$. As the choice of $\mathbb Q\in\mathcal P$ was arbitrary, we have $\mathcal{P}\subseteq\mathcal{P}_\mathrm{mink}$ and thus $v_\mathrm{P}\leq v_\mathrm{mink}$. 

Next, we prove that $v_\mathrm{stat}\leq v_\mathrm{mink}$. Setting $\bm{a}_k=\bm{0}$ for all $k\in[K]$ is feasible in~\eqref{prob:static_upper_bound}, and the objective value of this feasible solution amounts to $\sum_k t_k\sup_{\mathbb{Q}_k\in\mathfrak{B}_k} \mathbb{E}_{\mathbb{Q}_k}[\ell(\bm{\xi})]=v_\mathrm{mink}$. As~\eqref{prob:static_upper_bound}  minimizes over $\{\bm{a}_k\}_k$, its optimal value cannot exceed~$v_\mathrm{mink}$. Thus, we have $v_\mathrm{stat}\leq v_\mathrm{mink}$.

It remains to be shown that~$v_\mathrm{P}\leq v_\mathrm{stat}$. This inequality is implied by Corollary~\ref{cor:static_upper_bound} when~$\ell$  is continuous, but we provide here an alternative direct argument that requires no regularity on~$\ell$. Let $\mathbb{Q}\in\mathcal{P}$ and let $\{\bm{a}_k\}_k$ be feasible in~\eqref{prob:static_upper_bound}. Since $\mathbb{Q}\in\mathfrak{B}_k$ for every~$k\in[K]$, we obtain
\begin{align*}
    \sum_{k=1}^K\sup_{\mathbb{Q}_k\in\mathfrak{B}_k}\mathbb{E}_{\mathbb{Q}_k}
    \big[t_k\ell(\bm{\xi})+\bm{a}_k^\top\bm{\xi}\big]
    \geq\sum_{k=1}^K\mathbb{E}_{\mathbb{Q}}\big[t_k\ell(\bm{\xi})
    +\bm{a}_k^\top\bm{\xi}\big]\;\mathrel{\geq}\mathbb{E}_{\mathbb{Q}}[\ell(\bm{\xi})],
\end{align*}
where the second inequality uses the subadditivity of the expectation established in the proof of Proposition~\ref{prop:reformulation_main} together with the identities $\sum_{k=1}^K t_k=1$ and $\sum_{k=1}^K \bm{a}_k=\bm{0}$. Taking the supremum over $\mathbb{Q}\in\mathcal{P}$ on the right and the infimum over feasible $\{\bm{a}_k\}_k$ on the left yields $v_\mathrm{P}\leq v_\mathrm{stat}$.
\end{proof}

\begin{proof}[Proof of Proposition~\ref{prop:mix}]
To show $\mathcal{P}_\mathrm{mink}\subseteq\mathcal{P}_\mathrm{mix}$, select any $\mathbb{Q} \in\mathcal{P}_\mathrm{mink}$, and recall from Section~\ref{sec:minkowski} that $\mathbb{Q}=\sum_{k=1}^K t_k\mathbb{Q}_k$ for some $\mathbb{Q}_k\in\mathfrak{B}_{\varepsilon_k}^{d_k}(\hat{\mathbb{P}}_k)$, $k\in[K]$. By Assumption~\ref{assump:aggregating-discrepancy}, we then have
\begin{align*}
    d(\mathbb{Q},\mathbb{P}_\mathrm{mix})\leq\sum_{k=1}^K t_k\, d_k(\mathbb{Q}_k,\hat{\mathbb{P}}_k)\leq\sum_{k=1}^K t_k\varepsilon_k=\varepsilon, 
\end{align*}
which implies that~$\mathbb{Q}\in\mathcal{P}_\mathrm{mix}$.  As the choice of $\mathbb{Q} \in\mathcal{P}_\mathrm{mink}$ was arbitrary, we have $\mathcal{P}_\mathrm{mink} \subseteq\mathcal{P}_\mathrm{mix}$ and thus $v_\mathrm{mink}\leq v_\mathrm{mix}$. The inequalities $v_\mathrm{P}\leq v_\mathrm{stat}\leq v_\mathrm{mink}$ follow from Proposition~\ref{prop:mink}.
\end{proof}

\begin{proof}[Proof of Proposition~\ref{prop:elastic-net}]
We first derive~\eqref{eq:elastic-static-value}. To this end, fix vectors $\{\bm{a}_k\}_k$ with $\sum_{k=1}^K\bm{a}_k=\bm{0}$, define $g_k(\bm{\xi})\coloneqq t_k\ell(\bm{\xi})+\bm{a}_k^\top\bm{\xi}$ for each $k\in[K]$, and evaluate each worst-case expectation in~\eqref{prob:static_upper_bound} in closed form. Note that each $g_k$ is a real-valued convex function on $\mathbb{R}^m$ and is thus continuous. Since $\Xi=\mathbb{R}^m$ and $\mathbb E_{\hat{\mathbb P}}[\|\bm \xi\|_2]<\infty$, \citet[theorem~6.3]{mohajerin2018data} gives
\begin{equation}
   \label{eq:k-th-worst-case-expectation}
    \sup_{\mathbb{Q}_k\in\mathfrak{B}_k}\mathbb{E}_{\mathbb{Q}_k}[g_k(\bm{\xi})]=\mathbb{E}_{\hat{\mathbb{P}}}[g_k(\bm{\xi})]+\varepsilon_k\,\kappa_k (g_k),
\end{equation}
where $\kappa_k (g_k)$ is the Lipschitz modulus of~$g_k$ with respect to \hbox{$\|\cdot\|_{p_k}$}. By \citep[corollary~13.3.3]{R97:convex_analysis}, we have $\kappa_k (g_k)=\sup\{\|\bm{s}\|_{p_k^*}:g_k^*(\bm{s})<\infty\}$. Computing~$\kappa_k (g_k)$ thus requires identifying $\dom(g_k^*)$. To this end, assume temporarily that~$t_k\bm{\beta}\neq\bm{0}$. By \citep[theorem~12.3]{R97:convex_analysis}, we have $g_k^*(\bm{s})=(t_k\ell)^*(\bm{s}-\bm{a}_k)$, which implies that $\dom(g_k^*)=\{\bm{a}_k\}+\dom((t_k\ell)^*)$. If~$\bm{z}$ is not proportional to~$\bm{\beta}$, there exists~$\bm{w}\in\mathbb R^m$ with $\bm{\beta}^\top\bm{w}=0$ and $\bm{z}^\top\bm{w}>0$. In this case, setting $\bm{\xi}=s\bm{w}$ and letting~$s$ tend to infinity shows that $(t_k\ell)^*(\bm{z})=\infty$ and consequently that~$\bm{z}\notin\dom((t_k\ell)^*)$. If~$\bm z$ is proportional to~$\bm \beta$, on the other hand, then there exists~$\theta\in\mathbb R$ with~$\bm z=t_k\theta\bm{\beta}$ thanks to our assumption that~$t_k\bm \beta\neq\bm 0$. In this case, we have
\[
	(t_k\ell)^*(\bm z)=\sup_{\bm{\xi}} \left\{t_k\theta\bm{\beta}^\top\bm{\xi}-t_kL(\bm{\beta}^\top\bm{\xi}) \right\}=t_kL^*(\theta),
\]
which is finite if and only if $\theta\in\dom(L^*)$. In summary, the above reasoning implies that
\[
	\dom((t_k\ell)^*)=\{t_k\theta\bm{\beta}:\theta\in\dom(L^*)\} \quad \implies \quad \dom(g_k^*)=\{t_k\theta\bm{\beta}+\bm{a}_k:\theta\in\dom(L^*)\}.
\]
This in turn ensures that
\begin{align}
    \label{eq:kappa_k-formula}
    \kappa_k (g_k)=\sup_{\theta\in\dom(L^*)}\|t_k\theta\bm{\beta}+\bm{a}_k\|_{p_k^*}=\max\Big\{\|t_k\underline{\theta}\bm{\beta}+\bm{a}_k\|_{p_k^*},\|t_k\overline{\theta}\bm{\beta}+\bm{a}_k\|_{p_k^*}\Big\},
\end{align}
where the second equality holds because $\|t_k\theta\bm{\beta}+\bm{a}_k\|_{p_k^*}$ is convex in~$\theta$ such that its supremum over the compact interval $[\underline{\theta},\overline{\theta}]=\mathrm{cl}(\dom(L^*))$ is attained at one of the two endpoints. If~$t_k\bm{\beta}=\bm{0}$, then~$g_k$ is linear and $\dom(g_k^*)=\{\bm a_k\}$. In this case, we have~$\kappa_{k}(g_k)=\|\bm{a}_k\|_{p_k^*}$, and thus~\eqref{eq:kappa_k-formula} remains valid for all possible values of~$t_k$ and~$\bm\beta$. As $\mathbb E_{\hat{\mathbb P}}[\|\bm \xi\|_2]<\infty$ and each~$g_k$ is real-valued, exploiting the identities $\sum_{k=1}^K t_k=1$ and $\sum_{k=1}^K\bm{a}_k=\bm{0}$ to simplify $\sum_{k=1}^K\mathbb{E}_{\hat{\mathbb{P}}}[g_k(\bm{\xi})]=\mathbb{E}_{\hat{\mathbb{P}}}[\ell(\bm{\xi})]$ finally yields~\eqref{eq:elastic-static-value}. The inequality $v_\mathrm{P}\leq v_\mathrm{stat}$ follows from Proposition~\ref{prop:mink}, and $v_\mathrm{stat}\leq v_\mathrm{elastic}$ follows by combining $v_\mathrm{stat}\leq v_\mathrm{mink}$ from Proposition~\ref{prop:mink} with $v_\mathrm{mink}=v_\mathrm{elastic}$ established in Section~\ref{sec:elastic-net}. The recovery of Problem~\eqref{prob:elastic-net} at $\bm{a}_k=\bm{0}$ holds because $\max\{|\underline{\theta}|,|\overline{\theta}|\}=\mathrm{Lip}(L)$.
\end{proof}

\begin{proof}[Proof of Proposition~\ref{prop:sgd_wass}]
The proof proceeds in four steps.

\paragraph{Step 1: Well-Definedness.} Since $\ell(\bm{\beta};\bm{\xi})$ is $L$-Lipschitz in~$\bm{\beta}$ uniformly across all~$\bm{\xi}$ and upper semicontinuous in~$\bm{\xi}$ for every fixed~$\bm{\beta}$, it is jointly upper semicontinuous in~$\bm \beta$ and~$\bm \xi$. In addition, the support set~$\Xi$ is closed. Hence, for any fixed~$k\in[K]$, the function $-f_k(\bm{\theta};\bm{\xi},\hat{\bm{\xi}}_k) +\iota_\Xi(\bm{\xi})$, with~$\iota_\Xi$ denoting the indicator function of~$\Xi$, is jointly lower semicontinuous in~$\bm \theta$, $\bm\xi$ and $\hat{\bm \xi}_k$. By \citep[example~14.31]{rockafellar1998variational}, it thus constitutes a normal integrand. As~$\Xi$ is nonempty and compact, the optimal solution map $ \arg\max_{\bm \xi\in\Xi}f_k(\bm{\theta};\bm{\xi},\hat{\bm{\xi}}_k)$ thus admits a Borel measurable selector~$\bm{\xi}_k^\star$ \citep[corollary~14.6, theorem~14.37]{rockafellar1998variational}. Likewise,~$f_k(\bm{\theta};\bm{\xi},\hat{\bm{\xi}}_k)$ represents a normal integrand because it is continuous in~$\bm{\theta}$ and Borel measurable in $(\bm{\xi},\hat{\bm{\xi}}_k)$ \citep[example~14.29]{rockafellar1998variational}. Thus, a Borel measurable selector $\bm g_k(\bm{\theta};\bm{\xi},\hat{\bm{\xi}}_k)\in\partial_{\bm{\theta}}f_k(\bm{\theta};\bm{\xi},\hat{\bm{\xi}}_k)$ exists by \citep[corollary~14.6, theorem~14.56]{rockafellar1998variational}, and $\bm G$ in~\eqref{eq:stoch-subgrad} is a well-defined random vector. Finally, for any fixed $\bm{\beta}\in\mathbb R^n$, the loss function $\ell(\bm{\beta};\bm{\xi})$ is bounded above in~$\bm{\xi}$ by upper semicontinuity and compactness and bounded below by assumption. Consequently, one readily verifies that $\max_{\bm{\xi}\in\Xi}f_k(\bm{\theta};\bm{\xi},\hat{\bm{\xi}}_k)$ is bounded in $\hat{\bm{\xi}}_k$ on~$\Xi$ and that~$F(\bm{\theta})$ is finite.

\paragraph{Step 2: Weak Convexity.} For any fixed~$\bm\xi$ and~$\hat{\bm \xi}_k$, the function $t_k\ell(\bm{\beta};\bm{\xi})$ is $t_k\rho$-weakly convex in~$\bm\beta$, while all other terms of~$f_k(\bm{\theta};\bm{\xi},\hat{\bm{\xi}}_k)$ are affine in~$\bm{\theta}$. Hence, $f_k(\bm{\theta};\bm{\xi},\hat{\bm{\xi}}_k)$ is $t_k\rho$-weakly convex in~$\bm{\theta}$. Pointwise maxima and expectations of $t_k\rho$-weakly convex functions are $t_k\rho$-weakly convex because adding $\frac{t_k\rho}{2}\|\bm \beta\|_2^2$ commutes with these operations, both of which preserve convexity. Summing over~$k$ and using $\sum_{k=1}^K t_k=1$ shows via~\eqref{prob:wass-dual} that $F(\bm \theta)$ is $\rho$-weakly convex.

\paragraph{Step 3: Second-Moment Bound.} The components of the subgradient $\bm g_k(\bm{\theta};\bm{\xi},\hat{\bm{\xi}}_k)$ corresponding to the parameters $\bm{a}_k$, $\bm{A}_k$ and $\gamma_k$ are $\bm{\xi}$, $\bm{\xi}\bm{\xi}^\top$ and $\varepsilon_k^p-\|\bm{\xi}-\hat{\bm{\xi}}_k\|^p$, respectively. They are bounded uniformly on the compact set $\Xi\times\Xi$. The component corresponding to $\bm{\beta}$ has Euclidean norm at most $t_kL$ because every subgradient of an $L$-Lipschitz function has norm at most~$L$. As~$\bm{\xi}_k^\star (\bm{\theta},\hat{\bm{\xi}}_k) \in\Xi$, the norm $\| \bm G(\bm{\theta};\hat{\bm{\xi}})\|_2$ is thus bounded by a constant independent of~$\bm{\theta}$ and~$\hat{\bm{\xi}}$. Consequently, the second moment of $\bm G(\bm{\theta};\hat{\bm{\xi}})$ is bounded uniformly across all~$\bm \theta\in\Theta$.

\paragraph{Step 4: Unbiasedness.} Fix $k\in[K]$ and $\hat{\bm{\xi}}_k\in\Xi$, and write $f^\star_k(\bm{\theta};\hat{\bm{\xi}}_k)\coloneqq\max_{\bm{\xi}\in\Xi}f_k(\bm{\theta};\bm{\xi},\hat{\bm{\xi}}_k)$. As the effective loss $f_k(\bm{\theta};\bm{\xi},\hat{\bm{\xi}}_k)$ is $t_k\rho$-weakly convex in~$\bm{\theta}$, every $\bm{g}\in\partial_{\bm{\theta}}f_k(\bm{\theta};\bm{\xi},\hat{\bm{\xi}}_k)$ satisfies
\[
	f_k(\bm{\theta}';\bm{\xi},\hat{\bm{\xi}}_k)\geq f_k(\bm{\theta};\bm{\xi},\hat{\bm{\xi}}_k)+\langle\bm{g},\bm{\theta}'-\bm{\theta}\rangle-\frac{t_k\rho}{2}\|\bm{\theta}'-\bm{\theta}\|_2^2 \quad \forall \bm{\theta}' \in \Theta_0.
\]
 Evaluating this inequality at $\bm g = \bm g_k(\bm{\theta};\bm{\xi},\hat{\bm{\xi}}_k)$ and $\bm{\xi}=\bm{\xi}_k^\star(\bm{\theta},\hat{\bm{\xi}}_k)$, and using the basic relations
 \[
 	f^\star_k(\bm{\theta}';\hat{\bm{\xi}}_k)\geq f_k(\bm{\theta}';\bm{\xi}_k^\star(\bm{\theta},\hat{\bm{\xi}}_k),\hat{\bm{\xi}}_k) \quad \text{and} \quad f^\star_k(\bm{\theta};\hat{\bm{\xi}}_k)=f_k(\bm{\theta};\bm{\xi}_k^\star (\bm{\theta},\hat{\bm{\xi}}_k) ,\hat{\bm{\xi}}_k)
\]
gives
 \[
 	f^\star_k(\bm{\theta}';\hat{\bm{\xi}}_k)\geq f^\star_k(\bm{\theta};\hat{\bm{\xi}}_k)+\langle \bm g_k(\bm{\theta};\bm{\xi}_k^\star (\bm{\theta},\hat{\bm{\xi}}_k) ,\hat{\bm{\xi}}_k),\bm{\theta}'-\bm{\theta}\rangle-\frac{t_k\rho}{2}\|\bm{\theta}'-\bm{\theta}\|_2^2 \quad \forall \bm{\theta}' \in \Theta_0.
\]
Taking expectations over $\hat{\bm{\xi}}_k\sim\hat{\mathbb{P}}_k$ and summing over $k\in[K]$ yields 
\[
	F(\bm{\theta}')\geq F(\bm{\theta})+\langle\mathbb{E}_{\hat{\mathbb P}} [\bm G(\bm{\theta};\hat{\bm{\xi}})],\bm{\theta}'-\bm{\theta}\rangle-\frac{\rho}{2}\|\bm{\theta}'-\bm{\theta}\|_2^2 \quad \forall \bm{\theta}' \in \Theta_0,
\]
where $F(\bm{\theta})$ and $F(\bm{\theta}')$ are finite by the reasoning in Step~1, and the expectation of $\bm G(\bm{\theta};\hat{\bm{\xi}})$ is well-defined by the reasoning in Step~3. We thus conclude that~$\mathbb{E}[\bm G(\bm{\theta};\hat{\bm{\xi}})]\in\partial F(\bm{\theta})$.
\end{proof}

\section{Direct Restrictions of the Dual Functions}\label{app:dual_lambda_approx}
By Proposition~\ref{prop:affine_upper_bound}, restricting the decision rules $\bm{\alpha}_k(\bm{\xi}',\bm{\xi})$ of Section~\ref{sec:main} to the affine (or static) class amounts to restricting the dual functions~$\lambda_k(\bm{\xi})$ in~\eqref{eq:reformulation_dual} to the sum of the weighted loss~$t_k\ell(\bm{\xi})$ and a quadratic (or linear) function of~$\bm{\xi}$. These dual functions thus depend on the loss. A natural alternative is to restrict the dual functions directly to a loss-independent class, such as the quadratic or linear polynomials. We show that this is never better and in general strictly worse. The reason is elementary. Feasibility in~\eqref{eq:reformulation_dual} only requires $\sum_{k=1}^K\lambda_k(\bm \xi) \geq\ell(\bm \xi)$ on~$\Xi$, and the objective is monotonically nondecreasing in each~$\lambda_k$. Thus, any slack $\sum_{k=1}^K\lambda_k(\bm \xi)-\ell(\bm \xi)$ can be subtracted without harm. Subtracting the share $t_k$ of the slack from~$\lambda_k(\bm \xi)$ produces exactly the loss-adjusted functions of Section~\ref{sec:main}. To formalize this idea, let $\mathcal{L}\subseteq C(\Xi)$ be any linear space of continuous functions. Denote by $v_{\mathcal{L}}^\lambda$ the optimal value of~\eqref{eq:reformulation_dual} with the restriction $\lambda_k\in\mathcal{L}$, and by $v_{\mathcal{L}}^\ell$ its optimal value with the restriction $\lambda_k=t_k\ell+\psi_k$, where $\psi_k\in\mathcal{L}$ and $\sum_{k=1}^K\psi_k=0$ on~$\Xi$. When $\mathcal{L}$ is the space of all affine or quadratic polynomials, $v_{\mathcal{L}}^\ell$ coincides with $v_\mathrm{stat}$ or $v_\mathrm{aff}$, respectively, because the constant terms in the $\psi_k$ sum to zero and drop out of the objective.

\begin{proposition}\label{prop:lambda_general}
Suppose that Assumption~\ref{assump:general}\,\ref{assump:general-support} holds and that $\ell$ is continuous. Let $\mathcal{L}\subseteq C(\Xi)$ be a linear space, and fix $\bm{t}\in\Delta(K)$. Then, $v_\mathrm{P}\leq v_{\mathcal{L}}^\ell\leq v_{\mathcal{L}}^\lambda$. If additionally $\ell\in\mathcal{L}$, then $v_{\mathcal{L}}^\ell=v_{\mathcal{L}}^\lambda$.
\end{proposition}

\begin{proof}
We show that any solution feasible in~\eqref{eq:reformulation_dual} under the restriction $\lambda_k\in\mathcal{L}$ can be converted into a solution that is feasible under the restriction $\lambda_k=t_k\ell+\psi_k$ with $\psi_k\in\mathcal{L}$ and $\sum_{k=1}^K\psi_k=0$, has no larger objective value, and is in particular feasible in~\eqref{eq:reformulation_dual} without any restriction. The three inequalities $v_\mathrm{P}\leq v_\mathrm{D}\leq v_{\mathcal{L}}^\ell\leq v_{\mathcal{L}}^\lambda$ then follow, the first being weak duality (Proposition~\ref{prop:reformulation_main}).

Let $\{\lambda_k\}_k\subseteq\mathcal{L}$ be feasible in~\eqref{eq:reformulation_dual}, that is, $\sum_{k=1}^K\lambda_k(\bm{\xi})\geq\ell(\bm{\xi})$ for all $\bm{\xi}\in\Xi$, and define $\psi_k\coloneqq\lambda_k-t_k\sum_{j=1}^K\lambda_j$ and $\lambda'_k\coloneqq t_k\ell+\psi_k$ for $k\in[K]$. Then $\psi_k\in\mathcal{L}$ because $\mathcal{L}$ is a linear space, $\sum_{k=1}^K\psi_k=0$ because $\sum_{k=1}^Kt_k=1$, and $\lambda'_k\in C(\Xi)$ because $\ell$ is continuous. The functions $\lambda'_k$ satisfy $\sum_{k=1}^K \lambda'_k=\ell$, so they are feasible in~\eqref{eq:reformulation_dual}, and they satisfy $\lambda'_k(\bm{\xi})=\lambda_k(\bm{\xi})-t_k\big(\sum_{j=1}^K\lambda_j(\bm{\xi})-\ell(\bm{\xi})\big)\leq\lambda_k(\bm{\xi})$ for all $\bm{\xi}\in\Xi$ by feasibility of $\{\lambda_k\}_k$. Since all functions involved are continuous on the compact set $\Xi$ and hence bounded, the expectations below are finite, and monotonicity of the expectation gives $\sup_{\mathbb{Q}_k\in\mathfrak{B}_k}\mathbb{E}_{\mathbb{Q}_k}[\lambda'_k(\bm{\xi})]\leq\sup_{\mathbb{Q}_k\in\mathfrak{B}_k}\mathbb{E}_{\mathbb{Q}_k}[\lambda_k(\bm{\xi})]$ for every $k\in[K]$. Summing over $k$ shows that the objective value of $\{\lambda'_k\}_k$ does not exceed that of $\{\lambda_k\}_k$, which establishes the claim. Finally, if $\ell\in\mathcal{L}$, then $t_k\ell+\psi_k\in\mathcal{L}$ for every admissible $\psi_k$, so the loss-adjusted solutions are themselves admissible under the restriction $\lambda_k\in\mathcal{L}$, and thus we find $v_{\mathcal{L}}^\lambda\leq v_{\mathcal{L}}^\ell$.
\end{proof}

We now set $\mathcal{L}$ equal to the space of quadratic polynomials on $\Xi$. With this choice, $v_{\mathcal{L}}^\ell$ is the optimal value of~\eqref{eq:reformulation_dual} obtained by restricting each dual function to the weighted loss plus a quadratic polynomial, where the quadratic, linear and constant coefficients each sum to zero across $k\in[K]$. As the constants cancel in the objective, this value matches $v_\mathrm{aff}$, the affine decision rule approximation of Proposition~\ref{prop:affine_upper_bound}. The value $v_{\mathcal{L}}^\lambda$, in turn, is what one obtains by restricting the dual functions directly to quadratic polynomials. Writing $v^\lambda_\mathrm{quad}\coloneqq v_{\mathcal{L}}^\lambda$, Proposition~\ref{prop:lambda_general} thus implies that~$v_\mathrm{aff}\leq v_\mathrm{quad}^\lambda$. Put differently, affine decision rules for $\bm{\alpha}_k$ are never worse than quadratic decision rules for $\lambda_k$, even though the latter have more free parameters. The following corollary makes the direct quadratic restriction explicit by eliminating its constant terms.

\begin{corollary}\label{cor:lambda_quad}
Suppose that Assumption~\ref{assump:general} holds and that $\ell$ is continuous. Restricting the dual functions in~\eqref{eq:reformulation_dual} to quadratic polynomials $\lambda_k(\bm{\xi})=\langle\bm{B}_k,\bm{\xi}\bm{\xi}^\top\rangle+\bm{b}_k^\top\bm{\xi}+c_k$ yields
\begin{align}\tag{D$_\mathrm{quad}^\lambda$}\label{prob:lambda_quad}
    v_\mathrm{quad}^\lambda\coloneqq \left\{ \begin{array}{cl} \displaystyle \inf_{\{\bm{B}_k\in\mathbb{S}^m\}_k} \inf_{\{\bm{b}_k\in\mathbb{R}^m\}_k} & \displaystyle \sum_{k=1}^K \sup_{\mathbb{Q}_k\in\mathfrak{B}_k}\mathbb{E}_{\mathbb{Q}_k} \big[\langle\bm{B}_k,\bm{\xi}\bm{\xi}^\top\rangle+\bm{b}_k^\top\bm{\xi}\big] \\
    & \quad \displaystyle +\sup_{\bm{\xi}\in\Xi}\Big\{\ell(\bm{\xi})-\sum_{j=1}^K\big(\langle\bm{B}_j,\bm{\xi}\bm{\xi}^\top\rangle+\bm{b}_j^\top\bm{\xi}\big)\Big\},
    \end{array} \right.
\end{align}
and $v_\mathrm{P}\leq v_\mathrm{aff}\leq v_\mathrm{quad}^\lambda$ for every $\bm{t}\in\Delta(K)$, with equality $v_\mathrm{aff}=v_\mathrm{quad}^\lambda$ whenever $\ell$ is quadratic.
\end{corollary}

\begin{proof}
Fix any $\bm{B}_k\in\mathbb S^m$ and $\bm{b}_k\in\mathbb R^m$, $k\in[K]$. Each $\mathfrak{B}_k$ is nonempty by Assumption~\ref{assump:general}\,\ref{assump:general-nonempty-intersection}, and the quadratic polynomials are bounded on the compact set $\Xi$. Thus, the objective of~\eqref{eq:reformulation_dual}~equals
\[
	\sum_{k=1}^K\sup_{\mathbb{Q}_k\in\mathfrak{B}_k}\mathbb{E}_{\mathbb{Q}_k}[\langle\bm{B}_k,\bm{\xi}\bm{\xi}^\top\rangle+\bm{b}_k^\top\bm{\xi}]+\sum_{k=1}^Kc_k
\]
and is increasing in $\sum_{k=1}^Kc_k$, while the semi-infinite constraint can be represented as 
\[
	\sup_{\bm{\xi}\in\Xi} \bigg\{ \ell(\bm{\xi})-\sum_{j=1}^K(\langle\bm{B}_j,\bm{\xi}\bm{\xi}^\top\rangle+\bm{b}_j^\top\bm{\xi}) \bigg\} \leq\sum_{k=1}^Kc_k,
\]
where the supremum is finite because $\ell$ is continuous and $\Xi$ is compact. Minimizing over the constants therefore makes the constraint tight, and substituting it into the objective gives~\eqref{prob:lambda_quad}. The inequalities as well as the equality for quadratic loss functions follow from Proposition~\ref{prop:lambda_general} with $\mathcal{L}$ the space of quadratic polynomials, for which $v_{\mathcal{L}}^\ell=v_\mathrm{aff}$ as noted above.
\end{proof}

Setting $\bm{B}_k=\bm{0}$ for all $k\in[K]$ in~\eqref{prob:lambda_quad} yields the direct restriction of the dual functions to affine polynomials, with optimal value $v_\mathrm{aff}^\lambda$. Setting $\mathcal{L}$ to the space of affine polynomials gives $v_\mathrm{P}\leq v_\mathrm{stat}\leq v_\mathrm{aff}^\lambda$, with equality whenever $\ell$ is affine. Static decision rules for $\bm{\alpha}_k$ thus dominate affine decision rules for $\lambda_k$. The following example shows that the hierarchy $v_\mathrm{aff}\leq v_\mathrm{quad}^\lambda\leq v_\mathrm{aff}^\lambda$ can be strict at every step, and that it is so already for $K=1$, where $v_\mathrm{aff}=v_\mathrm{P}$.

\begin{example}\label{ex:lambda_gap}
Assume that $K=1$, $\Xi=[-1,1]$, and $\mathfrak{B}_1=\{\hat{\mathbb{P}}\}$, where $\hat{\mathbb{P}}$ is the uniform distribution on $\Xi$. Set $\ell(\xi)=|\xi|$. Since $K=1$ forces $a_1=A_1=0$ and $t_1=1$, both decision rule approximations are exact, that is, $v_\mathrm{P}=v_\mathrm{stat}=v_\mathrm{aff}=\mathbb{E}_{\hat{\mathbb{P}}} [|\xi|]=1/2$. As $\mathbb{E}_{\hat{\mathbb{P}}}[\xi]=0$ and $\Xi$ is symmetric around~$0$, it is optimal to set $b=0$ in~\eqref{prob:lambda_quad}. As $\mathbb{E}_{\hat{\mathbb{P}}}[\xi^2]=1/3$, we further have
\[
	v_\mathrm{aff}^\lambda=\inf_{b\in\mathbb R} \max\{1-b,1+b\}=1\quad \text{and} \quad v_\mathrm{quad}^\lambda=\inf_{B\in\mathbb R} \left\{\frac{B}{3}+\sup_{|\xi|\leq 1} \left(|\xi|-B\xi^2 \right) \right\} =\frac{1}{\sqrt{3}}.
\]
The last equality holds because the supremum equals $1-B$ for $B\leq 1/2$ and $1/(4B)$ for $B\geq 1/2$, and the infimum is attained at $B=\sqrt{3}/2$. Thus $v_\mathrm{aff}^\lambda=1>v_\mathrm{quad}^\lambda=1/\sqrt{3}>1/2=v_\mathrm{aff}=v_\mathrm{stat}$.
\end{example}

Problem~\eqref{prob:lambda_quad} is convex because its objective function is given by a pointwise supremum of functions affine in $\{\bm{B}_k,\bm{b}_k\}_k$. The $K$ embedded worst-case expectations and the robust term involving the supremum over $\Xi$ are tractable, for instance, if each $\mathfrak B_k$ is a $2$-Wasserstein ball for every $k\in[K]$, $\ell$ is a quadratic loss function and $\Xi=\mathbb{R}^m$, in which case they admit semidefinite reformulations by the arguments familiar from the proof of Lemma~\ref{lem:affine-sdp-equivalence}. In that case $v_\mathrm{quad}^\lambda=v_\mathrm{aff}$ by arguments familiar from the proof of Corollary~\ref{cor:lambda_quad}, and the exactness result of Proposition~\ref{prop:affine-tightness} carries over to the direct quadratic restriction of the dual functions~$\lambda_k$, $k\in[K]$.

\begin{remark}[Connection to Wasserstein Barycenters]\label{rem:barycenter}
If $\ell\equiv 0$, then our decision rule approximations are linked to Wasserstein barycenter computations in disguise. The static approximation of Corollary~\ref{cor:static_upper_bound}  recovers the barycenter's mean and the affine approximation of Proposition~\ref{prop:affine_upper_bound} its mean and covariance matrix, which explains both their exactness for Gaussian references and their blindness to higher-order structure. To see this, set $\ell\equiv 0$, and let each $\mathfrak{B}_k$ be a $p$-Wasserstein ball of radius $\varepsilon_k$. Then Problem~\eqref{prob:sup-primal} merely asks whether the $K$ balls intersect, while~$v_\mathrm{stat}$ and~$v_\mathrm{aff}$ are certificates for this feasibility question. Applying Wasserstein duality \citep[theorem~4.18]{kuhn2025distributionally} to each worst-case expectation in~\eqref{prob:static_upper_bound} introduces one nonnegative multiplier per ball, and these multipliers correspond, up to normalization, to the barycenter weights~$\bm{w}\in\Delta(K)$ in~\eqref{eq:barycenter}. After dualization, $v_\mathrm{stat}$ equals $-\infty$ if $\sum_{k=1}^K w_k\varepsilon_k^p$ is smaller than a lower bound on the value of the barycenter problem~\eqref{eq:barycenter}, and it equals~$0$ otherwise, where the lower bound is obtained by applying to~\eqref{eq:barycenter} the copy-and-dualize construction of Proposition~\ref{prop:reformulation_main} and Lemma~\ref{lem:adjustable_opt_trick} and restricting the resulting dual functions to linear functions, exactly as the static rules do in Proposition~\ref{prop:affine_upper_bound}. The same holds for~$v_\mathrm{aff}$ with quadratic dual functions. For $p=2$, $\Xi=\mathbb{R}^m$ and reference distributions with positive definite covariance matrices, the arguments of Lemma~\ref{lem:affine-sdp-equivalence} show that the static lower bound is the weighted least-squares problem for the reference means, whose solution is the weighted mean, and that the affine lower bound is the barycenter problem for the Gelbrich distance~\eqref{eq:gelbrich-dist} in the space of mean-covariance pairs, whose solution is~\eqref{eq:barycenter-mean-cov}. As the Gelbrich distance is a lower bound on the $2$-Wasserstein distance that is attained on positive-semidefinite affine families (see Definition~\ref{def:psd-affine-family} and \citet[theorem~2]{nguyen2021mean}), the affine bound is exact for Gaussian references, in line with Proposition~\ref{prop:affine-tightness} and the closed-form Gaussian barycenter of Section~\ref{sec:barycenter}, and it discards everything beyond second moments otherwise. The exact counterpart of this correspondence, namely that the $K$ balls intersect if and only if a barycenter problem with free weights has value at most $\sum_{k=1}^K w_k\varepsilon_k^p$, is also alluded to in~\citep[p.~18]{rychener2024wasserstein}. The decision rule approximations of Section~\ref{sec:main} are exactly the moment-level relaxations of that barycenter problem.
\end{remark}

\section{\texorpdfstring{Details of Section~\ref{sec:SGD-experiment}}{Details of Section 5.3}}\label{app:sgd-experiment}
This appendix collects supplementary information on the experiment of Section~\ref{sec:SGD-experiment} omitted from the main text. The setting of that section, in which a decision maker receives noisy moments of the numerical features and a noisy histogram of the categorical features, arises naturally in practice. Surveys typically summarize numerical features by their means and covariances, because binning them into cells would discard information and inflate the number of cells, whereas they report categorical features as histograms of counts. Similarly, data holders that protect privacy commonly release several differentially private statistics about the same records. For example, the 2020 U.S.\ Census publications include a Supplemental Demographic and Housing Characteristics File with differentially private moment statistics such as average household and family sizes~\citep{sexton2024phsafe}, as well as a Noisy Measurement File with differentially private histograms~\citep{abowd2022topdown}. When both statistics describe the same records and each elementary ambiguity set contains their distribution with probability at least $95\%$, the union bound ensures that the fusion ambiguity set $\mathfrak{B}_1\cap\mathfrak{B}_2$ contains it with probability at least $90\%$.

\paragraph{Noisy moments.} Let $\bm{\varphi}(\bm{\xi})\in\mathbb{R}^{134}$ stack the distinct entries of $(\bm{x}_\mathrm{num},y)$ and of its outer product with itself, so that $\bm{m}_{\mathbb{Q}}\coloneqq\mathbb{E}_{\mathbb{Q}}[\bm{\varphi}(\bm{\xi})]$ is the vector of distinct first and second moments of~$(\bm{x}_\mathrm{num},y)$ under any distribution~$\mathbb{Q}$ on~$\Xi$. The radius~$\varepsilon_1$ of the moment set~$\mathfrak{B}_1$ is calibrated by the following principle. If the Euclidean ball of radius~$\varepsilon_1$ around the noisy estimate~$\hat{\bm{m}}$ contains the moments $\bm{m}_{\mathbb{Q}}$ of a distribution~$\mathbb{Q}$ with probability $95\%$, then $\mathfrak{B}_1$ contains~$\mathbb{Q}$ with the same probability. We describe two settings in which such radii are available. First, suppose that a survey reports the empirical moments~$\hat{\bm{m}}$ of $N_1$ independent samples from a population distribution~$\mathbb{P}$. By the central limit theorem, the error $\hat{\bm{m}}-\bm{m}_{\mathbb{P}}$ is asymptotically Gaussian with a covariance matrix that can be estimated from the survey. Its confidence regions are therefore ellipsoids, whose size asymptotically follows from a quantile of the chi-squared distribution with as many degrees of freedom as there are moments~\citep{hotelling1931generalization,delage2010distributionally}. Setting~$\varepsilon_1$ to the largest semi-axis of such an ellipsoid yields a Euclidean ball that contains the ellipsoid, so that $\mathfrak{B}_1$ contains~$\mathbb{P}$ with at least the prescribed asymptotic confidence. Second, suppose that a data holder possesses the $N$ records of the hidden empirical distribution~$\mathbb{P}_N$ and releases~$\hat{\bm{m}}$ under $(\varepsilon^\mathrm{DP},\delta^\mathrm{DP})$-differential privacy~\citep{dwork2014algorithmic}. This means that for any two datasets that differ in a single record, the probabilities of any set of outputs differ by at most a factor~$\exp(\varepsilon^\mathrm{DP})$ up to an additive slack~$\delta^\mathrm{DP}$, so that the release reveals little about any individual. A standard way to achieve this is the Gaussian mechanism, which releases $\hat{\bm{m}}=\bm{m}_{\mathbb{P}_N}+\sigma\bm{z}$ with $\bm{z}\sim\mathcal{N}(\bm{0},\bm{I})$~\citep{balle2018improving,biswas2020coinpress}. The standard deviation~$\sigma$ is proportional to the largest change in $\bm{m}_{\mathbb{P}_N}$ that a single record can cause. As all entries of $\bm{\varphi}(\bm{\xi})$ lie in $[-1,1]$, replacing one record changes each moment by at most $2/N$, so that the noise, and with it the radius below, decreases like $1/N$. Since $\|\hat{\bm{m}}-\bm{m}_{\mathbb{P}_N}\|_2/\sigma$ follows a chi distribution with $134$ degrees of freedom, setting $\varepsilon_1$ to $\sigma$ times its $95\%$ quantile ensures that $\mathfrak{B}_1$ contains~$\mathbb{P}_N$ with probability exactly $95\%$. Our experiments use this setting with $(\varepsilon^\mathrm{DP},\delta^\mathrm{DP})=(0.5,10^{-6})$.

\paragraph{Noisy histograms.} The cells that define the histogram set~$\mathfrak{B}_2$ group the values of~$(\bm{x}_\mathrm{cat},y)$ by the label, the sex, the education level, the marital status, the latest repayment category (current, one month late, at least two months late) and the number of late months among the five preceding months (zero, one or two, three to five), which yields $J=432$ cells. Let $\bm{\psi}(\bm{\xi})\in\{0,1\}^J$ be the indicator vector of the cell that contains~$(\bm{x}_\mathrm{cat},y)$, so that $\bm{p}_{\mathbb{Q}}\coloneqq\mathbb{E}_{\mathbb{Q}}[\bm{\psi}(\bm{\xi})]$ represents the histogram of any distribution~$\mathbb{Q}$ on~$\Xi$. The total variation distance $\frac12\|\hat{\bm{p}}-\bm{p}_{\mathbb{Q}}\|_1$ between a noisy histogram~$\hat{\bm{p}}$ and~$\bm{p}_{\mathbb{Q}}$ equals the largest error that $\hat{\bm{p}}$ makes in the probability of any group of cells. The radius~$\varepsilon_2$ of the histogram set~$\mathfrak{B}_2$ is calibrated by the same principle as~$\varepsilon_1$. If the total variation ball of radius~$\varepsilon_2$ around~$\hat{\bm{p}}$ contains~$\bm{p}_{\mathbb{Q}}$ with probability $95\%$, then $\mathfrak{B}_2$ contains~$\mathbb{Q}$ with the same probability. We consider the same two settings as for the moments. First, suppose that a survey reports the empirical histogram~$\hat{\bm{p}}$ of $N_2$ independent samples from a population distribution~$\mathbb{P}$. The cell counts $N_2\hat{\bm{p}}$ then follow a multinomial distribution,~and
\[
	\mathbb{P}^{N_2} \big(\|\hat{\bm{p}}-\bm{p}_{\mathbb{P}}\|_1\geq r \big) \leq(2^J-2)\exp(-N_2r^2/2) \quad \forall r>0;
\]
see \citet[theorem~2.1]{weissman2003inequalities}. Setting the bound to $5\%$ and halving the $\ell_1$-radius yields $\varepsilon_2= (\log((2^J-2)/0.05)/(2N_2))^{1/2}$, which ensures $\mathbb{P}\in\mathfrak{B}_2$ with probability at least~$95\%$. Second, suppose that the data holder releases a histogram of its $N$ records under $(\varepsilon^\mathrm{DP},0)$-differential privacy. The Laplace mechanism achieves this by adding independent Laplace noise with scale $2/\varepsilon^\mathrm{DP}$ to each cell count, as replacing one record changes two counts by one each~\citep{hay2010boosting}. To keep a histogram after adding noise, we post-process the noisy counts in two steps. As $\mathfrak{B}_2$ restricts the probabilities of the cells rather than their counts, we first divide the noisy counts by~$N$, which uses no private information because neighboring datasets have the same size. Replacing the resulting noisy proportions by the closest probability vector in Euclidean distance then yields the histogram~$\hat{\bm{p}}$. To restore a total probability of one at the smallest Euclidean cost, this projection onto the probability simplex subtracts the same constant from all cells and sets the cells that become negative to zero, and as post-processing it preserves differential privacy~\citep[proposition~2.1]{dwork2014algorithmic}. The resulting histogram satisfies $\frac12\|\hat{\bm{p}}-\bm{p}_{\mathbb{P}_N}\|_1\leq\frac12\|\bm{\eta}-\bar\eta\bm{1}\|_1$, where~$\bm{\eta}$ is the noise added to the cell proportions and $\bar\eta$ its average. Indeed, $\bm{w}\coloneqq\bm{p}_{\mathbb{P}_N}+\bm{\eta}-\bar\eta\bm{1}$ is the projection of the noisy proportions onto the hyperplane of vectors that sum to one, and $\hat{\bm{p}}$ is the projection of~$\bm{w}$ onto the probability simplex, that is, $\hat p_c=\max\{w_c-\tau,0\}$ for some $\tau\geq0$. Hence, $\max\{\hat p_c-p_c,0\}\leq\max\{w_c-p_c,0\}$ for every cell~$c$ and every probability vector~$\bm{p}$. As $\hat{\bm{p}}$, $\bm{w}$ and~$\bm{p}$ all sum to one, their $\ell_1$-distances equal twice the sums of these positive parts, and the claim follows. As this bound does not depend on the data, we set~$\varepsilon_2$ to its $95\%$ quantile, which we estimate by simulating the Laplace noise alone, so that $\mathfrak{B}_2$ contains~$\mathbb{P}_N$ with probability at least $95\%$. As each of the $J$ counts receives noise of constant size, $\varepsilon_2$ grows linearly in~$J$ and decreases like $1/(N\varepsilon^\mathrm{DP})$. This is why the cells group the values of~$(\bm{x}_\mathrm{cat},y)$, whose histogram would otherwise have more than two million entries, and it explains why \textsc{tv} is weak for small~$N$ but improves quickly as $N$ grows. Our experiments use this setting with $(\varepsilon^\mathrm{DP},\delta^\mathrm{DP})=(0.5,0)$, so that the moments and the histogram are jointly $(1,10^{-6})$-differentially private by basic composition~\citep[theorem~3.16]{dwork2014algorithmic}.

\paragraph{Subgradient method.} Each update of the subgradient method of Section~\ref{sec:sgd} requires, for each elementary ambiguity set~$\mathfrak{B}_k$ involved, the worst-case expectation of an integrand~$\lambda_k$ together with a worst-case distribution, from which Danskin's theorem yields a subgradient. For \textsc{moment} and \textsc{tv}, the integrand is the loss $\ell(\bm{\beta};\bm{\xi})$, and for \textsc{minkowski} it is $t_k\ell(\bm{\beta};\bm{\xi})$. For \textsc{affine}, the integrand is $\lambda_k(\bm{\xi})=t_k\ell(\bm{\beta};\bm{\xi})+\bm{a}_k^\top\bm{s}(\bm{\xi})+\langle\bm{A}_k,\bm{\chi}(\bm{\xi})\bm{\chi}(\bm{\xi})^\top\rangle$, which uses the lifting $\bm{s}(\bm{\xi})\coloneqq(\bm{\chi}(\bm{\xi}),\bm{\psi}(\bm{\xi}))$ in the sense of Remark~\ref{rem:lifted-rules}. Here, $\bm{\chi}(\bm{\xi})$ is obtained from $\bm{\xi}$ by replacing the categorical features with their one-hot encoding, and $\bm{\psi}(\bm{\xi})$ is the vector of cell indicators defined above. The quadratic terms thus involve only~$\bm{\chi}(\bm{\xi})$, while the linear terms include one offset per cell, which lets the decision rules adapt to the cells constrained by~$\mathfrak{B}_2$. \textsc{static} uses the same integrand with~$\bm{A}_k=\bm{0}$.

For the moment set~$\mathfrak{B}_1$, we replace the worst-case expectation of~$\lambda_1$ by its dual~\citep[theorem~4.5]{kuhn2025distributionally}, which minimizes $\bm{b}^\top\hat{\bm{m}}+\varepsilon_1\|\bm{b}\|_2+\max_{\bm{\xi}\in\Xi}\{\lambda_1(\bm{\xi})-\bm{b}^\top\bm{\varphi}(\bm{\xi})\}$ over the multipliers~$\bm{b}\in\mathbb{R}^{134}$ and equals the worst-case expectation by strong duality. We update~$\bm{b}$ together with the network weights and the decision rule coefficients. Given a maximizer~$\bm{\xi}^\star$, the subgradient of the dual objective is the gradient of~$\lambda_1$ at~$\bm{\xi}^\star$ with respect to the network weights and coefficients, and $\hat{\bm{m}}-\bm{\varphi}(\bm{\xi}^\star)+\varepsilon_1\bm{b}/\|\bm{b}\|_2$ with respect to~$\bm{b}$. For the histogram set~$\mathfrak{B}_2$, the worst case is available in closed form. Let $\Xi_c\subseteq \Xi$ denote the observations whose categorical features and label lie in cell~$c$, and let $M_c\coloneqq\max_{\bm{\xi}\in\Xi_c}\lambda_2(\bm{\xi})$ be attained at~$\bm{\xi}_c^\star$. As~$\mathfrak{B}_2$ restricts only the cell probabilities, the worst-case distribution places the probability of each cell~$c$ on~$\bm{\xi}_c^\star$. Its histogram~$\bm{p}^\star$ arises from~$\hat{\bm{p}}$ by moving probability~$\varepsilon_2$ from the cells with the smallest values~$M_c$, each giving up at most its probability under~$\hat{\bm{p}}$, to the cell with the largest value~\citep[proposition~6.13]{kuhn2025distributionally}. The worst-case expectation is $\sum_c p^\star_cM_c$, and its gradient is $\sum_c p^\star_c\nabla\lambda_2(\bm{\xi}_c^\star)$.

All maximization problems over~$\Xi$ and over the cells~$\Xi_c$ are nonconvex and mixed-integer, because the loss is not concave in~$\bm{\xi}$ and the categorical features are discrete. We approximate them by multi-start projected gradient ascent, and thus the resulting subgradients are approximate. Each ascent relaxes the one-hot encoding of every categorical feature to its probability simplex and projects the numerical features onto~$[-1,1]$. After $32$ steps, it rounds each relaxed encoding to the nearest valid one-hot encoding. For each cell~$c$, it also ensures that the rounded point lies in~$\Xi_c$. We use eight starting points. The first is the maximizer of the previous update, which remains a good candidate because consecutive updates change the network only slightly. Three more perturb this maximizer by Gaussian noise with standard deviation~$0.15$ on the numerical features~\citep{jin2017escape} and Gumbel noise with scale~$0.35$ on the one-hot encodings~\citep{jang2017categorical}. The last four are random points with uniformly distributed numerical features and random categories~\citep{madry2018towards}. Among the eight rounded end points of the ascents, we keep the one with the largest objective value. All robust methods use the constant step size~$0.05$ for the network weights, the multipliers~$\bm{b}$ and the coefficients of the cell indicators, and~$0.0125$ for the remaining decision rule coefficients. At each step, we enforce the zero-sum constraints by subtracting from every decision rule coefficient its average over $k\in\{1,2\}$. \textsc{oracle} uses the step size $0.05$ and minibatches of $128$ records drawn from the hidden empirical distribution~$\mathbb{P}_N$.

\end{document}